\documentclass[aos,MSNbibl,seceqn,dvips]{arximspdf}
\usepackage{dcolumn}
\usepackage{graphicx}

\doi{10.1214/14-AOS1222} 
\volume{42}
\issue{4}
\pubyear{2014}
\firstpage{1262}
\lastpage{1311}

\makeatletter
\newcolumntype{d}[1]{D{.}{.}{#1}}
\renewcommand{\mid}{|}
\newtheorem{theorem}{Theorem}
\newtheorem{lemma}{Lemma}
\newproclaim{remark}{Remark}
\newproclaim{algo}{Algorithm}
\newtheorem{proposition}{Proposition}
\def\Exp{\operatorname{Exp}}
\def\CDF{\mathrm{CDF}}
\def\FDR{\mathrm{FDR}}
\def\FDP{\mathrm{FDP}}
\def\fMRI{\mathrm{fMRI}}
\def\SIM{\mathrm{SIM}}
\def\P{\mathrm{P}}
\def\ID{\mathrm{I}}
\def\I{\mathrm{I}}
\def\II{\mathrm{II}}
\def\BH{\mathrm{BH}}
\def\OR{\mathrm{OR}}
\def\GO{\mathrm{GO}}
\def\CN{\mathrm{CN}}
\def\MFA{\mathrm{MFA}}
\def\bbeta{{\bolds{\beta}}}
\def\bmu{{\bolds{\mu}}}
\def\boldsymbolX{\mathbf{X}}
\def\boldsymbolW{\mathbf{W}}
\def\boldsymbolZ{\mathbf{Z}}
\def\boldsymbolP{\mathbf{P}}
\def\boldsymbolt{\mathbf{t}}
\def\boldsymbolS{\mathbf{S}}
\makeatother

\begin{document}
\begin{frontmatter}

\title{Single-index modulated multiple testing}
\runtitle{Single-index modulated multiple testing}

\begin{aug}
\author[a]{\fnms{Lilun}~\snm{Du}\ead[label=e1]{dulilun@stat.wisc.edu}}
\and
\author[b]{\fnms{Chunming}~\snm{Zhang}\corref{}\ead[label=e2]{cmzhang@stat.wisc.edu}\thanksref{t1}}
\runauthor{L. Du and C. Zhang}
\affiliation{University of Wisconsin--Madison, and\\ Nankai University and University of Wisconsin--Madison}
\address[a]{Department of Statistics\\
University of Wisconsin\\
Madison, Wisconsin 53706\\
USA\\
\printead{e1}}
\address[b]{School of Mathematical Sciences\\
Nankai University\\
Tianjin 300071\\
China\\
and\\
Department of Statistics\\
University of Wisconsin\\
Madison, Wisconsin 53706\\
USA\\
\printead{e2}}
\end{aug}
\thankstext{t1}{Supported by NSF Grants DMS-11-06586 and
DMS-13-08872, and Wisconsin Alumni Research Foundation.}

\received{\smonth{12} \syear{2013}}
\revised{\smonth{3} \syear{2014}}

%
\begin{abstract}
In the context of large-scale multiple testing, hypotheses are often
accompanied with certain prior information.
In this paper, we present a single-index modulated (SIM) multiple
testing procedure, which maintains
control of the false discovery rate while
incorporating prior information, by assuming the availability of a
bivariate \mbox{$p$-}value,
$(p_1, p_2)$, for each hypothesis, where $p_1$ is a preliminary
\mbox{$p$-}value from prior information and
$p_2$ is the primary \mbox{$p$-}value for the ultimate analysis.
To find the optimal rejection region for the bivariate \mbox{$p$-}value, we
propose a criteria based on the ratio of probability density functions
of $(p_1, p_2)$ under the true null and nonnull.
This criteria in the bivariate normal setting further motivates us to
project the bivariate \mbox{$p$-}value to a single-index, $p(\theta)$,
for a wide range of directions $\theta$. The true null distribution of
$p(\theta)$ is estimated via
parametric and nonparametric approaches, leading to two procedures for
estimating and controlling the false discovery rate.
To derive the optimal projection direction $\theta$, we propose a new
approach based on power comparison, which is further
shown to be consistent under some mild conditions.
Simulation evaluations indicate that the SIM multiple testing procedure
improves the detection power significantly while controlling the false
discovery rate.
Analysis of a real dataset will be illustrated.
\end{abstract}

%
\end{frontmatter}

%
\section{Introduction} \label{Sec-1}
Large-scale simultaneous hypothesis testing problems, with thousands or
even tens of thousands of cases considered together,
have become a familiar feature in scientific fields such as biology,
medicine, genetics, neuroscience, economics and finance.
For example, in genome-wide association study, testing for association
between genetic variation and a complex disease typically requires
scanning hundreds of thousands
of genetic polymorphisms; in functional magnetic resonance imaging
($\fMRI$), time-course measurements
over $10^4$--$10^5$ voxels in the brain are typically available to allow
investigators to determine
which areas of the brain are involved in a cognitive task.
Multiple testing procedures, especially the false discovery rate ($\FDR
$) control method \cite{BenjaminiHochberg1995},
have been widely used to screen the massive data sets to identify a few
interesting cases.

In many real-world applications, the tests are accompanied with a
scientifically meaningful structure.
In $\fMRI$, each test corresponds to a specific brain location;
in microarray studies, each test is related to a specific gene.
These types of structural information usually provide valuable prior
information. For example,
previous studies may suggest that some null hypotheses are more or less
likely to be false;
similarly, in spatially-structured problems, nonnull hypotheses are
more likely to be clustered than true nulls.
It is thus anticipated that exploiting structural prior information
will improve
the performance of conventional multiple testing procedures. Several
attempts have been made
in the literature to incorporate prior information. For instance,
methods that up-weight or down-weight hypotheses appeared in \cite{BenjaminiHochberg1997,Genoveseetal2006} and \cite{HuZhaoZhou2010}.
A comprehensive review of weighted hypothesis testing can be found in
\cite{RoederWasserman2009} and the references therein.
A different approach, based on a two-stage approach mainly arising from
the microarray literature
\cite{Bourgonetal2010,HackstadtHess2009,Lusaetal2008,McClintickEdenberg2006,Talloenetal2007,Tritchleretal2009}, extracted the prior information
to remove a subset of genes which seem to generate uninformative
signals in the \emph{filtering} stage,
followed by applying some multiple testing procedure to the remaining genes
which have passed the filter in the \emph{selection} stage.

Very little work, however, has been published on
theoretically quantifying the extent to which the pair of filter and
test statistics in the above two-stage procedure, as well as
the pair of random weight and test statistics in weighted hypothesis
testing affect $\FDR$ and power. This issue is critically important,
because arbitrarily
choosing a filter (or weight) statistic may lead to loss of type I
error control. To guarantee the validity of \emph{filtering}
in the two-stage multiple testing procedure, \cite{Bourgonetal2010}~recommended the use of a filter statistic (i.e., overall sample
variance) which is independent of the test statistic to reduce the
impact that multiple testing
adjustment has on detection power. Analogously, the weight and test
statistics are assumed to be
independent in the literature of weighted hypothesis testing. However,
questions always arise about
(I) the {adequacy of the independence assumption} between the filter
(or weight) and test statistics, and
(II) the {subjectiveness in setting the proportion of hypotheses to be
removed in the \emph{filtering} stage}.

We intend to incorporate the prior information into large-scale
multiple testing, via a proposed
single-index modulated ($\SIM$) multiple testing procedure.
This inspires us to study a bivariate \mbox{$p$-}value $(p_{i1},p_{i2})$ for
each of the $i$th hypothesis, $i=1, \ldots, m$, where
$m$ is the number of hypotheses,
$p_{i1}$ is the preliminary \mbox{$p$-}value from the prior information (e.g., the filter or weight), and
$p_{i2}$ is the primary \mbox{$p$-}value for the ultimate analysis (from the
test statistic).
Unlike \cite{Bourgonetal2010} and \cite{Genoveseetal2006}, we do not
impose the independence assumption between the filter (or weight) and
test statistics.
This greatly broadens the scope of filters (or weights) that can be chosen.
Moreover, we wish to point out that \cite{Chi2008} explored a $\FDR$
procedure which can achieve the control of $\FDR$ with
asymptotically maximum power through nested regions of multivariate
\mbox{$p$-}values of test statistics.
However, that approach assumed independence between components in each
multivariate \mbox{$p$-}value under true null hypotheses,
thus is not directly applicable to our
study.

In our approach, the bivariate \mbox{$p$-}value in multiple testing is
projected into a single-index, $p(\theta)$,
where the direction $\theta$ takes value in the interval $[0, \pi/2]$.
Due to the projection, the true null distribution of the
single-index $p(\theta)$ is no longer uniform and thus needs to be estimated.
We propose a parametric and a nonparametric
approach to estimate it. A data-driven estimator based on power
comparison is developed for the optimal projection direction $\theta$.
This estimator is further shown to be consistent under some mild conditions.
The resulting method leads to the estimation and control of $\FDR$ for
the $\SIM$ multiple testing procedure.
Compared with the conventional multiple testing procedure which ignores
the prior information,
the $\SIM$ multiple testing procedure can improve the detection
power substantially as long as components in the bivariate \mbox{$p$-}value
are not highly positively correlated.
Extensive simulation studies support
the validity and detection power of our approach.
Analysis of a real dataset illustrates the practical utility of the
proposed $\SIM$ procedure.

The rest of the paper is organized as follows.
Section~\ref{Sec-2} reviews the conventional multiple testing
procedure, and outlines the proposed $\SIM$ multiple testing procedure.
Section~\ref{Sec-3} supplies theoretical derivation of the $\SIM$
multiple testing procedure.
Section~\ref{Sec-4} presents methods for estimating and controlling
$\FDR$ used in the $\SIM$ multiple testing procedure
and Section~\ref{Sec-5} investigates their theoretical properties.
Section~\ref{Sec-6} evaluates the performance of the proposed
procedure in simulation studies.
Section~\ref{Sec-7} analyzes a real dataset.
Section~\ref{Sec-8} ends the paper with a brief discussion.
All technical proofs are relegated to Appendices \ref{appA} and \ref{appB}.

%
\section{Overview of the single-index modulated multiple testing procedure} \label{Sec-2}
%
\subsection{Review of the conventional multiple testing procedure}\label{Sec-21}
For the sake of discussion, we begin with a brief review of the
conventional multiple testing procedure.
For testing a family of null hypotheses, $\{H_0(i)\}_{i=1}^{m}$, with
the corresponding \mbox{$p$-}values
$\{p_1,\ldots,p_m\}$, Table~\ref{Table-1}
describes the outcomes when applying some significance rule, which
means rejecting null hypotheses
with corresponding \mbox{$p$-}values less than or equal to some threshold.
The false discovery rate
($\FDR$), $\FDR= E(\frac{V}{R\vee1})$, depicts\vspace*{1pt}
the expected proportion of incorrectly rejected
null hypotheses \cite{BenjaminiHochberg1995}, where $R\vee1=\max\{
R, 1\}$.
An empirical process definition of $\FDR$,
\begin{eqnarray*}
\FDR(t)&=&E \biggl\{\frac{V(t)}{R(t) \vee1} \biggr\},\qquad t\in[0,1],
\end{eqnarray*}
was introduced by \cite{Storeyetal2004},
where
$V(t) = \#\{ \mbox{true null } p_i\dvtx  p_i \le t\}$, and $R(t)= \#\{p_i\dvtx
p_i \le t\}$.

%
\begin{table}[t] 
\tabcolsep=0pt
\tablewidth=250pt
\caption{Outcomes from testing $m$ null hypotheses based on a significance rule}\label{Table-1}
\begin{tabular*}{\tablewidth}{@{\extracolsep{\fill}}@{}lccc@{}}\hline
& \textbf{Retain null} & \textbf{Reject null} & \textbf{Total}\\
\hline
Null is true & $U$ & $V$ & $m_0$ \\
Nonnull is true & $T$ & $S$ & $m_1$ \\
Total & $W$ & $R$ & $m$ \\
\hline
\end{tabular*}
\end{table}

Compared with the frequentist framework of $\FDR$, $\FDR$ methods
also have a
Bayesian rationale in terms of the two-groups model. Let $F_{0}(t)$ and
$F_{1}(t)$ be
the cumulative distribution functions ($\CDF$) of a \mbox{$p$-}value under
the true null and nonnull, respectively,
and define $F(t)=\pi_0 F_{0}(t)+\pi_1 F_{1}(t)$ as its marginal $\CDF$,
where $\pi_0=\P$(null is true) and $\pi_1=1-\pi_0$.
Then the Bayes formula yields the posterior probability,
%
\begin{eqnarray}
\label{21} \mathrm{Fdr}(t) &=& \P(\mbox{true null} \mid p \le t)=
\frac{\pi_0 F_{0}(t) }{
\pi_0F_{0}(t)+\pi_1 F_{1}(t) }=\frac{\pi_0 F_{0}(t) }{ F(t) },
\end{eqnarray}
of a null hypothesis being true
given that its \mbox{$p$-}value is less than or equal to some threshold $t$.

Assuming that \mbox{$p$-}values under the true null are independent (or weakly
dependent) and uniformly distributed on the interval $[0,1]$,
\cite{Storey2002} proposed a point estimate of $\FDR$ by
%
\begin{equation}
\label{22} \widehat{ \FDR}(t)= \frac{m \hat {\pi}_0 t }{R(t) \vee1} =\frac{ \hat {\pi}_0 t }{\{R(t) \vee1\}/m}.
\end{equation}
For a chosen level $\alpha$, a data-driven threshold for the
\mbox{$p$-}values is determined by
%
\begin{equation}
\label{23} t_{\alpha} ( \widehat{ \FDR} )= \sup\bigl\{ 0 \le t \le1\dvtx
\widehat{ \FDR}(t) \le\alpha\bigr\}.
\end{equation}
Reject a null hypothesis if its \mbox{$p$-}value is less than or equal to
$t_{\alpha} ( \widehat{ \FDR} )$.
Hereafter, we will refer to (\ref{22}) as the estimation approach for
$\FDR$ and (\ref{23}) as the controlling approach for $\FDR$.
%
\subsection{Outline of the single-index modulated multiple testing}\label{Sec-22}
Before describing the details of our proposed single-index modulated
multiple testing,
we outline the major idea and methodology.
\begin{longlist}[(a)]
\item[(a)]
For each bivariate \mbox{$p$-}value $(p_{i1}, p_{i2})$, $i=1,\ldots,m$,
project it into a sequence of single indices,
$\{p_i(\theta_l) \}_{l=1}^{L}$,
according to $p_i(\theta)=\Phi( \cos(\theta)\Phi^{-1}(p_{i1})
+\sin(\theta)\Phi^{-1}(p_{i2}))$,
where $\{\theta_l \}_{l=1}^{L}$ are equally spaced on the interval
$[0, \pi/2]$.

\item[(b)] For each $\theta_l$,
estimate the true null distribution function of $\{p_i(\theta_l)\dvtx
i=1,\ldots, m\} $
by $\widehat {F}_0(t, \theta_l)$ using either a parametric or
nonparametric approach.

\item[(c)]
For each $\theta_l$, calculate ${ R (\hat t^*_{\alpha' }(\theta
_l), \theta_l )} $,
where $R (t, \theta)={\#\{ p_i(\theta) \le t \} } $, and
$ \hat {t}^*_{\alpha' } (\theta_l)
=\sup\{ 0 \le t \le1\dvtx  { m \widehat {F}_0(t, \theta_l) } /{ \{ R (t,
\theta_l ) \vee1 \} } \le\alpha' \} $,
with $\alpha' \in(0,1)$.
Determine the data-driven optimal projection direction $\hat
{\theta}(\alpha')=\theta_{L^*}$, where $L^*=\arg\max_{1 \le l \le
L} R ( \hat {t}^*_{ \alpha' } (\theta_l), \theta_l)$.

\item[(d)]
Estimate the proportion $\pi_0$ of true null hypotheses by $\hat {\pi}_0$.

\item[(e)]
For the projected \mbox{$p$-}values $\{ p_{i}( \hat {\theta}(\alpha'))\dvtx
i=1,\ldots, m \}$, set the
threshold $\hat {t}_{\alpha}$ to be $\sup\{ 0 \le t \le1\dvtx  { m
\hat {\pi}_0 \widehat {F}_0(t, \hat {\theta}(\alpha') ) } /
{ \{ R (t, \hat {\theta}(\alpha') ) \vee1 \} } \le\alpha\}$,\vspace*{1pt}
where $\alpha \in(0,1)$. Reject a null hypothesis $H_0(i)$ if the
corresponding
$p_i( \hat {\theta}(\alpha') )$ is less than or equal to
$\hat {t}_{\alpha}$.
\end{longlist}

The idea of the single-index projection in part (a) is not
straightforward, evolving from Sections~\ref{Sec-31} and~\ref{Sec-32}, to Section~\ref{Sec-33}.
Section~\ref{Sec-31} starts with an intuitive idea of using a
rectangular shape of the rejection region for bivariate \mbox{$p$-}values;
Section~\ref{Sec-32} derives a general form of optimal rejection
region using local false discovery rate \cite{EfronTibshirani2002};
Section~\ref{Sec-33} is motivated from the bivariate normal setting, where
the optimal rejection region in Section~\ref{Sec-32} will lead to the
projected \mbox{$p$-}value, that is, the single-index $p(\theta)$.
The parametric and
nonparametric estimators in part (b) will be given in Section~\ref{Sec-42}.
Incorporating this, the estimator for the proportion of true null
hypotheses in part (d)
is derived in Section~\ref{Sec-43}. The optimal projection direction
in part (c)
is estimated by a novel approach given in Section~\ref{Sec-44}.
The procedure in part (e) for estimation and control of the false
discovery rate is provided in Section~\ref{Sec-45}.
%
\section{Optimal rejection region for bivariate \mbox{$p$-}values} \label{Sec-3}
Recall that for univariate \mbox{$p$-}values, the rejection region is an
interval $[0, t]$. In this section,
we will discuss the rejection region for bivariate \mbox{$p$-}values and its
optimal choice.
%
\subsection{Optimal rejection region based on a rectangle} \label{Sec-31}
Intuitively, the false discovery rate for the bivariate \mbox{$p$-}values can
be defined based on
a rectangular \mbox{rejection} region, $[0, t_1]\times [0, t_2]$.
For notational simplicity, let $\mathbf{p}=(p_1, p_2)$ denote the
bivariate \mbox{$p$-}value, and
define $F_{0}( \mathbf{p})$, $F_{1}( \mathbf{p})$ and $F( \mathbf
{p})$ to be the
true null joint distribution,
nonnull joint distribution and joint distribution of $\mathbf{p}$,
respectively.
Also, let $f_{0}( \mathbf{p})$, $f_{1}( \mathbf{p})$ and $f( \mathbf
{p})$ be the corresponding probability density functions~(p.d.f.).
Then the Bayesian $\mathrm{Fdr}$ for the bivariate \mbox{$p$-}value based on
a rectangular rejection region is formulated as
%
\begin{equation}
\mathrm{Fdr}( \boldsymbolt) = \P(\mbox{true null}\mid\mathbf{p}\le
\boldsymbolt) = \frac{\pi_0 F_{0}( \boldsymbolt)}{ \pi_0F_{0}( \boldsymbolt)+\pi
_1F_{1}(\boldsymbolt) }, \label{31}
\end{equation}
where $\boldsymbolt=(t_1, t_2)$ and $\{ \mathbf{p}\le\boldsymbolt\}
$ denotes the event $\{ p_1\le t_1, p_2 \le t_2\}$.
There are infinite choices of rejection regions $[0, t_1] \times [0,
t_2]$ such that $\mathrm{Fdr}( \boldsymbolt) \le\alpha$.
A~possible criteria to choose ${\boldsymbolt}^*=(t_1^*, t_2^*)$ for a
best rejection region
is based on power comparison. Specifically, that choice is
%
\begin{equation}
\boldsymbolt^*=\arg\max_{ \boldsymbolt} \bigl\{F_{1}(
\boldsymbolt)\dvtx  \mathrm{Fdr}( \boldsymbolt) \le\alpha\bigr\}. \label{32}
\end{equation}

\begin{remark}
The Bayesian $\mathrm{Fdr}$ formula~(\ref{31}) can also be derived
using conditional probability,
%
\begin{eqnarray}\label{33}
\mathrm{Fdr}( \boldsymbolt) &=& \frac{\P(\mbox{true  null} \mid p_1 \le t_1) \P(p_2 \le
t_2\mid p_1 \le t_1, \mbox{true null})}{\P(p_2 \le t_2\mid p_1
\le t_1)}
\nonumber\\[-8pt]\\[-8pt]
&=& \frac{ \mathrm{Fdr}_{p_1}(t_1) \P(p_2 \le t_2\mid p_1 \le t_1,
\mbox{true null})}{\P(p_2 \le t_2\mid p_1 \le t_1)},\nonumber
\end{eqnarray}
where $\mathrm{Fdr}_{p_1}(t_1)=\P(\mbox{true null}\mid p_1 \le t_1)$.
From formula (\ref{33}), the Bayesian $\mathrm{Fdr}$ for the
bivariate \mbox{$p$-}value based on a rectangular rejection region
is not simply the product of those with respect to the preliminary
\mbox{$p$-}value and
primary \mbox{$p$-}value, that is, $\mathrm{Fdr}( \boldsymbolt) \ne\mathrm
{Fdr}_{p_1}(t_1)\times \mathrm{Fdr}_{p_2}(t_2) $,
where $\mathrm{Fdr}_{p_2}(t_2)=\P(\mbox{true null}\mid p_2 \le t_2)$.
Furthermore, formula~(\ref{33}) provides an insight into the
two-stage multiple testing in
\cite{Bourgonetal2010} if $p_1$ is utilized as the filter in the
\emph{filtering} stage and $p_2$ is obtained from a test statistic
in the \emph{selection} stage.
Comparing (\ref{33}) with (\ref{21}), we find that
$\mathrm{Fdr}_{p_1}(t_1)$ in the \emph{filtering}
stage is the proportion of
the true null hypotheses
served in the \emph{selection} stage. On the one hand, in order to
improve the power in the \emph{selection} stage,
we can control $\mathrm{Fdr}_{p_1}(t_1)$ to be small.
On the other hand, increasing $\mathrm{Fdr}_{p_1}(t_1)$ will assure
that we do not screen out too many nonnull hypotheses from the \emph
{filtering} stage.
\end{remark}
%
\subsection{General form of optimal rejection region} \label{Sec-32}
In Section~\ref{Sec-31}, we observe that among infinite choices of
rectangular rejection regions $[0,t_1]\times[0, t_2]$
such that the Bayesian $\mathrm{Fdr}$ is less than or equal to $\alpha
$, there exists one ``best'' rectangle $[0,t_1^*]\times[0, t_2^*]$
with highest power.
In this section, we seek a general form of optimal rejection region, by
relaxing the shape of rejection region.
Let $\boldsymbolS$ denote a rejection region.
Following (\ref{31}), the Bayesian $\mathrm{Fdr}$ can be generalized to
%
\begin{equation}
\mathrm{Fdr}( \boldsymbolS)= \P(\mbox{true null} \mid \mathbf{p} \in
\boldsymbolS) = \frac{\pi_0 F_{0}( \boldsymbolS)}{\pi_0 F_{0}( \boldsymbolS)
+\pi_1 F_{1}( \boldsymbolS)}, \label{34}
\end{equation}
where $F_{j}( \boldsymbolS)=\int_\boldsymbolS f_{j}( \mathbf
{p})\,d\mathbf{p}$, $j=0,1$.
An optimal rejection region $\boldsymbolS^*$ is based on the following
definition:
%
\begin{equation}
\label{35} \boldsymbolS^*=\arg\max_{ \boldsymbolS}\bigl
\{F_{1}( \boldsymbolS)\dvtx  \mathrm{Fdr}( \boldsymbolS) \le\alpha\bigr\}.
\end{equation}
Note that (\ref{32}) is a special case of (\ref{35}), by
restricting $\boldsymbolS$ to be rectangular.
%
\begin{proposition} \label{Proposition-1}
Assume the two-groups model holds for the bivariate \mbox{$p$-}values and
let $\mathrm{fdr}(\mathbf{p})=\pi_0 f_{0}(\mathbf{p})/\{\pi_0
f_{0}(\mathbf{p})+\pi_1 f_{1}(\mathbf{p}) \}$ be the generalization
of local false discovery rate; see \cite{EfronTibshirani2002} and
\cite{Efronetal2001}.
Further suppose that for any constant~$C_0$,
%
\begin{equation}
\label{36} \mathrm{P} \bigl( \mathbf{p}\dvtx  \mathrm{fdr}(\mathbf{p})=C_0
\bigr)=0.
\end{equation}
Denote by $\boldsymbolS_{\OR}$ the rejection region to be formed by
$\boldsymbolS_{\OR}=\{ \mathbf{p}\dvtx  \mathrm{fdr}(\mathbf{p}) \le C
\}$,
where $C$ is a constant such that $\mathrm{Fdr}(\boldsymbolS_{\OR}
)=\alpha$.
Then for any rejection region $\boldsymbolS$ satisfying $\mathrm
{Fdr}(\boldsymbolS) \le\alpha$,
we have $F_{1}( \boldsymbolS) \le F_{1}(\boldsymbolS_{\OR} )$.
\end{proposition}

From Proposition~\ref{Proposition-1}, the general form of optimal
rejection region (\ref{35}) can be equivalently described as follows: within
the rejection region $\boldsymbolS^*$, the local false discovery rate
$\mathrm{fdr}( \mathbf{p})$
should be less than or equal to some threshold, which is equivalent to setting
${f_{1}( \mathbf{p})}/{f_{0}( \mathbf{p})}$ to be larger than or
equal to some threshold.\vadjust{\goodbreak}
Thus, we propose the optimal rejection region (\ref{35}) to be formed by
%
\begin{equation}
\label{37} \boldsymbolS^*=\bigl\{ \mathbf{p}\dvtx  f_{1}(
\mathbf{p})/f_{0}(\mathbf{p}) \ge C\bigr\},
\end{equation}
where $C$ is a constant such that $\mathrm{Fdr}( \boldsymbolS^*
)=\alpha$.
%
\begin{remark}
In traditional hypothesis testing,
the Neyman--Pearson lem\-ma indicates that the rejection region of the
uniformly most powerful
(UMP) test is in the form of
likelihood ratio of test statistics if both the null and nonnull
hypotheses are simple.
Hence, the form of the optimal rejection region $\boldsymbolS^*$ using
local false discovery rate
is similar to that derived from the UMP test.
(\ref{37}) is also a homogeneous version of the optimal discovery
procedure proposed by \cite{Storey2007}, where
the null and nonnull distributions across the tests are less
homogeneous and strongly correlated.
\end{remark}
%
\subsection{Optimal rejection region under bivariate normality} \label{Sec-33}
In this subsection, we will first derive the true null and nonnull
distributions of a bivariate \mbox{$p$-}value under bivariate normality,
followed by approximating the shape of the optimal rejection region
using criteria~(\ref{37}).

Efron \cite{Efron2007} introduced a $z$-value $ ( \Phi^{-1}(p)  )$
into traditional multiple testing problem and assumed that the
empirical null distribution of $z$-value is normal with mean $\mu$ and
standard deviation $\sigma$.
To derive an explicit form of the true null distribution of $\mathbf{p}$,
we borrow the idea of empirical null distribution in \cite{Efron2007} and
make extension to the case of bivariate \mbox{$p$-}values, assuming the
bivariate normality as follows:
\begin{longlist}
\item[(N1)] Under the true null hypothesis,
the transformed \mbox{$p$-}value $ ( \Phi^{-1}(p_1),\break  \Phi^{-1}(p_2)
 )$ follows a bivariate normal distribution $\mathcal{N}(\bmu
_0, \Sigma_0)$, where
%
\begin{equation}
\label{38} \bmu_0 =( \mu_{0; 1}, \mu_{0; 2}
)^T,\qquad \Sigma_0= \pmatrix{ \sigma_{0; 1}^2
& \rho_0 \sigma_{0; 1} \sigma_{0; 2}
\vspace*{3pt}\cr
\rho_0 \sigma_{0; 1} \sigma_{0; 2} &
\sigma_{0; 2}^2 }.
\end{equation}
\item[(N2)]
Under the nonnull, the transformed \mbox{$p$-}value
$(\Phi^{-1}(p_1), \Phi^{-1}(p_2))$
also follows a bivariate normal distribution $\mathcal{N}(\bmu_1,
\Sigma_1)$, where
%
\begin{equation}
\label{39} \bmu_1 =( \mu_{1;1}, \mu_{1;2}
)^T,\qquad \Sigma_1= \pmatrix{ \sigma_{1;1}^2
& \rho_1 \sigma_{1;1} \sigma_{1; 2}
\vspace*{3pt}\cr
\rho_1 \sigma_{1;1} \sigma_{1; 2} &
\sigma_{1; 2}^2 }.
\end{equation}
\end{longlist}

%
\begin{remark}
The assumption \textup{(N1)} is strictly satisfied if the components
of bivariate \mbox{$p$-}value are independent under the true null. For the
dependence case, this assumption is approximately true.
As a specific example, \textup{(N1)} holds if the preliminary test
statistic and primary test statistic
(bivariate test statistic) under the true null follows a bivariate
normal distribution for one-sided hypotheses; see~(\ref{B2}) in
Appendix~\ref{appB}. The assumption \textup{(N2)} is not required
for the general theory in Section~\ref{Sec-5} and only serves as a
motivation for developing the proposed rejection region~(\ref{313}).
\end{remark}

If \textup{(N1)} and \textup{(N2)} hold, some algebraic calculations
yield the densities of $\mathbf{p}$ under the true null and nonnull,
%
\begin{eqnarray}\label{310}
\qquad f_{0}( \mathbf{p}) &=& \frac{1} {\sigma_{0; 1} \sigma_{0; 2} \sqrt{1-\rho_0^2}} \exp \biggl(
\frac{ \{\Phi^{-1}(p_1)\}^2+ \{\Phi^{-1}(p_2)\}^2 }{2} \biggr)\nonumber
\\
&&{}\times
\exp \biggl( - \biggl( \biggl\{ \frac{\Phi^{-1}(p_1)-\mu_{0; 1} }{\sigma
_{0; 1}} \biggr\}^2 +
\biggl\{ \frac{\Phi^{-1}(p_2)-\mu_{0; 2} }{\sigma_{0; 2}} \biggr\}^2\nonumber
\\
&&\hspace*{49pt}{} -2 \rho_0
\biggl\{\frac{ \Phi^{-1}(p_1)-\mu_{0; 1} }{\sigma_{0; 1}}\biggr\}\biggl\{ \frac{\Phi
^{-1}(p_2)-\mu_{0; 2} }{\sigma_{0; 2} }\biggr\} \biggr)\nonumber
\\
&&\hspace*{197.5pt} {} \Big/ \bigl({2\bigl(1-\rho_0^2\bigr)}\bigr)\biggr),
\nonumber\\[-8pt]\\[-8pt]
f_{1}( \mathbf{p}) &=& \frac{1} {\sigma_{1;1} \sigma_{1; 2} \sqrt{1-\rho_{1}^2}} \exp \biggl(
\frac{ \{\Phi^{-1}(p_1)\}^2+ \{\Phi^{-1}(p_2)\}^2 }{2} \biggr) \nonumber
\\
&&{} \times
\exp \biggl( - \biggl( \biggl\{ \frac{\Phi^{-1}(p_1)-\mu_{1;1} }{\sigma
_{1;1}} \biggr\}^2 +\biggl\{ \frac{\Phi^{-1}(p_2)-\mu_{1; 2} }{\sigma_{1; 2}} \biggr\}^2\nonumber
\\
&&\hspace*{49pt}{}  -2 \rho_1
\biggl\{\frac{ \Phi^{-1}(p_1)-\mu_{1;1} }{\sigma_{1;1}}\biggr\}\biggl\{ \frac{\Phi
^{-1}(p_2)-\mu_{1; 2} }{\sigma_{1; 2}}\biggr\} \biggr)\nonumber
\\
&&\hspace*{198pt} {} \Big/ \bigl(2\bigl(1-\rho_{1}^2\bigr)\bigr) \biggr).\nonumber
\end{eqnarray}

By combining (\ref{310}) with the criteria~(\ref{37}), the optimal
rejection region under bivariate normality
takes the form
%
\begin{equation}
\boldsymbolS^* = \bigl\{ \mathbf{p}\dvtx  {\boldsymbolZ}^T \bbeta\ge C
\bigr\}, \label{311}
\end{equation}
with a constant $C$ such that $\mathrm{Fdr}(\boldsymbolS^*)=\alpha$, where
\[
{\boldsymbolZ}= \bigl( \bigl\{\Phi^{-1}(p_1)\bigr
\}^2, \bigl\{\Phi^{-1}(p_2)\bigr\}^2,
\Phi^{-1}(p_1) \Phi^{-1}(p_2),
\Phi^{-1}(p_1), \Phi^{-1}(p_2) \bigr)
\]
and $\bbeta$ is the corresponding vector of coefficients determined by
$\bmu_0$, $\bmu_1$, ${\Sigma_0}$ and~${\Sigma_1}$.

If\vspace*{1pt} the covariance matrices satisfy ${\Sigma_0}= { \Sigma_1}$, the
optimal rejection region
in~(\ref{311}) can be formulated in term of a single-index $\beta
_1\Phi^{-1}(p_1)+\beta_2 \Phi^{-1}(p_2)$, where
$(\beta_1,\beta_2)$ is determined by $\bmu_0$, $\bmu_1$, ${\Sigma_0}$.
This is more intuitive than the form (\ref{311}) from two perspectives.
From dimension reduction viewpoint, researchers always prefer reducing
the number of variables to choosing ${\boldsymbolZ}$.
From principal component analysis aspect, the transformed \mbox{$p$-}value
$(\Phi^{-1}(p_1), \Phi^{-1}(p_2) )$ can be visualized from two
orthogonal directions.
Instead of searching for the eigenvectors of common covariance matrix
${ \Sigma_0 }$, our goal is to find
a direction $(\beta_1, \beta_2)$, such that the projected points
corresponding to the true null
hypotheses deviate from those corresponding to the true nonnull as far
as possible.
Then (\ref{37}) will prompt us to introduce a ``single-index \mbox{$p$-}value,''
%
\begin{equation}
\label{312} p(\theta)=\Phi \bigl( \cos(\theta) \Phi^{-1}(p_1)+
\sin(\theta) \Phi^{-1}(p_2) \bigr),
\end{equation}
where $0 \le\theta\le{ \pi}/{2}$ acts as a tuning parameter.
This in turn yields our proposed rejection region [which is optimal
under \textup{(N1)} and \textup{(N2)}] defined as
%
\begin{equation}
\label{313} \boldsymbolS^{*}(\theta)=\bigl\{ \mathbf{p}\dvtx  p(\theta)
\le t\bigr\},
\end{equation}
where the threshold $t$ is chosen to control $\FDR$. We call this the
``single-index modulated ($\SIM$) multiple testing procedure.''

%
\begin{figure}[t] 

\includegraphics{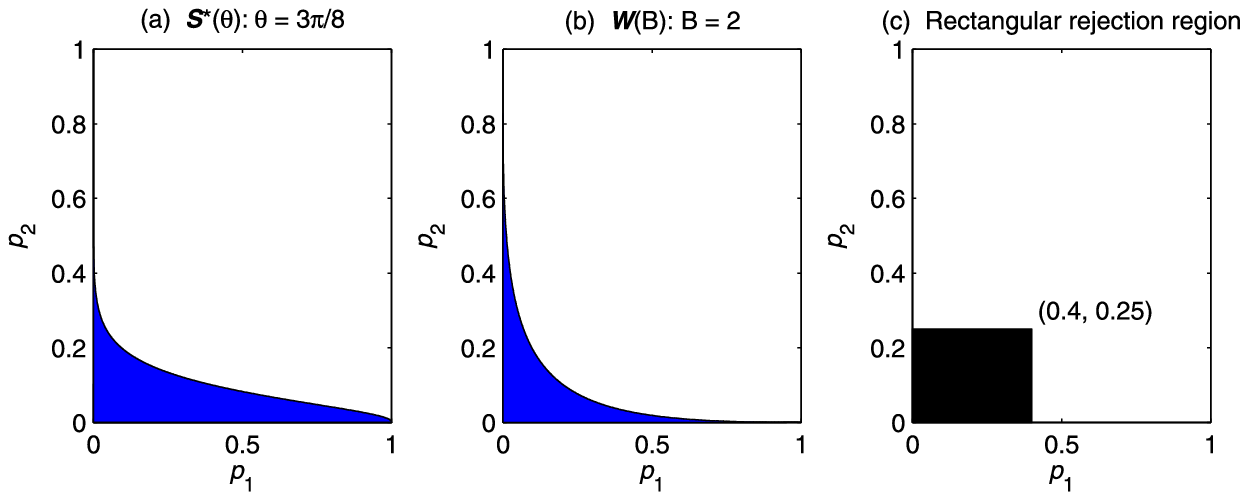}

\caption{Compare shapes of rejection regions $\boldsymbolS
^*(\theta)$, $\boldsymbolW(B)$ and a rectangle. Here, the threshold
$t=0.1$ is used in $\boldsymbolS^*(\theta)$ and $\boldsymbolW(B)$.}\label{Figure-1}
\end{figure}

As a comparison, the shape of the rejection region $\boldsymbolS
^*(\theta)$ is different from the rectangle used in the two-stage
multiple testing procedure of \cite{Bourgonetal2010}; see Figure~\ref{Figure-1}.
In addition, the philosophy underlying the two procedures varies.
For the two-stage procedure, a multiple testing procedure is
only applied to the subset of hypotheses survived from the \emph
{filtering} stage.
In contrast, the proposed $\SIM$ procedure does not screen any
hypotheses out, but projects
the bivariate \mbox{$p$-}value into a single-index $p(\theta)$. After that,
methods in Section~\ref{Sec-4} for estimation and
control of $\FDR$ are implemented using all the $m$ hypotheses.

To draw connection to the weighted multiple testing procedure of \cite
{Genoveseetal2006}, we first generate the weights from the preliminary
\mbox{$p$-}values and then combine the primary \mbox{$p$-}values with the weights. To
be specific, in the first stage, we generate \emph{cumulative
weights} \cite{Roederetal2006} proportional to $\{ v_i=\Phi( \Phi
^{-1}(1-p_{i1})-B )\dvtx  i=1,\ldots,m \}$, where $B$ is a tuning parameter.
Because the weights are constrained to
have mean~1, $w_i=v_i/\bar{v}_m$ is a valid choice, where
$\bar{v}_m=\sum_{i=1}^{m} \Phi( \Phi^{-1}( 1-p_{i1})-B )/m$.
In the second stage, standard $\BH$ procedure \cite{BenjaminiHochberg1995} is applied to the weighted \mbox{$p$-}values, that
is, $\{ {p_{i2} } / {w_i}\dvtx  i=1,\ldots,m\}$.
The rejection region of the weighted multiple testing procedure is
formed by
$\boldsymbolW(B)=\{\mathbf{p}\dvtx  p_2 \le\frac{ \Phi( \Phi
^{-1}(1-p_1)-B ) }{\bar{v}_m } t \}$; see\vspace*{2pt} Figure~\ref{Figure-1}
for the graphical illustration. Surprisingly, the $\SIM$ multiple
testing procedure and the weighted multiple testing procedure share
similar patterns of rejection.

\section{Estimation and control of $\FDR$ for the $\SIM$ procedure}\label{Sec-4}
In this section, we will first investigate properties of the
single-index $p(\theta)$,
followed by utilizing these properties to estimate and control the
false discovery rate.
For each possible direction $\theta$, denote by
$\{p_i(\theta)\dvtx  i=1, \ldots, m \}$ the sequence of projected
\mbox{$p$-}values. Let $F_0(t, \theta)$, $F_1(t, \theta)$
and $F(t, \theta)$ be the true null distribution,
nonnull distribution and marginal distribution of $p(\theta)$,
respectively. Similarly,
$f_0(t,\theta)$, $f_1(t, \theta)$ and $f(t, \theta)$ are their
corresponding density functions.
Following notations in Section~\ref{Sec-21}, the frequentist $\FDR$
and Bayesian $\mathrm{Fdr}$ for the projected \mbox{$p$-}values
$\{ p _i(\theta)\dvtx  i=1, \ldots, m \}$ are defined by
\begin{eqnarray*}
\FDR(t, \theta) &=& E \biggl\{ \frac{V(t, \theta)}{ R(t, \theta)\vee1 } \biggr\},
\\
\mathrm{Fdr}(t,
\theta) &=& \frac{\pi_0 F_0(t, \theta)} {F(t, \theta)},
\end{eqnarray*}
respectively,
where $V(t, \theta)=\#\{\mbox{true null } p_i(\theta)\dvtx  p_i(\theta)
\le t\}$ and
$R(t, \theta)=\break \#\{p_i(\theta)\dvtx  p_i(\theta) \le t\}$ are the number
of hypotheses erroneously rejected and
the number of hypotheses rejected, based on some significance rule for
the sequence of projected \mbox{$p$-}values.
%
\subsection{Property of the single-index \texorpdfstring{$p(\theta)$}{p(theta)}} \label{Sec-41}
The true null distribution of $p(\theta)$ in~(\ref{312}) plays an
important role in estimating the
false discovery rate. From the theory of statistics,
the theoretical true null distributions of $p_1$ and $p_2$ are uniform.
In the special case where $p_1$ and $p_2$ are independent under the
true null,
it is straightforward to show that $p(\theta)$ under the true null
also follows a uniform distribution.
In general, the assumption \textup{(N1)} with
$\bolds{\mu}_{0}=\mathbf{0}$ facilitates us to derive the
$\CDF$ of $p(\theta)$ under the true null hypothesis.
To be specific,
\begin{eqnarray*}
F_0(t, \theta) &=& \P\bigl( p(\theta) \le t \mid \mbox{true null}
\bigr) = \Phi \biggl(\frac
{ \Phi^{-1}(t) }{ \sigma_0(\theta) } \biggr),
\end{eqnarray*}
where
%
\begin{equation}\label{41}
\qquad \sigma_0 (\theta) =\sqrt{ \bigl\{\cos(
\theta)\bigr\}^2\sigma_{0;1}^2+ \bigl\{\sin(
\theta)\bigr\}^2 \sigma_{0;2}^2 + 2
\rho_0 \sigma_{0;1} \sigma_{0;2} \cos(\theta) \sin(
\theta) },
\end{equation}
and $\sigma_{0;1}$, $\sigma_{0;2}$ and $\rho_0$ are as defined in
(\ref{38}).
The following two categories summarize some properties of $ p(\theta) $.
\begin{longlist}[(II)]
\item[(I)]
If $\sigma_0 (\theta) =1$, $ p(\theta)$ under the true null hypothesis
follows a standard uniform distribution.

\item[(II)]
If $\sigma_0 (\theta) \ne1$, the true null distribution of $
p(\theta) $
is not uniform but symmetric with respect to ${1}/{2}$.
\end{longlist}

If $p_1$ and $p_2$ are both uniformly distributed under the true null,
the expression of $\sigma_0(\theta)$ can
be further simplified to $\sqrt{1+\rho_0 \sin(2\theta)}$. Under the
independence assumption (i.e., $\rho_0=0$),
$p(\theta)$ is uniform for all $\theta$, which belongs to category (I).
The case of negative correlation (i.e., $\rho_0<0$) implies $\sigma
_0(\theta)<1$,
shrinking most of the projected points corresponding to the true null
concentrating around the point $1/2$. Consequently, this case has
better potential to be powerful.
The positive correlation worsens the structure of \mbox{$p$-}values a little,
shifting some of the combined \mbox{$p$-}values corresponding to the true null
to the area adjacent to 0 or 1, but it is still
symmetric with respect to $1/2$. \cite{ZhangFanYu2011}
employed a $p^{*}$-value, the median of \mbox{$p$-}values in the neighborhood
of the original \mbox{$p$-}value, to capture the geometric feature in brain
imaging. The true null distribution of $p^*$ is beta, which is
symmetric with respect to ${1}/{2}$. Thus, the pair $(p, p^*)$ of
\mbox{$p$-}values belongs to category (II).

Although the assumption \textup{(N1)}
is imposed when deriving the specific form of $F_0(t,\theta)$, we
could relax
the normality assumption by assuming that the true null distribution of
$p(\theta)$ is symmetric about $1/2$ for all $\theta$. The symmetry
property assumption can be equivalently stated as:
\begin{longlist}[(N3)]
\item[(N3)]
The probability density function of $\mathbf{p}$ under the true null
is centrally symmetric with respect to the point $(1/2,1/2)$, that is,
$f_0(p_1, p_2)=f_0(1-p_1,1-p_2)$.
\end{longlist}
\textup{(N3)} provides flexibility in accommodating a wider range of
distributions for $\mathbf{p}$. For example, \textup{(N3)} holds if
the bivariate test statistic under the true null follows a bivariate
$t$ distribution for one-sided hypotheses; see (\ref{B3}) in Appendix~\ref{appB}.
In addition to estimating the parameter $\sigma_0(\theta)$,
Section~\ref{Sec-42} will develop an adaptive data-driven estimator
for $F_0(t, \theta)$
using a nonparametric approach based on \textup{(N3)}.
While this relaxed assumption causes certain loss in efficiency for
estimating $F_0(t, \theta)$, it achieves a gain in robustness.
%
\subsection{Estimating the true null distribution of \texorpdfstring{$p(\theta)$}{p(theta)}}\label{Sec-42}
Recall the properties of $ p(\theta)$ in Section~\ref{Sec-41}. If
the normality assumption \textup{(N1)} holds, one can estimate the
true null distribution of $p(\theta)$ using the following parametric approach:
%
\begin{equation}
\label{42} \widehat F_0^{ \I} (t, \theta)=\Phi \biggl(
\frac{ \Phi^{-1}(t)}{
\hat{\sigma}_0(\theta) } \biggr),
\end{equation}
where $ \hat{\sigma}_0(\theta)$ stands for some parametric
estimator of $\sigma_0(\theta)$.
Here, we will provide a simple and efficient estimator in the following
procedure:
\begin{longlist}[(a)]
\item[(a)]
Select a constant $c \ge0$, such that $z$-values, $z(\theta) = \Phi
^{-1}(p(\theta))$, from $(-c,\infty]$
are more likely to come from the true null hypothesis.
\item[(b)]
Split the data $\{z_i(\theta)\dvtx  i=1,\ldots,m\}$
into\vspace*{1pt} three parts, that is,
$\widetilde{Z}_{[-\infty,-c]}$,
$\widetilde{Z}_{(-c,c]}$, and
$\widetilde{Z}_{(c,\infty]}$, where the notation $\widetilde{Z}_{I}$
denotes the sample
from interval $I$. Here, $I$ can be a closed, open or half-open interval.
\item[(c)]
Drop the sample $\widetilde{Z}_{[-\infty, -c]}$ and impute
$-\widetilde{Z}_{ [c, \infty] }$ into
the interval\break  $[-\infty, -c]$. $\hat{\sigma}_0(\theta) $ is\vspace*{1pt} the
standard error of the newly constructed data
$ \widetilde{Z}{}^{*}=\{ -\widetilde{Z}_{ [c, \infty) }, \widetilde
{Z}_{ (-c, c] },
\widetilde{Z}_{ (c, \infty) }\}$.
\end{longlist}

If the normality assumption \textup{(N1)} is violated, we provide a
nonparametric estimator based on
the assumption \textup{(N3)}.
The nonparametric approach follows the idea of \cite{ZhangFanYu2011}.
To be specific, $F _0(t, \theta)$ can be estimated by the empirical
distribution function,
%
\begin{eqnarray}\label{43}
\widehat{ F} _0^{\II} (t, \theta) &=& \cases{
\displaystyle\frac{ \sum_{i=1}^{m} \ID\{ p_i(\theta) \ge(1-t)\}}{
2\sum_{i=1}^{m} \ID\{ p_i(\theta) >0.5\}+ \sum_{i=1}^{m} \ID\{
p_i(\theta) =0.5\}},
\vspace*{5pt}\cr
\qquad\mbox{if $0 \le t \le0.5$},
\vspace*{6pt}\cr
\displaystyle1-\frac{ \sum_{i=1}^{m} \ID\{ p_i(\theta) \ge t\} }{
2\sum_{i=1}^{m} \ID\{ p_i(\theta) >0.5\}+ \sum_{i=1}^{m} \ID\{
p_i(\theta) =0.5\}},
\vspace*{5pt}\cr
\qquad\mbox{if $0.5 < t\le1$.}}
\end{eqnarray}
%
\subsection{Estimating the proportion \texorpdfstring{$\pi_0$}{pi0} of true null hypotheses} \label{Sec-43}
There is an active research pursued in estimating $\pi_0$
(e.g., \cite{BenjaminiHochberg2000,Benjaminietal2006,HochbergBenjamini1990,LiangNettleton2012,SchwederSpjtvoll1982,Storey2002,Storeyetal2004}).
\cite{Storey2002} and \cite{Storeyetal2004} proposed an estimator
$\hat {\pi}_0 (\lambda)=\{m-R(\lambda) \}/\{ (1-\lambda)m \}$
with a tuning parameter
$\lambda$ in $[0,1)$ to be specified.
\cite{LiangNettleton2012} summarized many adaptive and dynamically
adaptive procedures
for estimating $\pi_0$ and proposed a unified dynamically adaptive procedure.
In this paper, we follow the same principle in \cite
{LiangNettleton2012} and propose
two estimators of $\pi_0$ dynamically according to two estimators of
the true null distribution of $p(\theta)$
proposed in (\ref{42}) and (\ref{43}), respectively,
%
\begin{eqnarray}\label{44}
\hat{\pi}{}^ {\I} _0(\theta) &=&\frac{ m- R( \hat{\lambda} {}^{ \I} (\theta), \theta) }{
\{1-\widehat F_0^{\I} ( \hat{\lambda} {}^{ \I} (\theta),
\theta)\}m },
\nonumber\\[-8pt]\\[-8pt]
\hat{\pi}{}^{ \II }_0( \theta) &=&\frac{ m- R( \hat{\lambda} {}^{\II} (\theta), \theta) }{
\{1-\widehat F_0^{ \II}( \hat{\lambda}{}^{ \II} (\theta),
\theta)\}m },\nonumber
\end{eqnarray}
where $ \hat{ \lambda} {}^{ \I} (\theta)$ and $ \hat{
\lambda} {}^{ \II} (\theta)$
are dynamically chosen as in the algorithm below.
\begin{algo*}[(For choosing $\lambda$)]
For a sequence of values $0 \equiv \lambda_0 < \lambda_1<\cdots
<\lambda_n \le1/2$,
$ \hat { \lambda} (\theta) $ is chosen to be $ \lambda_{ I^{*}
}$, where
$ I^* =\min\{ 1 \le j \le n-1\dvtx\break   \hat {\pi}_0 (\lambda_j,
\theta)
\ge \hat {\pi}_0 ( \lambda_{j-1}, \theta)\}$
if $ \hat {\pi}_0 ( \lambda_j, \theta) \ge \hat { \pi}_0
(\lambda_{j-1}, \theta)$ for some
$ j=1,\ldots,n-1 $ and $\lambda_{I^*} = \lambda_n$ otherwise. Here,
$\hat {\pi}_0 (\lambda, \theta) $ is\vspace*{1pt} defined as
${ \sum_{i=1}^{m} \ID\{ p_i (\theta) >\lambda\} } /\break  { [\{ 1-
\widehat F_0 (\lambda, \theta)\} m ]}$,
where the estimator $\widehat F_0$ can be either (\ref{42}) or (\ref
{43}) for
the $\CDF$ of $p(\theta)$ under the true null hypothesis.
\end{algo*}
%
\begin{remark}
We make the remarks concerning the algorithm.
\begin{itemize}
\item
The range $(0, 1/2]$ of the sequence of values $\{\lambda_i\dvtx  i=1,
\ldots, n\}$ is different from that in the right boundary procedure proposed
by \cite{LiangNettleton2012}, where $\lambda$ can be loosely selected
from $[0, 1)$. We restrict the
range to ${1}/{2}$ from two perspectives. On the one hand, it can be
verified that $\hat {\pi}_0^{ \II}(\lambda, \theta)$
is a constant for all $\lambda\ge1/2$ and~$\theta$. On the other
hand, condition \textup{(C5)} in\vspace*{1pt}
Appendix~\ref{appA} that $F_1(1/2, \theta)=1 $ for all~$\theta$,
guaranteeing the consistency of $\widehat F^{\II}_{0}(t, \theta)$,
enables us to search for $\lambda$ in a narrower range, which will be
more efficient in practice.

\item
Theoretically, it is equivalent to get $\lambda(\theta) $
as $\lambda(\theta) =\inf_{ 0 \le t \le{1}/{2} } \{ t\dvtx  F_1(t,
\theta)=1\}$.
If $t \le\lambda(\theta) $, there is an upward-bias for estimating
$\pi_0$, that is,
$\pi_1 \times \frac{1-F_1(t, \theta) }{1-F_0(t, \theta)}$;
if $t> \lambda(\theta) $, the variance of $\hat{\pi}_0 (t,
\theta)$ is proportional
to $ {1}/[ \{1-F_0(t, \theta)\}^2m ]$.
Instead of estimating $F_1(t, \theta)$, the algorithm described in the
algorithm paragraph provides
a rough but simple approach to estimate $\lambda(\theta)$.
Here, we would like to point out that fixing $\lambda$ is not
applicable to our approach,
since $\lambda(\theta)$ varies with the tuning parameter $\theta$.
\end{itemize}
\end{remark}
%
\subsection{Selection of projection direction \texorpdfstring{$\theta$}{theta}} \label{Sec-44}
\begin{figure}

\includegraphics{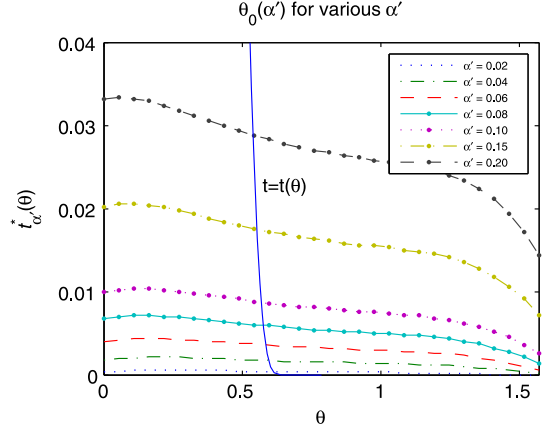}

\caption{Illustrate the optimal projection direction $\theta
_{0}(\alpha')$ in (\protect\ref{45}) for various choices of $\alpha
'$ when $( \Phi^{-1}(p_1), \Phi^{-1}(p_2) )$ follows (\protect\ref
{38}) with $\bmu_0=\mathbf{0}$, $\sigma_{0;1}=\sigma_{0;2}=1$,
$\rho_0=0.2$ under the true null,
and follows (\protect\ref{39}) with $\bmu_1=(-2, -1.5)^{T}$,
$\sigma_{1;1}=\sigma_{1;2}=1$, $\rho_1=0.6$ under nonnull, respectively.
The solid line is the implicit curve $t=t(\theta)$ satisfying \protect
(\protect\ref{46}) in Proposition~\protect\ref{Proposition-2}.
From the proof of Proposition~\protect\ref{Proposition-2}, the
$x$-coordinates of the intersection points are $\theta_0(\alpha')$.}\label{Figure-2}
\end{figure}

A specific $\theta$ corresponds to a projection direction, $( \cos
(\theta), \sin(\theta) )$ in (\ref{312}), for the transformed \mbox{$p$-}value
$( \Phi^{-1}(p_1), \Phi^{-1}(p_2) )$. The choice of $\theta=0$
amounts to utilizing $p_1$ alone,
whereas setting $\theta=\pi/2$ is equivalent to making inference with
the information from $p_2$ alone.
This indicates that our method indeed generalizes the conventional
multiple testing.
Recalling the shape of rejection region (\ref{313}) and the criteria
(\ref{35}), different values of
$\theta$ correspond to different shapes of rejection regions and the
one with the highest power is preferred.
Denote by $\theta_0(\alpha')$ the optimal value of $\theta$, that is,
%
\begin{equation}
\label{45} \theta_0\bigl(\alpha'\bigr)=\arg\max
_{ 0 \le\theta\le\pi/2 } F_1\bigl(t^*_{\alpha' }(\theta), \theta
\bigr),
\end{equation}
where $t^{*}_{\alpha' }(\theta) = \sup\{ 0 \le t \le1\dvtx   { F_0(t,
\theta)} / { F(t, \theta) } \le\alpha' \}$ and $0<\alpha'<1$.
The threshold $t^*_{\alpha' }(\theta)$ in criteria~(\ref{45}) is
chosen such that $\mathrm{Fdr}$ with respect to
various $\theta$ is controlled at level $\pi_0 \alpha'$.
%
\begin{proposition} \label{Proposition-2}
Suppose\vspace*{2pt} that $F_0(t, \theta)$ and $F_1(t, \theta)$ are continuously
differentiable and
$\frac{\partial F_1(t,\theta)}{\partial t}- \beta\frac{\partial
F_0(t, \theta)}{\partial t} \ne0$ with
$\beta=(1/\alpha'-\pi_0)/\pi_1$, for any\vspace*{1pt} interior point $(t, \theta, \alpha')$ in $[0,1]\times[0, \pi/2]\times[0, 1/\pi_0]$.
Then $\theta_0(\alpha')$ in criteria~(\ref{45}) is constant for all
$0 < \alpha' < 1/\pi_0$, if and only if
the solution $\theta$ of $t$ of the equation
%
\begin{equation}
\label{46} \frac{\partial F_1(t, \theta)}{\partial t} \Big/ \frac{\partial
F_1(t, \theta)}{\partial\theta}= \frac{\partial F_0(t, \theta)}{\partial t} \Big/
\frac{\partial
F_0(t, \theta)}{\partial\theta}
\end{equation}
is unique and equals a constant. Particularly, the above condition is
satisfied under assumptions \textup{(N1)} and
\textup{(N2)} with $\Sigma_0=\Sigma_1$.
\end{proposition}
Proposition~\ref{Proposition-2} implies
that $\theta_0(\alpha')$ does not depend on $\alpha'$ when $( \Phi
^{-1}(p_1),\break  \Phi^{-1}(p_2) )$
is bivariate normally distributed with identical covariance matrix
under the true null and nonnull.
For bivariate normal models with unequal covariance matrices,
Figure~\ref{Figure-2} shows that $\theta(\alpha')$ varies slightly
with $\alpha'$.
Numerical studies in Section~\ref{Sec-6} further confirm that $\theta
_0(\alpha')$ is robust to other bivariate distributions.
Hence, the selection of $\alpha'$ can be quite flexible except that
only mild restriction needs
to be imposed to make $\theta_0(\alpha')$ identifiable based on
conditions \textup{(C7)} to \textup{(C10)} in Appendix~\ref{appA}.
In particular, setting $\alpha'=\alpha/\pi_0$ will ensure that the
$\mathrm{Fdr}$ for various $\theta$ be controlled exactly at $\alpha$.

The Bayesian $\mathrm{Fdr}$ formula is equivalent to $F_1(t, \theta
)=\frac{1-\pi_0 \alpha' } {1- \pi_0 } F(t, \theta) $,\break  implying
that the criteria in (\ref{45})
can be replaced by $\theta_0( \alpha' )=\break \arg\max_{ 0 \le\theta\le
\pi/2} F( t^*_{\alpha' }(\theta), \theta)$. In Section~\ref{Sec-42},
we have two types
of estimators for $F_0(t, \theta)$, which can be used to develop
estimation approach for $\theta$.
Denoting $\widehat F _0 (t, \theta)$ to be either type of estimator,
the plug-in method for choosing the optimal direction $\theta_0(\alpha
')$ is thus given by
%
\begin{equation}
\label{47} \hat {\theta} \bigl(\alpha'\bigr)= \arg\max
_{ 0 \le\theta\le\pi/2 } \frac{R( \hat {t}^*_{\alpha' } (\theta), \theta)} m,
\end{equation}
where $\hat{t}^*_{\alpha' }(\theta)=
\sup\{ 0 \le t \le1\dvtx  { m \widehat F_0 ( t, \theta) } / { \{ R(t,
\theta) \vee1 \} } \le\alpha' \}$.
For notational clarity, we denote by $\{ \hat{t}{}_{\alpha'}^{ *
\I }(\theta), \hat { \theta}{}^{ \I}(\alpha') \}$ and $ \{
\hat {t}_{ \alpha' }^{ * \II }(\theta), \hat {\theta} {}^ {
\II}(\alpha') \}$
the estimators of $\{ t^*_{\alpha' }(\theta), \theta_0(\alpha') \}$
obtained by the parametric and nonparametric approaches, respectively.
%
%
\subsection{Procedures for estimating and controlling $\FDR$} \label{Sec-45}
For each fixed $\theta$, we provide two methods for $\FDR$ estimation
with respect to the projected \mbox{$p$-}values
$\{p_i(\theta)\dvtx  i=1,\ldots,m\}$ according to two estimators of
$F_0(t, \theta)$ proposed
in Section~\ref{Sec-42}.
\begin{longlist}
\item[\textit{Method} I.]
Incorporating the parametric approach for estimating $F_0 (t, \theta)$
and $\pi_0$
leads to a procedure for estimation and control of $\FDR$.
Combining (\ref{42}) and (\ref{44}), we propose
%
\begin{equation}
\label{48} \widehat{ \FDR}{}^{ \I}(t, \theta)= \frac{ \hat{\pi}_0^{ \I}( \theta) \widehat F_0^{ \I} (t,
\theta) }{
 \{ R(t, \theta)\vee1 \} /m }
\end{equation}
for our $\FDR$ estimation. A conservative $\FDR$ estimator naturally
leads to a procedure for controlling $\FDR$.
Similar to (\ref{23}), the data-driven threshold for the projected \mbox{$p$-}values
$\{p_i(\theta)\dvtx  i=1,\ldots,m\}$ is determined by
%
\begin{equation}
\label{49} t_{\alpha} \bigl( \widehat{ \FDR}{}^{\I} (:,\theta)
\bigr) = \sup\bigl\{ 0 \le t \le1\dvtx  \widehat{ \FDR}{}^{\I} (t,\theta) \le
\alpha\bigr\}.
\end{equation}
A null hypothesis is rejected if the corresponding $p(\theta)$ is less
than or equal to the threshold
$t_{\alpha} ( \widehat{ \FDR}{}^{\I} (:,\theta))$. The data-driven
threshold (\ref{49})
together with the point estimation method (\ref{48}) for the false
discovery rate
comprises the first $\FDR$ procedure, denoted by $\FDR^{\I}$.

\item[\textit{Method} II.]
The nonparametric approach proposed for estimating $F_0 (t, \theta)$
and $\pi_0$
can substitute the parametric counterpart in method~I. Similar to (\ref
{48}) and (\ref{49}),
the procedure for the estimation and control of $\FDR$ is given by
%
\begin{eqnarray}
\widehat{ \FDR}{}^{ \II}(t, \theta) &=& \frac{ \hat{\pi}_0^{ \II}( \theta) \widehat F_0^{ \II}
(t, \theta) }{
 \{ R(t, \theta)\vee1 \}/m }, \label{410}
\\
t_{\alpha} \bigl( \widehat{ \FDR}{}^{\II} (:,\theta)\bigr) &=&\sup
\bigl\{ 0 \le t \le1\dvtx  \widehat{ \FDR}{}^{\II} (t,\theta) \le \alpha\bigr\}.
\label{411}
\end{eqnarray}
The second $\FDR$ procedure, denoted by $\FDR^{\II}$, consists of
(\ref{410}) and (\ref{411}).
\end{longlist}
%
\begin{remark}
Incorporating $\hat{\theta} {}^{ \I} (\alpha' ) $ and $\hat {\theta} {}^{ \II} (\alpha' ) $
obtained from Section~\ref{Sec-44}
into $\FDR^{\I}$ and $\FDR^{\II}$, respectively, we obtain our final
procedure
for estimating and controlling $\FDR$.
\end{remark}
%
\subsection{Issue on stability and power for the $\SIM$ procedure} \label{Sec-46}
In this subsection, we first investigate the stability of the $\SIM$
procedure when the preliminary \mbox{$p$-}value is not accurate.
Suppose that the bivariate \mbox{$p$-}value $(p_1, p_2)$ is calculated from
the bivariate test statistic
$(X_1, X_2)$ with\vspace*{1pt} marginal true null $\CDF$s $F_{0; X_1}$ and $F_{0; X_2}$.
Due to some perturbation on $X_1$, w\vspace*{1pt}e observe a contaminated version
$\widetilde X_1$ with the true null $\CDF$ $F_{0; \widetilde X_1}$.
By using the incorrect true null $\CDF$ $F_{0; X_1}$,
the preliminary $p_1$ is incorrectly calculated as $\tilde p_1$. A
natural question is how
sensitive our $\SIM$ methods are if $X_1$ carries some wrong information.
%
\begin{proposition} \label{Proposition-3}
Suppose $(X_1,X_2)$ are the preliminary and primary test statistics for
one-sided hypotheses, where $F_{0; X_1}$ and $F_{0; X_2}$ are their
marginal $\CDF$s under the true null, respectively. Assume the
classical errors-in-va\-riables model on $X_1$, that is, $\widetilde
X_1=X_1+\eta$, where $\eta$ is independent of $(X_1,X_2)$ and the
p.d.f.s of $X_1$ under the true null and $\eta$ are both
symmetric with respect to~$0$. If the joint p.d.f. of $(p_1,
p_2)$ under the true null, where $(p_1,p_2)= (F_{0; X_1}(X_1),
F_{0; X_2}(X_2) )$ for left-sided hypotheses or $(p_1, p_2)=
(1-F_{0; X_1} (X_1),\break  1-F_{0; X_2} (X_2) )$ for right-sided
hypotheses, is centrally symmetric with respect to $(1/2,1/2)$, then
the joint p.d.f. of $(\tilde p_1, p_2)$ under the true
null is also centrally symmetric with respect to $(1/2,1/2)$, where
$\tilde p_1=F_{0; X_1}(\widetilde X_1)$ for left-sided hypotheses
or $\tilde p_1=1-F_{0; X_1}(\widetilde X_1)$ for right-sided hypotheses.
\end{proposition}

Proposition~\ref{Proposition-3} indicates that $\FDR$ of method~$\II$
can still be controlled even if the preliminary test statistic is
measured with classical additive error \cite{Carrolletal2010}.
Although our discussion is restricted to the situation where the
p.d.f. of preliminary test statistic under the true null is
symmetric about 0, it indeed includes a large class of distributions,
for example, normal distribution and $t$ distribution.
In general, it can be verified that method~$\II$ is valid if
%
\begin{equation}
\label{412} f_{0; (\tilde p_1, p_2)}(\tilde p_1, p_2)
\le f_{0;
(\tilde p_1, p_2) }(1-\tilde p_1, 1-p_2),
\end{equation}
where $f_{0; (\tilde p_1, p_2)}(\tilde p_1, p_2)$ is the
p.d.f. of $(\tilde p_1, p_2)$ under the true null and
$\tilde p_1+p_2 \le1$. Under (\ref{412}), the probability mass
under the true null in the upper-right tail of $(\tilde p_1, p_2)$
is no less than that in the lower-left tail, resulting in some
conservative procedure. To simplify the argument, we only consider the
case where $\tilde p_1$ and $p_2$ are independent, which
simplifies the sufficient condition (\ref{412}) to
%
\begin{equation}
\label{413} f_{0; \tilde p_1 }(\tilde p_1 ) \le
f_{0; \tilde p_1}(1-\tilde p_1 ),
\end{equation}
where $f_{0; \tilde p_1}( \tilde p_1 )$ is the
p.d.f. of $\tilde p_1$ under the true null
and $0 \le \tilde p_1 \le1/2$. Some pairs of asymmetric
distributions of ${X_1}$ and
${\widetilde X_1}$, satisfying the condition (\ref{413}), are
summarized below:
\begin{itemize}
\item
${X_1} \sim\Exp( \bar {\lambda}_1)$ and ${\widetilde X_1}
\sim\Exp( \bar {\lambda}_2)$ with $\bar {\lambda
}_1>\bar {\lambda}_2>0$, where
$\Exp({\lambda} )$ denotes the exponential distribution with
parameter ${\lambda}$.

\item
${X_1} \sim\chi^2_{r}$ and ${\widetilde X_1} \sim\chi^2_{s}$ with $r<s$.

\item\emph{Chi-square versus weighted chi-square distribution} \\
${X_1} \sim\chi^2_{r}$ and ${\widetilde X_1} \sim\sum_{i=1}^{r}
\omega_i Z_i^2$,
where $\{Z_i\}_{i=1}^{r} \stackrel{\mathrm{i.i.d.}}{\sim} N(0,1)$
and $\omega_i \ge1$, $i=1,\ldots,r$.

\item \emph{$F$ versus generalized $F$ distribution} \\
${X_1} \sim F(r,s)$ and ${\widetilde X_1} \sim\frac{ (\sum_{i=1}^{r} \omega_i Z_i^2 ) /r }{\chi^2_{s}/s}$,
where $\{Z_i\}_{i=1}^{r} \stackrel{\mathrm{i.i.d.}}{\sim} N(0,1)$,
$\sum_{i=1}^{r} \omega_i Z_i^2 $ is independent of $\chi^2_{s}$, and
$\omega_i \ge1$, $i=1,\ldots,r$.
\end{itemize}

Having established that the $\SIM$ procedure controls $\FDR$ when the
preliminary \mbox{$p$-}values carry some wrong information, we next turn to
theoretically justify why the current way of combination of the
bivariate \mbox{$p$-}value achieves a higher power.
Let $t_{\alpha}(\theta)$ denotes the threshold such that
${\pi_0 F_0(t, \theta)}/\break { F(t, \theta)}=\alpha$.
Then the power function can be formulated by $F_1(t_{\alpha}(\theta
),\theta)=\beta'F_0(t_{\alpha}(\theta),\theta)$, with $\beta
'=(1/\alpha-1)\pi_0/\pi_1$.
Our goal is to quantify how much power can be improved via combining
the bivariate \mbox{$p$-}value.
From the Bayesian $\mathrm{Fdr}$ formula, the ratio of power of the
$\SIM$ procedure to conventional
multiple testing procedure using $p_2$ alone ($\theta=\pi/2$) can be
derived as
\begin{eqnarray*}
\frac{F_{1}(t_{\alpha}(\theta), \theta) }{ F_1(t_{\alpha}(\pi/2),
\pi/2 )} &=& \frac{F_{0}(t_{\alpha}(\theta), \theta) }{ F_0(t_{\alpha}(\pi
/2), \pi/2 )}
\\
&=& 1+\frac{ ({\partial}/({\partial\theta}))  \{
F_{0}(t_{\alpha}(\theta), \theta)  \} |_{\theta={\pi/2}}
(\theta-\pi/2) }{
 F_0(t_{\alpha}(\pi/2), \pi/2 ) }
 \\
&&{}  +O \bigl(
(\theta-\pi/2)^2 \bigr)
\\
&=& 1+\Delta(\theta)+O \bigl( (\theta-
\pi/2)^2 \bigr).
\end{eqnarray*}
More derivations in Appendix~\ref{appB} yield that the ratio of power improved
when $\theta$ is close to $\pi/2$ is approximated by
%
\begin{eqnarray}\label{414}
1+\Delta(\theta) &=& 1+\frac{ \phi[ \Phi^{-1}\{ t_{\alpha}(\pi/2) \} ] f_{1; p_2}(
t_{\alpha}(\pi/2) ) f_{0; p_2}( t_{\alpha}(\pi/2) ) }{
 f_{1; p_2}(t_{\alpha}(\pi/2) )-\beta' f_{0; p_2}(t_{\alpha}(\pi
/2) ) }
\nonumber\\[-8pt]\\[-8pt]
&&\hspace*{19pt} \times \frac{ (\theta-\pi/2) }{ F_0(t_{\alpha}(\pi/2), \pi/2 )}
\times \I(p_1),\nonumber
\end{eqnarray}
where $\I(p_1)=E_{H_0} \{ \Phi^{-1}(p_1) \mid p_2=t_{\alpha}(\pi/2)
\} -E_{H_1} \{ \Phi^{-1}(p_1) \mid p_2=t_{\alpha}(\pi/2) \} $,
$f_{0; p_2}$ and $f_{1; p_2}$ are the p.d.f.s of $p_2$ under
true null and nonnull, respectively. If the alternative distribution of
$p_2$ is strictly concave, similar argument in \cite
{GenoveseWasserman2002} yields that $ f_{1; p_2}(t_{\alpha}(\pi/2)
)-\beta' f_{0; p_2}(t_{\alpha}(\pi/2) )<0 $.
The term $\I(p_1)$ is positive, provided that the preliminary
\mbox{$p$-}values have some potential to detect the power. Combining these, we
have $1+\Delta(\theta)>1$.

Under assumptions \textup{(N1)} and \textup{(N2)}, $\I(p_1)$ has an
explicit form
%
\begin{eqnarray}\label{415}
&& \bigl[\mu_{0;1}+\rho_0
\sigma_{0;1}/\sigma_{0;2}\bigl\{\Phi^{-1}
\bigl(t_{\alpha
}(\pi/2)\bigr) -\mu_{0;2} \bigr\}\bigr]
\nonumber\\[-8pt]\\[-8pt]
&&\qquad {} - \bigl[
\mu_{1;1}+\rho_1 \sigma_{1;1}/\sigma_{1;2}
\bigl\{\Phi^{-1}\bigl(t_{\alpha
}(\pi/2)\bigr) -\mu_{1;2}
\bigr\}\bigr].\nonumber
\end{eqnarray}
From (\ref{415}), the correlation ($\rho_0$) between components of
the bivariate \mbox{$p$-}value under the true null and that ($\rho_1$) under
nonnull play different roles in improving power. The $\SIM$ procedure
using prior information and primary \mbox{$p$-}values that are negatively
correlated under the null hypothesis but positively correlated under
the alternative is a general approach that can substantially increase
power in practice.
%
\section{Asymptotic justification} \label{Sec-5}
In many applications such as biology, med\-icine, genetics, neuroscience,
economics and finance,
tens of thousands of hypotheses are tested simultaneously. It is hence
natural to investigate the behavior
of the two approaches we proposed for the large number $m$ of
hypotheses. In this section, we focus on the asymptotic
properties of the nonparametric estimator, $\widehat {\FDR}{}^{\II}(t,
\hat {\theta}{}^{\II} (\alpha' ) )$. All
theorems presented in this section can be derived similarly for the\vspace*{2pt}
parametric approach as long as the bivariate normality
for $(\Phi^{-1}(p_1), \Phi^{-1}(p_2))$ is satisfied.\vspace*{2pt}

Theorem~\ref{Theorem-1} below establishes the consistency of
$\hat {\theta}{}^{\II}(\alpha') $.
Intuitively, $\hat {\theta}{}^{\II}(\alpha') $ is analogous to an
$\mathrm{M}$-estimator such as least-squares estimators
and many maximum-likelihood estimators. However, typical proof of consistency
of $\mathrm{M}$-estimators is not applicable to $\hat {\theta
}(\alpha')$ because the $\CDF$
involved in (\ref{47}) is not differentiable. Hence,
the theoretical derivation is nontrivial and challenging. We will
provide Lemmas~\ref{Lemma-1}--\ref{Lemma-3} in Appendix~\ref{appA},
which are necessary for proving Theorem~\ref{Theorem-1}.
%
\begin{theorem} \label{Theorem-1}
Assume conditions \textup{(C1}) to \textup{(C9)} in Appendix~\ref{appA}. Then
$\hat{\theta}{}^{ \II} (\alpha' ) $ converges to $\theta
_0(\alpha' )$ almost surely.
\end{theorem}

Theorem~\ref{Theorem-2} below reveals that the proposed estimator
$\widehat {\FDR}{}^{\II}$
not only controls the $\FDR$ simultaneously for all $t\ge\delta$ and
$\delta>0$ for fixed $\theta$, but also
provides simultaneous and conservative control when incorporating the
data-driven estimator $\hat {\theta}{}^{\II}(\alpha' ) $.
%
\begin{theorem} \label{Theorem-2}
Assume conditions \textup{(C1}) to \textup{(C10)} in Appendix~\ref{appA}. Then
$\widehat { \FDR}{}^{ \II } ( t, \hat {\theta}{}^{ \II}(\alpha')
) $ provides
simultaneously conservative control of $\FDR(t,\break  \theta_0(\alpha') )$
in the sense that
\begin{eqnarray*}
\lim_{m \to \infty} \inf_{t \ge\delta} \bigl\{ \widehat {
\FDR}{}^{ \II} \bigl(t, \hat {\theta}{}^{ \II}\bigl(
\alpha'\bigr) \bigr) - \FDR \bigl(t, \theta_0\bigl(
\alpha'\bigr) \bigr)\bigr\} &\ge& 0,
\\
\lim_{m \to \infty}
\inf_{t \ge\delta} \biggl\{ \widehat { \FDR }{}^{ \II } \bigl(t,
\hat {\theta}{}^{ \II}\bigl(\alpha'\bigr) \bigr) -
\frac{ V(t, \theta_0(\alpha') ) }{ R( t, \theta_0(\alpha') )\vee
1} \biggr\} &\ge& 0
\end{eqnarray*}
with probability $1$.
\end{theorem}

To show that the proposed estimator $\widehat {\FDR}{}^{\II}(t,
\hat {\theta}{}^{ \II}(\alpha') )$
provides strong control of
$\FDR(t, \theta_0(\alpha') )$ asymptotically, we define
\begin{eqnarray*}
&&\widehat {\FDR}_{\lambda}^{\infty}(t, \theta)=\frac{ \{\pi_0+
\pi_1 (({1-F_1(\lambda, \theta)})/({1-F_0(\lambda, \theta)}))\}
F_0(t, \theta) }{ F(t, \theta) },
\end{eqnarray*}
which is a pointwise limit of $\widehat {\FDR}_{\lambda}^{ \II
}(t,\theta)=
\frac{ \hat {\pi}_0^{\II}(\lambda, \theta) \widehat F_0^{\II
}(t, \theta) }{ \{ R(t,\theta)\vee1 \} /m }$
under conditions \textup{(C1}) and \textup{(C2)} and Lemma~\ref{Lemma-2} in Appendix~\ref{appA}. The notations $\hat {\pi}_0^{ \I
}(\lambda, \theta)$
and $\hat {\pi}_0^{\II}(\lambda, \theta)$
are defined in a way similar to those in the algorithm of
Section~\ref{Sec-43}.

\begin{theorem} \label{Theorem-3}
Assume conditions \textup{(C1}) to \textup{(C10)} in Appendix~\ref{appA}.
Also, suppose that the sequence of values $\{\lambda_j\dvtx  j=1,\ldots,n\} \in(0,1/2]^n$ and
$n$ is a fixed finite integer. If for each $\lambda_j$, there is
$t_j \in(0,1]$ such that $\widehat { \FDR}{}^{\infty}_{\lambda_j}
(t_j,\break  \theta_0(\alpha') ) < \alpha$, then
\begin{eqnarray*}
&&\limsup_{m \to\infty} \FDR \bigl( t_{\alpha'} \bigl( \widehat {
\FDR}{}^{ \II}\bigl(:, \hat {\theta }{}^{ \II}\bigl(
\alpha'\bigr) \bigr) \bigr), \hat {\theta}{}^{ \II}\bigl(
\alpha' \bigr) \bigr) \le\alpha.
\end{eqnarray*}
\end{theorem}
%
\section{Numerical studies} \label{Sec-6}
In this section, we carry out simulation studies to evaluate the
performance of the $\SIM$ procedure in the aspects of controlling
$\FDR$ and detection power,
using the two proposed methods under various bivariate models for the
preliminary and primary test statistics.
The sequence of values $\{\lambda_j\dvtx  j=1, \ldots, n\}$ in the algorithm of Section~\ref{43} is
$\{0.02,  0.04,  0.06,  0.08,  0.1,  0.125,  0.15,  \ldots,
0.5\}$.
For simplicity, the constant $c$ in Section~\ref{Sec-42} is set to be
0. Unless otherwise stated,
$\alpha'$ is simply set to be $\alpha$ throughout this section,
following Proposition~\ref{Proposition-2}.
All simulations are based on $500$ replications.

The following procedures are compared:
\begin{itemize}
\item
Conventional $\FDR$ procedure: the $\FDR$ method using (\ref{22})
and (\ref{23}) with $\pi_0$ dynamically selected by the algorithm in Section~\ref{Sec-43}.
\item
Weighted multiple testing procedure: the weighted multiple testing
procedure proposed by \cite{Genoveseetal2006}, where the weighting
scheme is determined automatically by the preliminary \mbox{$p$-}values; refer
to the cumulative weights with $B=2$ in Section~\ref{Sec-33} for detail.
\item
Two-stage multiple testing procedure: the two-stage procedure defined
by \cite{Bourgonetal2010} with the first stage being preliminary
\mbox{$p$-}values filtering. The proportion of hypotheses to be removed in the
\emph{filtering} stage is set to be 50\%.
\end{itemize}

Note that the ``50\% variance filter'' in \cite{Bourgonetal2010} shares
the same spirit as the ``two-stage multiple testing procedure'' except that
the overall sample variance serves as the \emph{filter} statistic.

%
\subsection{Example 1: Bivariate normal model} \label{Sec-61}
This example comes from hypothesis testing of mean shift in normal
models, that is, $X=\mu+\varepsilon$
with $\varepsilon\sim N(0,1)$.
We perform $m={}$10,000 independent right-sided hypotheses testing for
$H_0\dvtx  \mu=0$
versus $H_1\dvtx  \mu>0$.
Among all the null hypotheses, a proportion $\pi_0$ of them are from
the true null hypotheses.
For the $i$th test, we generate a bivariate test statistic $(x_{i1},
x_{i2})$ from a bivariate normal distribution $\mathcal{N}(\bolds{\mu}, \Sigma)$ where $\Sigma=(\sigma_{ij})_{2 \times2}$
with $\sigma_{11}=1$, $\sigma_{12}=\sigma_{21}=\rho$ and $\sigma_{22}=1$.
We set $\bolds{\mu}=\mathbf{0}$ under the true null and
$\bolds{\mu}=(\mu_1,\mu_2)$
under nonnull. The marginal \mbox{$p$-}values for the $i$th test are
$p_{i1}=1-\Phi(x_{i1})$ and
$p_{i2}=1-\Phi(x_{i2})$, for $i=1,\ldots,m$.

\begin{figure}

\includegraphics{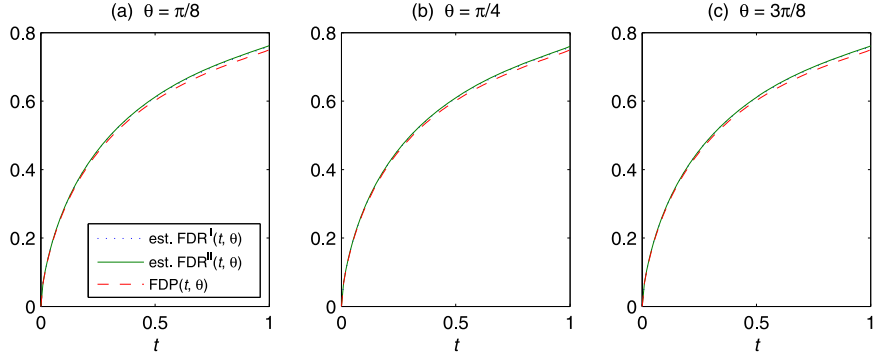}

\caption{Estimated $\FDR$ for methods $\I$~and~$\II$
and the corresponding true $\FDP$ as a function of threshold $t$ and
$\theta$ in Example $1$. Here, $\mu_1=\mu_2=2 $, $\pi_0=0.75$ and
$\rho=0.2$.}\label{Figure-3}
\end{figure}

To evaluate the overall performance of the estimated $\FDR(t,\theta)$
of methods~I~and~II
at the same threshold $t \in[0,1]$, we consider the
scenario where $\mu_1=\mu_2=2$, $\pi_0=0.75$ and $\rho=0.2$.
For notational convenience, denote by
$\FDP (t,\theta)=V(t,\theta) / \{ R(t,\theta) \vee1\}$
the false discovery proportion at threshold $t$ with respect to $\{
p_i(\theta)\dvtx  i=1, \ldots,m \}$.
Figure~\ref{Figure-3} compares\vspace*{2pt}
the average values of $\widehat{ \FDR}{}^{\I}(t,\theta)$,
$\widehat{ \FDR}{}^{\II}(t,\theta)$ and $\FDP(t,\theta)$
for $\theta=\pi/8$, $\pi/4$, $3\pi/8$. For each case, these two
types of estimators are very close to true $\FDP$,
lending support to the parametric and nonparametric estimation
procedures in Section~\ref{Sec-4}.

%

%
\begin{figure}[b] 

\includegraphics{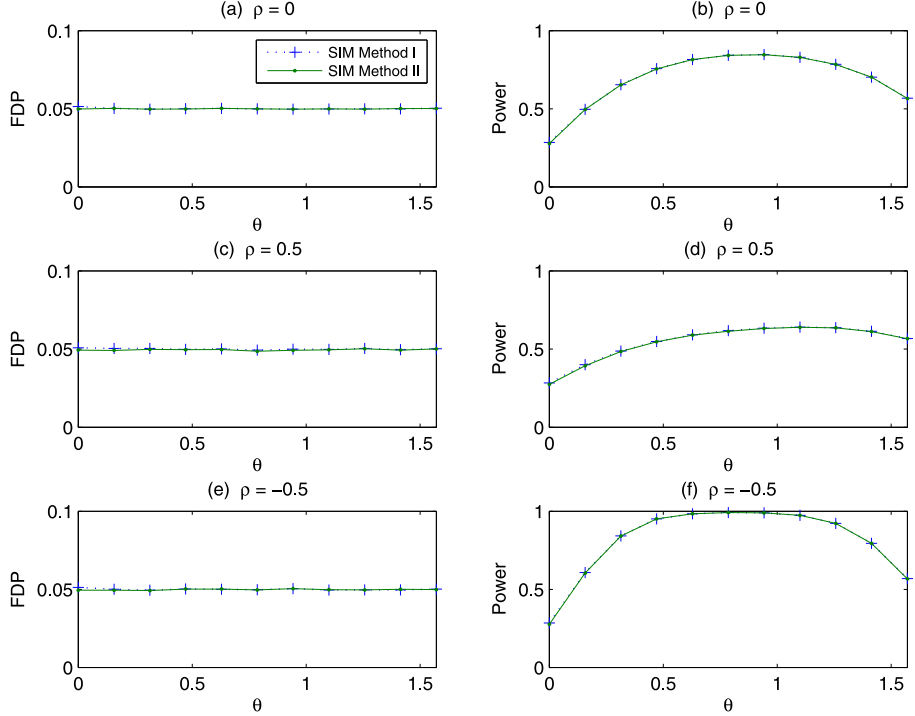}

\caption{Calculated $\FDP$ and power as a function of
$\theta$ and $\rho$ in Example $1$.
Here, $\mu_1=2$, $\mu_2=2.5$, $\pi_0=0.75$ and $\alpha=0.05$.}\label{Figure-4}
\end{figure}

To illustrate the role of $\theta$ for detecting power in our proposed
procedure,
a sequence of values $\{ \theta_l=(l-1)/10 \times \pi/2\dvtx  l=1,\ldots,11\}$ are designed. For simplicity, we
consider the scenario where $\mu_1=2, \mu_2=2.5$, $\pi_0=0.75$,
$\alpha=0.05$ and $\rho=\{0, 0.5, -0.5 \}$.
Figure~\ref{Figure-4} corresponds to the calculated $\FDP$ [i.e.,
$\FDP(\hat {t}_{\alpha})$] and the calculated power [i.e.,
$S(\hat {t}_{\alpha})/m_1$]
as a function of $\theta$, for $\rho=0, 0.5, -0.5$, respectively.
In either case, we observe that the average values of the calculated
$\FDP$ for both $\FDR^{\I}$ and $\FDR^{\II}$ are almost controlled
at $\alpha=0.05$ for
all $\theta$, and by appropriately choosing $\theta$, the $\SIM$
methods outperform
the conventional $\FDR$ procedure using $p_2$ alone (with $\theta=\pi/2$).
The correlation between the components of the bivariate \mbox{$p$-}value
sensitively affects the optimal power. Negative correlation
distinguishes $p_1$ and $p_2$ most significantly,
thus it is expected that this case can improve the power most via
combining the bivariate \mbox{$p$-}value.
As a comparison, positive correlation diminishes the detection
slightly. However, the power is still improved significantly
when comparing to the conventional $\FDR$ procedure using $p_2$ alone.

\begin{figure}

\includegraphics{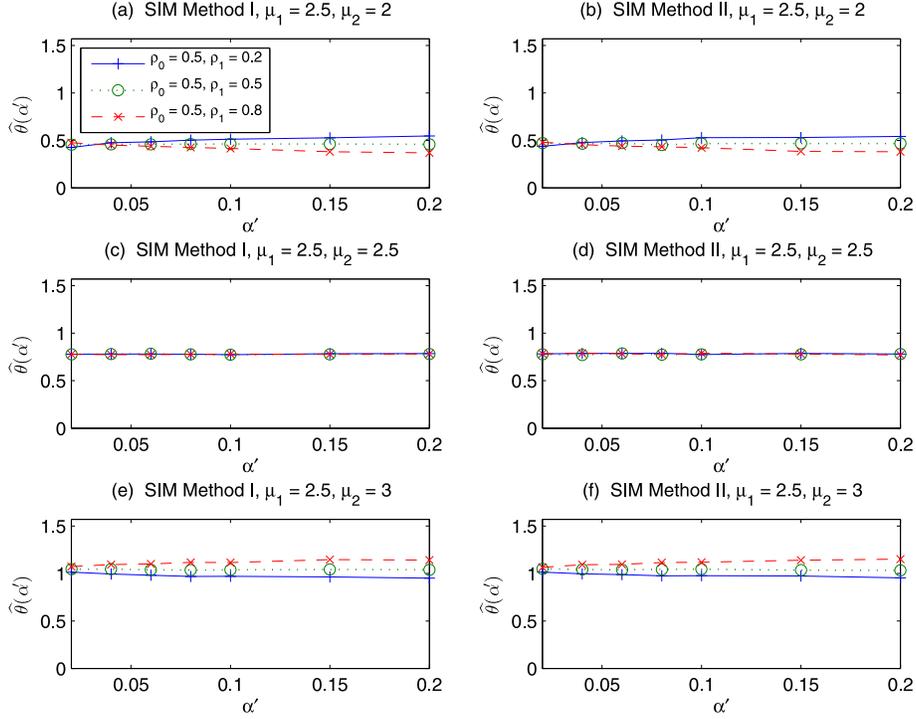}

\caption{$\hat {\theta}(\alpha')$ as a function of
$\alpha'$ for various combination of $(\mu_1, \mu_2, \rho_0, \rho
_1)$ in bivariate
normal models for $ (x_{i1}, x_{i2} )$. Here, $\pi_0=0.75$.}\label{Figure-5}
\end{figure}

To confirm the consistency of $\hat{\theta}(\alpha')$, we
compare 10 scenarios, where $\pi_0=0.75$,
$\rho=0.2$, $\alpha'=\{0.05, 0.1\}$ and $(\mu_1, \mu_2)$ takes five
different pair-values.
From Proposition~\ref{Proposition-2}, the optimal value $\theta
_0(\alpha')$ is constant for
different $\alpha'$, denoted by~$\theta_0$. Table~\ref{Table-2}
compares the average value of $\hat \theta(\alpha')$
and its standard error of methods~I~and~II with the optimal value
$\theta_0$.
In all situations, estimators are very close to the optimal value
$\theta_0$ except that the standard error
of $\hat{\theta}(\alpha')$ by method~II is slightly larger than
that by method~I.
This phenomenon is not surprising,\vspace*{1pt} since the nonparametric fit for
$F_0(t, \theta)$ and $F(t, \theta)$
contaminates the estimator $\hat{\theta}{}^{\II}(\alpha')$.
For unequal covariance matrices in bivariate normal models for $
(x_{i1}, x_{i2} )$ with the correlation
coefficients $\rho_0$ and $\rho_1$ in the\vspace*{1pt} true null and nonnull, respectively,
Figure~\ref{Figure-5} shows the stability of $ \hat {\theta
}(\alpha')$ for various choices of $\alpha'$
using both methods~I~and~II.

%
\begin{table}[t] 
\tabcolsep=0pt
\caption{Mean and standard error of $\hat {\theta}(\alpha')$ by $\FDR^{\I}$ and $\FDR^{\II}$
for $10$ scenarios, where $\pi_0=0.75$, $\rho=0.2$, $\alpha'=\{0.05,
0.10\}$ and $(\mu_1,\mu_2)$ are set to be $(2, 1)$, $(2, 1.5)$, $(2, 2)$,
$(2, 2.5)$, $(2, 3)$, respectively}\label{Table-2}
\begin{tabular*}{\tablewidth}{@{\extracolsep{\fill}}@{}lccccc@{}}
\hline
& \multicolumn{2}{c}{$\bolds{\alpha'=0.05}$} & \multicolumn{2}{c}{$\bolds{\alpha'=0.10}$}\\[-6pt]
& \multicolumn{2}{c}{\hrulefill} & \multicolumn{2}{c}{\hrulefill}
\\
$\bolds{(\mu_1, \mu_2)}$ & $\bolds{\FDR^{\I}}$& $\bolds{\FDR^{\II}}$& $\bolds{\FDR^{\I}}$&$\bolds{\FDR^{\II}}$& $\bolds{\theta_0}$\\
\hline
 {(2, 1)}  &  0.3231 (0.07)  &  0.3273 (0.13)  & 0.3190 (0.07)  &  0.3192 (0.10)  &  0.3218 \\
 {(2, 1.5)}  &  0.5706 (0.07)  &  0.5713 (0.12)  & 0.5684 (0.07)  &  0.5687 (0.10)  &  0.5743 \\
 {(2, 2)}  &  0.7785 (0.06)  &  0.7828 (0.11)  & 0.7813 (0.07)  &  0.7808 (0.09)  &  0.7854 \\
 {(2, 2.5)}  &  0.9523 (0.06)  &  0.9525 (0.09)  & 0.9436 (0.07)  &  0.9490 (0.09)  &  0.9505 \\
 {(2, 3)}  &  1.0732 (0.06)  &  1.0734 (0.09)  & 1.0720 (0.08)  &  1.0755 (0.11)  &  1.0769 \\
 \hline
\end{tabular*}
\end{table}

%

\begin{figure}

\includegraphics{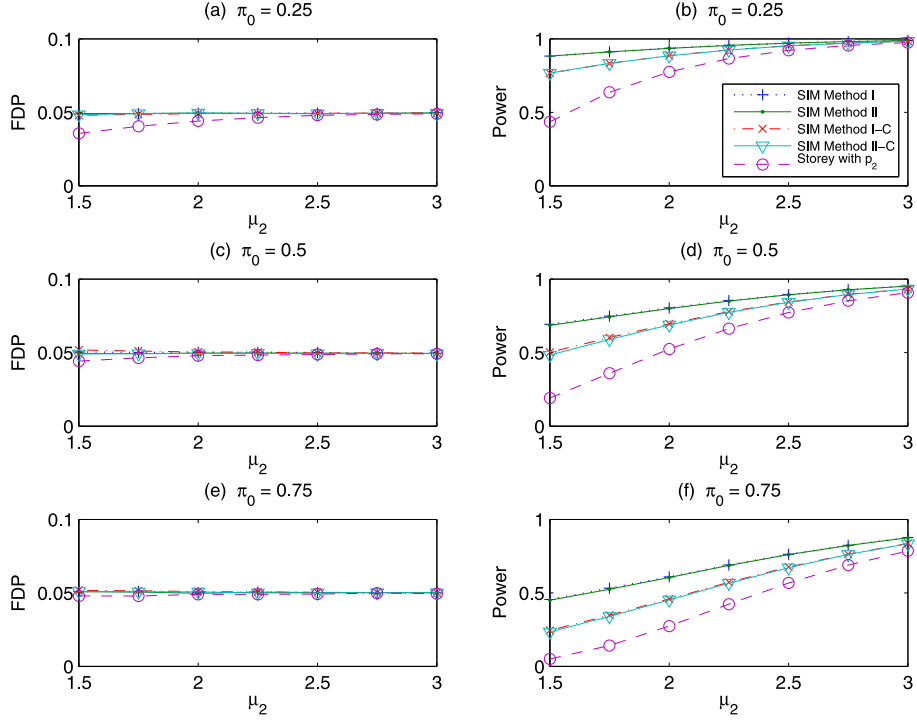}

\caption{Calculated $\FDP$ and power as a function of $\mu
_2$ and $\pi_0$
for the $\SIM$ methods,
their contaminated versions $(\SIM$ method~$\I$-\textup{C}, $\SIM$ method
$\II$-\textup{C}) and the conventional
$\FDR$ procedure (storey with $p_2)$ in Example $1$. The
contamination scheme is $\widetilde X_1=X_1+\eta$, where $\widetilde
X_1$ is the observable preliminary test statistic and $\eta$ is a
standard normal noise independent of the unobservable one $X_1$. Here,
$\mu_1=2$, $\alpha=0.05$ and $\rho=0.2$.}
\label{Figure-6}
\end{figure}

In the previous simulation results, we have demonstrated that for a
fixed value of $\theta$,
$\widehat {\FDR}{}^{\I}$ and $\widehat {\FDR}{}^{\II}$ provide
simultaneous and conservative control of $\FDR$; and that power can
improve significantly by appropriately choosing $\theta$.
Does the conclusion continue to hold for random $\hat{ \theta
}(\alpha)$? Figure~\ref{Figure-6}
examines the control of $\FDR$ as well as power comparison of the
$\SIM$ methods, their corresponding contaminated
versions and the conventional $\FDR$ procedure for various
combinations of $(\mu_2, \pi_0)$.
The left panels of Figure~\ref{Figure-6} compare the calculated $\FDP
$ of all settings.
Clearly, the calculated $\FDP$ for the $\SIM$ methods and their
contaminated versions is controlled at the prespecified $\alpha=0.05$,
confirming that the $\SIM$ methods are still valid when the
preliminary test statistics carry some wrong information.
The right panels correspond to the power of all the approaches.
We\vspace*{2pt} observe that the average values of power
of $\FDR^{{\I } }(t, \hat {\theta}(\alpha) )$ and $\FDR
^{{\II } }(t, \hat{\theta}(\alpha) )$
are consistently higher
than that of the conventional $\FDR$ procedure using $p_2$ alone.
Remarkably, the power of the contaminated versions of the $\SIM$
methods is not adversely affected, but between that of the $\SIM$
methods and the conventional $\FDR$ procedure.
%

%
\begin{figure}

\includegraphics{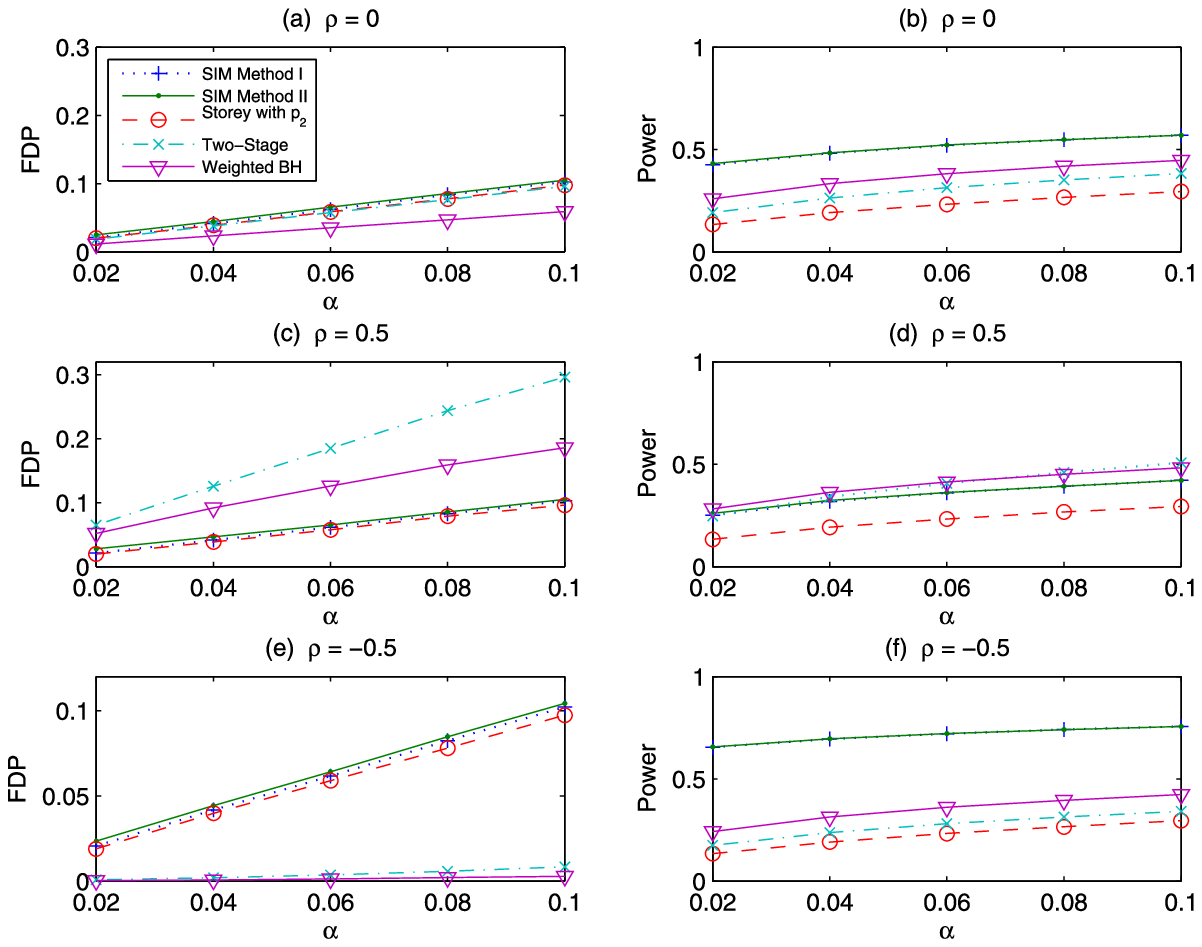}

\caption{Calculated $\FDP$ and power as a function of
$\alpha$ and $\rho$ for
the $\SIM$ methods, the weighed multiple testing procedure (weighted
$\BH)$, the two-stage multiple testing procedure
(two-stage) and the conventional $\FDR$ procedure (storey with $p_2)$ in Example $1$.
Here, the nonnull is a mixture of three bivariate normal distributions
with signals $(1,1)$, $(2,2)$ and $(3,3)$, respectively, and $\pi_0=0.9$.}\label{Figure-7}
\end{figure}

To further illustrate the advantage of the $\SIM$ methods, Figure~\ref{Figure-7} compares them with the weighted multiple testing procedure
and the two-stage multiple testing procedure which virtually use the
same amount of information from preliminary \mbox{$p$-}values and primary \mbox{$p$-}values
for various levels $\alpha$ and $\rho$ when the nonnull is a mixture
of three bivariate normal distributions with small, moderate and strong signals.
When the preliminary \mbox{$p$-}value and primary \mbox{$p$-}value are independent,
all the approaches are valid but the $\SIM$ methods outperform the
weighted multiple testing procedure and the two-stage multiple testing
procedure for all significant levels $\alpha$. Note that both the
weighted multiple testing procedure and the two-stage multiple testing
procedure are out of control if the components of bivariate \mbox{$p$-}value
are positively correlated and much more conservative under negative
dependence. In contrast, the $\SIM$ methods consistently estimate the
$\FDR$ under any dependence structure between the components of
bivariate \mbox{$p$-}value, providing much flexibility to choose \emph
{filters} or \emph{weights} in practice.

%
\subsection{Example 2: Bivariate $t$ distribution} \label{Sec-62}
In this example, we consider a set-up similar to Example 1 except that
the datasets are generated
from a bivariate $t$ distribution. To be specific, $\{ ( x_{i1}, x_{i2}
)\dvtx  i=1,\ldots,m\}$,
are sampled independently from a bivariate $t$ distribution with $3$
degrees of freedom and covariance matrix
identical to that in Example 1. Among all the null hypotheses, a
proportion $\pi_0$ of them
come from the true null hypotheses with mean zero, while the rest are
coming from nonnull
hypotheses with mean vector $\bmu=(\mu_1, \mu_2 )$.

\begin{figure}

\includegraphics{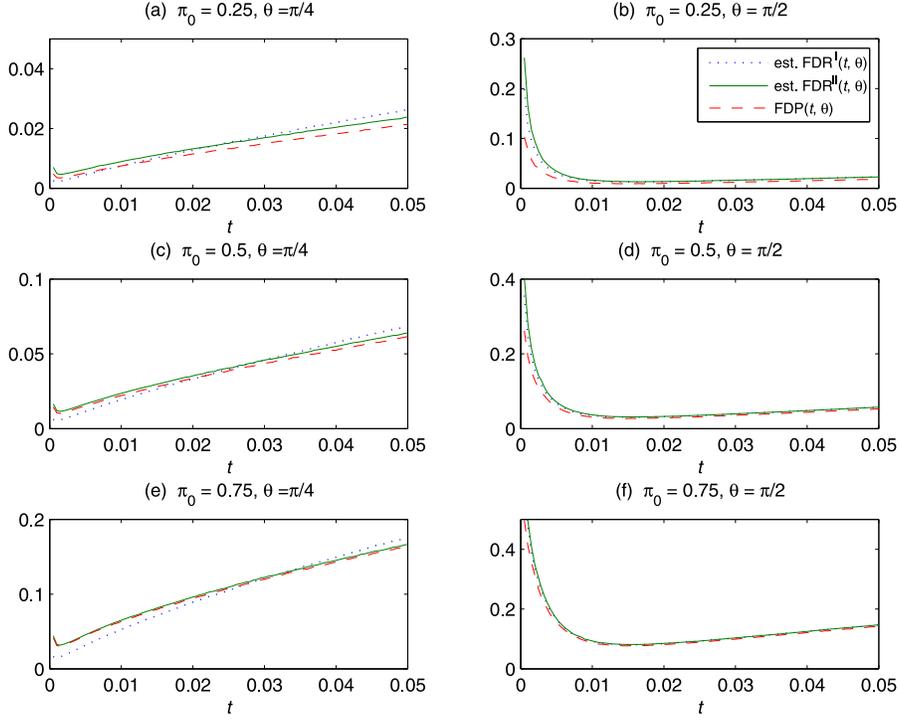}

\caption{Estimated $\FDR$ for methods~$\I$~and~$\II$ and
the corresponding true $\FDP$ as a function of $t$ in Example $2$.
Here, $\mu_1=\mu_2=4$ and $\rho=0.2$.}\label{Figure-8}
\end{figure}

Figure~\ref{Figure-8} compares\vspace*{2pt} the average values of the true $\FDP$,
$\widehat {\FDR} {}^{ \I }(t, \theta)$ and $\widehat {\FDR}{}^{\II}
(t, \theta)$ in a zoomed-in region of
$t \in[0, 0.05] $ for different combinations of $(\pi_0, \theta)$.
On the right panels where $\theta=\pi/2$ (using $p_2$), both methods~I~and~II provide conservative estimates of $\FDR$. For the case
$\theta=\pi/4$ on the left panels,
method~II provides conservative estimation of $\FDR$ and is less
conservative as
$\pi_0$ increases. Unlike method~$\II$, method~$\I$ underestimates
the true $\FDR$ for small $t$ and overestimates
it for large $t$, which makes the $\FDR$ out of control for small
$\alpha$.
This is not surprising, since the bivariate $t$ distribution with very
low degrees of freedom
violates the normality assumption.

%
\begin{figure}[b] 

\includegraphics{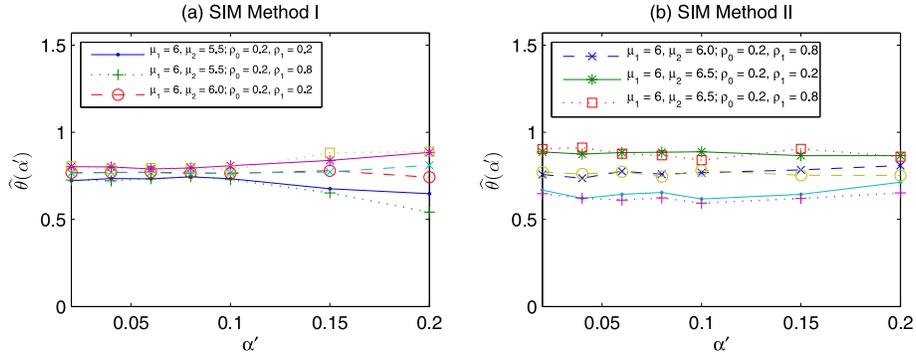}

\caption{$\hat {\theta}(\alpha')$ versus $\alpha'$
for various combination of
$(\mu_1, \mu_2, \rho_0, \rho_1)$ for bivariate $t$ distributions.
Here, $\pi_0=0.75$.}\label{Figure-9}
\end{figure}
%
\begin{figure}

\includegraphics{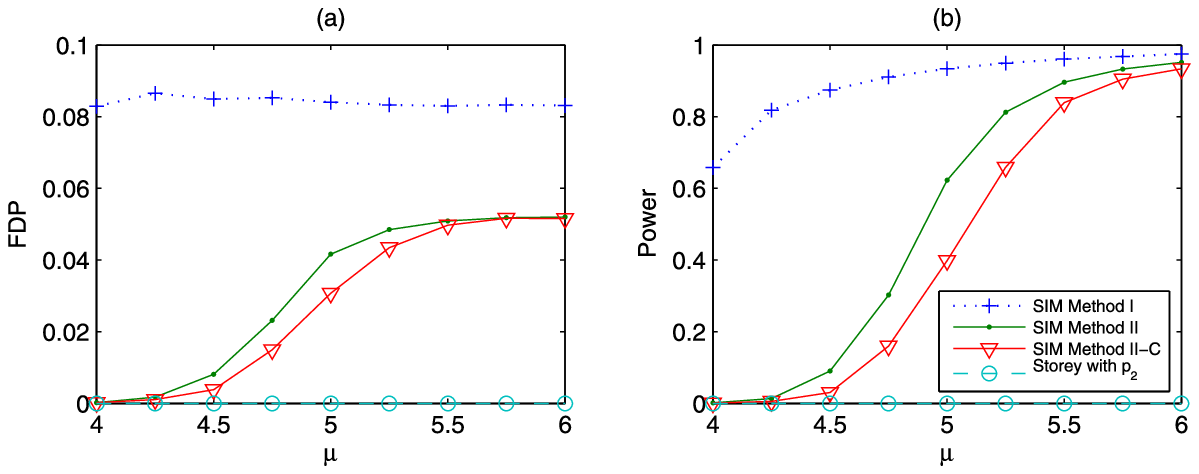}

\caption{Calculated $\FDP$ and power as a function of $\mu
$ $(\mu=\mu_1=\mu_2)$ for the $\SIM$ methods, the contaminated
version of method~$\II$ $(\SIM$ method~$\II$-\textup{C}) and the
conventional $\FDR$ procedure (storey with $p_2)$ in Example $2$.
The contamination scheme is $\widetilde X_1=X_1+\eta$, where
$\widetilde X_1$ is the observable preliminary test statistic and
$\eta$ is a standard normal noise independent of the unobservable one
$X_1$ in Example $2$.
Here, $\alpha=0.05$, $\rho=0.2 $, $\pi_0=0.9$ and $\mathrm{df}=3$.}\label{Figure-10}
\end{figure}

Before assessing the performance of the $\SIM$ methods incorporating
random $\hat {\theta}(\alpha)$, we first demonstrate that
$\hat {\theta}(\alpha')$
is robust to $\alpha'$ for various bivariate $t$ distributions in
Figure~\ref{Figure-9}, which lends support to
setting $\alpha'=\alpha$ when choosing the optimal projection direction.
Based on this setting, Figure~\ref{Figure-10}
summarizes the average values of the calculated $\FDP$ and power
of the $\SIM$ methods, the contaminated version of
method~$\II$ and the conventional $\FDR$ procedure for various
combinations of $(\mu_1, \mu_2)$.
We observe that the conventional $\FDR$ procedure lacks the ability to
detect statistical significance
for various signals even when $\alpha=0.05$. Nonetheless, by
incorporating the prior information from $p_1$ into $p_2$, method~$\II
$ improves the power while controlling the $\FDR$. Similar to the
previous case (Figure~\ref{Figure-6}), the calculated $\FDP$ for the
contaminated version of method~$\II$ is controlled at $\alpha=0.05$
and the corresponding power
is very close to that of method~$\II$. This illustrates the stability
of method~$\II$
when the preliminary \mbox{$p$-}value is not accurate.
Note that, even if method~$\I$ appears more powerful than method~$\II
$, the calculated $\FDP$ for method~$\I$ is out of control at level
higher than $\alpha=0.05$. The uncontrolled performance of method~I
indicates that the nonparametric approach has certain advantage in
accommodating a larger class of bivariate distributions,
and hence is practically more applicable.

\begin{figure}

\includegraphics{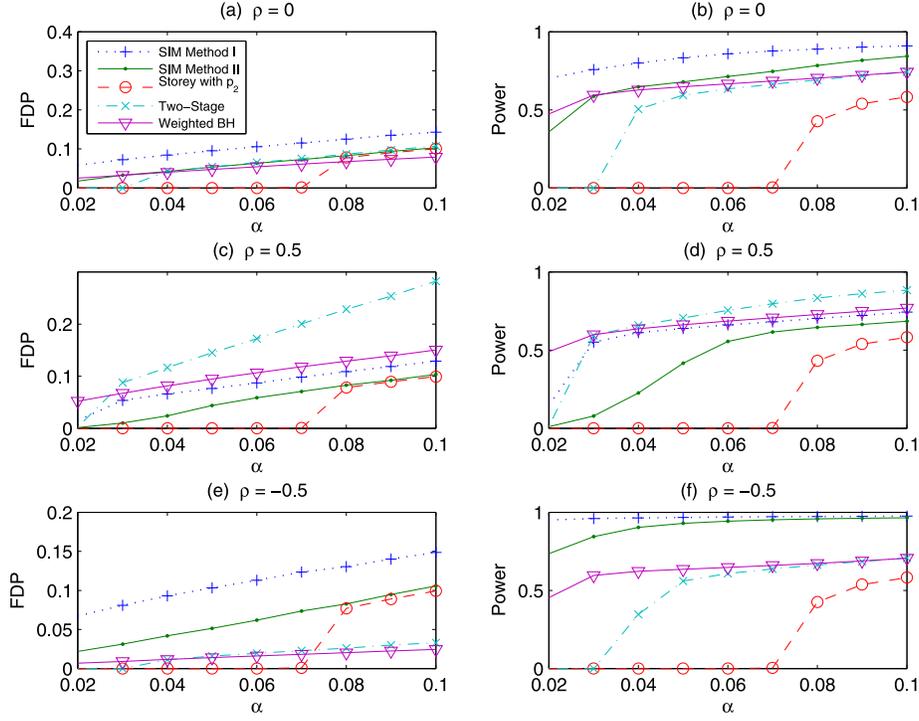}

\caption{Calculated $\FDP$ and power as a function of
$\alpha$ and $\rho$ for
the $\SIM$ methods, the weighed multiple testing procedure (weight
$\BH)$, the two-stage multiple testing procedure (two-stage) and
the conventional $\FDR$ procedure (storey with $p_2)$ in Example
$2$. Here, the nonnull is a mixture of three bivariate $t$
distributions with signals $(3,3)$, $(6,6)$ and $(8,8)$, respectively,
and $\pi_0=0.9$.}\label{Figure-11}
\end{figure}

Under a mixture of three bivariate $t$ distributions on the nonnull,
the comparison of the $\SIM$ methods with
the weighted multiple testing procedure and
the two-stage multiple testing procedure is demonstrated in Figure~\ref{Figure-11}. The story of bivariate $t$ distributions is similar to
that of bivariate normal models in
Figure~\ref{Figure-7} except that method~$\I$ loses its validity for
controlling $\FDR$ in all settings. In summary, method~$\II$ has the
merit of correctly and efficiently incorporating the prior information,
such as \emph{filters} in the two-stage multiple testing procedure
and \emph{weights} in the weighted multiple testing procedure, into
the conventional $\FDR$ procedure under any dependence structure
($\rho$).
%
%
\subsection{Example 3: Multiple testing with serially clustered signals} \label{Sec-63}
In practice, nonnull hypotheses are typically clustered.
Thus, we can take a preliminary \mbox{$p$-}value $p_{i1}$ to be the local
aggregation of $p_{j2}$, for $j$ located
in the neighborhood of the $i$th hypothesis,
where $\{ p_{i2}\dvtx  i=1,\ldots,m\}$ are the primary \mbox{$p$-}values.
The new pairs $\{ ( p_{i1}, p_{i2})\dvtx  i=1,\ldots,m\}$ consist of the
bivariate \mbox{$p$-}values.
In this example, we mimic the situation of serially clustered signals
to evaluate the performance of the $\SIM$ methods.
To be specific, we perform $m={}$10,000 one-sided hypotheses testing
independently, where test statistics follow $N(0,1)$ and $N(\mu, 1)$
for the true null and nonnull, respectively, for $\mu$ randomly chosen
from $\{1.5,2,2.5\}$.
The serial structure is designed as follows: the nonnull hypotheses
consist of three
clusters, that is, $\mathcal{C}_1=\{i=1001,\ldots,2000\}, \mathcal
{C}_2=\{ i=5001,\ldots,6000\}$
and $\mathcal{C}_3=\{i=8001,\ldots,9000\}$. There are various types
of preliminary \mbox{$p$-}values we can take, such as
the mean or median of the \mbox{$p$-}values in the neighborhood of the
original hypothesis; refer to \cite{ZhangFanYu2011} for details.
For simplicity, the \mbox{$p$-}values in the neighborhood of $p_{i2}$ is
chosen as
$\{p_{i-1,2}, p_{i+1,2}\}$ and the preliminary \mbox{$p$-}value is defined as
$p_{i1}=({ p_{i-1,2}+p_{i+1,2} })/{2}$, for $i=1,\ldots,m$.
Besides the conventional $\FDR$ procedure,
the mean filter, $p_i^{*}=({p_{i-1,2}+p_{i,2}+p_{i+1,2} })/{3}$
proposed by \cite{ZhangFanYu2011}, also serves as a competitor.
The results are shown in Table~\ref{Table-3}.
Method~II, the mean filter using $p_i^{*}$ and the conventional $\FDR$
procedure using $p_2$ provide conservative control of $\FDR$,
whereas $\FDR$ of method~I is slightly out of control for small
$\alpha$. This is reasonable as the normality assumption is not
strictly satisfied for the transformed \mbox{$p$-}value
$( \Phi^{-1}(p_{i1}), \Phi^{-1} (p_{i2}) )$.
In general, by utilizing the structural information of the primary \mbox{$p$-}values,
both method~II and the mean filter using $p_i^{*}$ are more powerful
than the conventional $\FDR$ procedure using $p_2$ alone. Rather than
giving the same weight $(1/3)$ to the neighborhood $(p_{i-1,2},
p_{i,2}, p_{i+1,2})$ in the mean filter $p_i^{*}$, the data-driven
procedure for selecting $\theta$ based on power comparison
for method~II adjusts different weights to the bivariate \mbox{$p$-}value
$(p_{i1}, p_{i2})$ according to their corresponding potential for
detecting power.
Consequently, method~II outperforms the mean filter using $p_i^{*}$ for
all possible $\alpha$.

%
\begin{table}
\tabcolsep=3pt
\caption{Calculated $\FDP$ and power comparison of methods~$\I$~and~$\II$, the mean filter using $p_i^{*}$ and the
conventional $\FDR$ procedure (storey with $p_2)$ in Example $3$}\label{Table-3}
\begin{tabular*}{\tablewidth}{@{\extracolsep{\fill}}@{}lcccccccc@{}} \hline
& \multicolumn{2}{c}{$\bolds{\FDR^{\I}}$ \textbf{using} $\bolds{(p_1, p_2)}$} & \multicolumn{2}{c}{$\bolds{\FDR^{\II}}$ \textbf{using} $\bolds{(p_1, p_2)}$}
& \multicolumn{2}{c}{\textbf{Mean filter using $\bolds{p_i^{*}}$}} & \multicolumn{2}{c@{}}{\textbf{Storey with $\bolds{p_2}$}}
\\[-6pt]
& \multicolumn{2}{c}{\hrulefill} & \multicolumn{2}{c}{\hrulefill} &\multicolumn{2}{c}{\hrulefill} & \multicolumn{2}{c@{}}{\hrulefill}
\\
$\bolds{\alpha}$ & $\bolds{\FDP}$ & \textbf{Power} & $\bolds{\FDP}$ & \textbf{Power} & $\bolds{\FDP}$ & \textbf{Power} & $\bolds{\FDP}$ & \textbf{Power}\\
\hline
0.01 &0.013 &0.617 &0.010 &0.578 &0.010 &0.505 &0.010 &0.059 \\
0.02 &0.024 &0.708 &0.020 &0.684 &0.020 &0.616 &0.019 &0.115 \\
0.03 &0.034 &0.759 &0.030 &0.742 &0.030 &0.682 &0.029 &0.164 \\
0.04 &0.044 &0.794 &0.040 &0.782 &0.040 &0.728 &0.038 &0.208 \\
0.05 &0.053 &0.820 &0.050 &0.811 &0.050 &0.763 &0.048 &0.247 \\
0.06 &0.063 &0.841 &0.060 &0.834 &0.060 &0.791 &0.058 &0.283 \\
0.07 &0.073 &0.858 &0.070 &0.852 &0.070 &0.813 &0.067 &0.317 \\
0.08 &0.082 &0.872 &0.079 &0.867 &0.080 &0.832 &0.077 &0.348 \\
0.09 &0.092 &0.884 &0.089 &0.881 &0.090 &0.849 &0.087 &0.377 \\
0.10 &0.101 &0.894 &0.099 &0.891 &0.100 &0.864 &0.096 &0.404 \\
0.20 &0.196 &0.952 &0.199 &0.953 &0.199 &0.945 &0.192 &0.610 \\
0.30 &0.293 &0.977 &0.299 &0.978 &0.299 &0.978 &0.288 &0.748 \\
\hline
\end{tabular*}
\end{table}

%
\subsection{Example 4: Two-sample $t$ test} \label{Sec-64}
\begin{figure}

\includegraphics{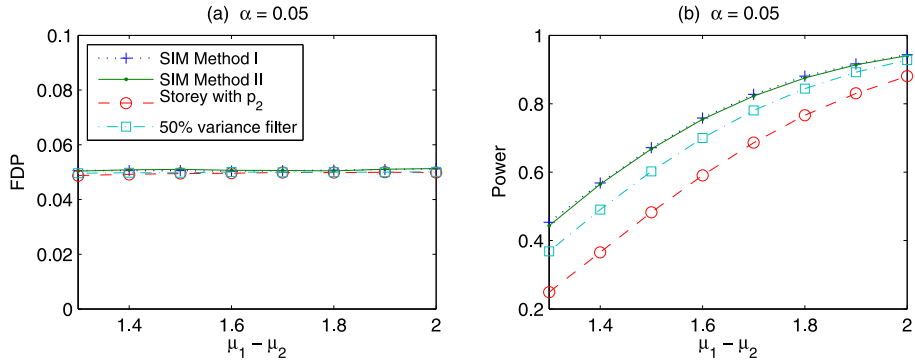}

\caption{Calculated $\FDP$ and power of methods $\I$ and $\II$, the $50\%$ variance filter and the conventional $\FDR$ procedure
(storey with $p_2)$ for various size effect $\mu_1-\mu_2$ in
Example $4$.}\label{Figure-12}
\end{figure}

In this example, we mimic the microarray experiment, where two-sample
$t$ test is performed to
detect differentially expressed genes for two classes comparison.
Suppose $m={}$10,000 genes are
\mbox{examined} independently, among which $10 \%$ are from the nonnull.
For the $i$th gene, let $\{ x_{i,1}, x_{i,2},\ldots, x_{i,10}\}$ and
$\{ y_{i,1}, y_{i,2},\ldots, y_{i,10}\}$ be two independent samples
from $N(\mu_1, 1)$ and $N(\mu_2, 1)$, respectively,
where $\mu_1=\mu_2$ is for nondifferentially expressed genes and $\mu
_1 \ne \mu_2$ is for differentially expressed genes.
The primary \mbox{$p$-}value, $p_{i2}$, is obtained by the standard two-sample
$t$ test.
To get the preliminary \mbox{$p$-}value, $p_{i1}$, the sum of squared error of
the two samples,
which has a chi-square distribution with $19$ degrees of freedom
and independent of $t$ statistic in the standard two-sample $t$ test
under the true null, can be utilized.
In this scenario, the independence between the components of bivariate
\mbox{$p$-}value implies that the true null distribution of the combined \mbox{$p$-}value
is uniform for all $\theta$.
To make a comprehensive comparison, the $50\%$ variance filter proposed
in \cite{Bourgonetal2010} is also considered.
Figure~\ref{Figure-12} shows that the performance of methods~I~and~II is almost the same and the corresponding power is improved
for different size effect $\mu_1-\mu_2$ for $\alpha=0.05$.
Particularly, our method is superior to the $50\%$ variance filter for
all cases. This is due to the fact that we employ
a data-driven procedure for choosing the tuning parameter $\theta$,
whereas the fraction $50 \%$ in the variance
filtering procedure is subjectively fixed.
%
%
\section{Integrative analysis on prostate cancer data} \label{Sec-7}
Genomic DNA copy number ($\CN$) alterations are key genetic events in
the development and progression of human cancers.
In parallel, microarray gene expression (GE) measurements of mRNA level
provide an alternative for detecting some significant genes which
contribute to certain cancer diseases.
As discussed by the previous study \cite{Kimetal2007}, the amplified
gene section was enriched with transcript overexpression, and the
deleted section was enriched with mRNA downregulation. Hence,
integration of $\CN$ aberration and GE to identify DNA $\CN$
alterations that induce changes in the expressional levels of the
associated genes is a common task in cancer studies. To this end,
several authors have explored integrative analysis of these two
heterogeneous data sources to reveal higher levels of interactions that
cannot be detected based on individual observations; see \cite
{Lahtietal2013} and the references therein.

\begin{figure}

\includegraphics{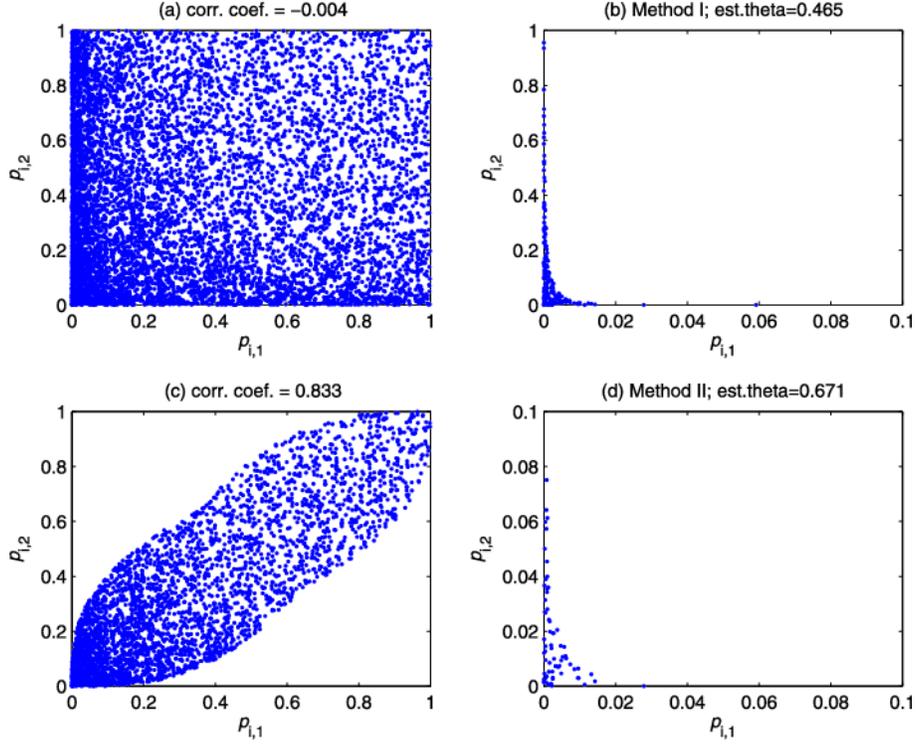}

\caption{\textup{(a)}: Scatter plot of bivariate \mbox{$p$-}values $(p_1, p_2)$, where the correlation coefficient of $p_1$ and $p_2$ is $-0.004$;
\textup{(b)}: geometric locations of the rejected\vspace*{1.5pt} genes using method~$\I$ with the significance level $\alpha=0.01$. Here, the projection
              direction is $\hat {\theta}{}^{\I}=0.465$;
\textup{(c)}: scatter plot of bivariate \mbox{$p$-}values $(p_1, p_2)$ of the trimmed genes, where the correlation coefficient of $p_1$ and $p_2$ is $0.833$;
\textup{(d)}: geometric locations of the rejected genes using method $\II$ for the trimmed genes with the significance level $\alpha=0.01$.
              Here, the projection direction is $\hat {\theta}{}^{\II}=0.671$.}\label{Figure-13}
\end{figure}

To demonstrate the practical utility of the $\SIM$ procedure, we
applied it to data produced by \cite{Kimetal2007} in a study on
prostate cancer progression. This study used an array comparative
hybridization (aCGH) to profile genome-wide $\CN$ changes through the
isolation of pure cell populations representing entire spectrum of
prostate disease using laser capture microdissection (LCM) and OmniPlex
Whole Genomic (WGA) Application. Data on $\CN$ alterations and GE were
matched for $m=7534$ genes using prostate cell populations from
low-grade ($n_1=27$) and high-grade samples ($n_2=17$) of cancerous
tissue. We calculated two-sided $t$ statistics ($t_1,t_2$) and their
\mbox{$p$-}values ($p_1,p_2$) for GE and $\CN$ aberrations for each of $7534$ genes.
Here, the primary \mbox{$p$-}value $p_2$ was obtained from the copy number in
DNA level and its transcriptional gene expression served as the
preliminary \mbox{$p$-}value $p_1$.
Panel (a) of Figure~\ref{Figure-13} shows the scatter plot of gene
expression and copy number \mbox{$p$-}values, where the sample correlation
coefficient of $p_1$
and $p_2$ is $-0.004$. This motivates us to apply our $\SIM$ method~$\I$ to target the genes evidencing
statistical significance in either DNA or mRNA level. Using the
significance level $\alpha=0.01$, our $\SIM$ method~$\I$ detects
$174$ rejections with their geometric locations showing in panel (b) of
Figure~\ref{Figure-13}.
The projection direction is estimated as $\hat {\theta}{}^{\I
}=0.465$, supporting that the preliminary \mbox{$p$-}value from GE is informative.

Note that our $\SIM$ procedure is valid for testing the conjunction of
null hypotheses to favor genes
with DNA copy number alterations \textit{or} differential expressions
under the alternative. Some genes are amplified or deleted in DNA level
but have insignificant GE in mRNA level, which can be accounted for by
the inappropriate use of ``methylation;'' while some upstream
``transcription factor'' genes found differentially expressed with
activation (or suppression) function will up (or down) downstream
genes. To further identify candidate genes with genetic alterations
that accompany corresponding transcriptomic changes, we utilized a
weight function, a product in DNA/RNA-Significance Analysis of
Microarrays (DR-SAM) \cite{Salarietal2010}, to screen out the genes
which are significant only in DNA or mRNA level. Specifically, the
weight function, which is defined as $w=\min\{ \frac{t_{1} }{t_{2}},
\frac{t_{2} }{t_{1}} \}$ ($0\le w \le1$), is the ratio of two
$t$-scores. Small weight is applied to favor genes with unbalanced
contributions on copy number and gene expression. Based on this
rationale, the genes with weights larger than a
threshold will serve as candidates for detecting concordantly altered
genes. Given a threshold, the scatter plot of genes passing the
threshold under the true null violates the normality and symmetry
property assumptions. Fortunately, the genes with points above the line
$p_1+p_2>1$ seldom come from the alternative. Hence, we modified the
weight function on the area with $p_1+p_2>1$ as
$w'(p_1,p_2)=w(1-p_1,1-p_2)$ such that the genes passing the threshold
satisfy the symmetry property assumption.
A~small threshold will enrich the alternative with some genes being
significant only in DNA or mRNA level, increasing the false discovery
rate; while a large threshold will screen out some genes exhibiting
concordant changes, resulting in low power. Based on this perspective,
the selection of the threshold using the modified weight function is
\textit{fdr-power} trade-off. For simplicity, we set the threshold
such that $50\%$ of the genes will be screened out. Panel (c) of
Figure~\ref{Figure-13} presents the scatter plot of the trimmed genes,
which will be used for testing. At $\alpha=0.01$, our $\SIM$ method
$\II$ estimates the projection direction as $\hat {\theta}{}^{\II
}=0.671$ and selects $62$ genes, as shown in panel (d) of Figure~\ref{Figure-13}. To make comprehensive comparisons, Table~\ref{Table-4}
shows the numbers of rejected genes by applying the $\SIM$ methods,
and the three competing procedures as used in our numerical studies.
In summary, all the three competing procedures with either $p_1$ or
$p_2$ as primary \mbox{$p$-}values, are more conservative than our $\SIM$ procedures.
%
\begin{table}[t] 
\tabcolsep=0pt
\tablewidth=253pt
\caption{Compare the numbers of rejections by the $\SIM$
methods, the conventional $\FDR$ procedure (storey), the two-stage
multiple testing procedure (two-stage), and the weighted multiple
testing procedure (weighted $\BH)$ when $\alpha=0.01$. Here, $(p_i
\mid p_j)$ indicates that $p_i$ is used as the primary \mbox{$p$-}value while
$p_j$ serves as the preliminary \mbox{$p$-}value}\label{Table-4}
\begin{tabular*}{\tablewidth}{@{\extracolsep{\fill}}@{}ld{3.0}@{}}
\hline
\textbf{Methods} & \multicolumn{1}{c@{}}{\textbf{Number of rejections}}\\
\hline
$\SIM$ method  $\I$ with the whole data & 174 \\
$\SIM$ method~$\II$ with the trimmed data & 62 \\[3pt]
Storey with $p_1$ & 31 \\
Storey with $p_2$ & 0 \\[3pt]
Two-stage with $(p_1\mid p_2)$ &16 \\
Two-stage with $(p_2 \mid p_1)$ & 1 \\[3pt]
Weighted $\BH$ with $(p_1\mid p_2)$ & 14 \\
Weighted $\BH$ with $(p_2 \mid p_1)$ & 0 \\
\hline
\end{tabular*}
\end{table}

%
\begin{table}[b] 
\tabcolsep=0pt
\caption{Summary of Gene Functional Classification from Gene Ontology $(\GO)$. $9$ $\GO$ terms are inferred to be active using
$\MFA$. Here, $\mathrm{P.MFA}$
represents the marginal posterior probability of activation, and basic
statistics on these terms are provided in the ``size'' column
(\#prostate cancer-associated genes/set size)} \label{Table-5}
\begin{tabular*}{\tablewidth}{@{\extracolsep{\fill}}@{}llcc@{}}
\hline
$\bolds{\GO}$ \textbf{ID} & \textbf{Gene set ($\bolds{\GO}$ term)} & \textbf{P.MFA} & \textbf{Size}\\
\hline
GO:0007031 & Peroxisome organization & 0.7909023 & 2${}/{}$58 \\
GO:0070307 & Lens fiber cell development & 0.7225289 & 1${}/{}$12 \\
GO:0001569 & Patterning of blood vessels & 0.7094174 & 2${}/{}$35 \\
GO:0001517 & N-acetylglucosamine 6-O-sulfotransferase activity &0.7036159 & 1${}/{}$6\phantom{0} \\
GO:0008455 & Alpha-1, 6-mannosylglycoprotein & 0.6962325 & 1${}/{}$1\phantom{0} \\
GO:0043190 & ATP-binding cassette (ABC) transporter complex & 0.6593146& 1${}/{}$6\phantom{0} \\
GO:0008332 & Low voltage-gated calcium channel activity & 0.6440800 &1${}/{}$3\phantom{0} \\
GO:0030612 & Arsenate reductase (thioredoxin) activity & 0.6339682 &1${}/{}$1\phantom{0} \\
GO:0004464 & Leukotriene-C4 synthase activity & 0.6276624 & 1${}/{}$2\phantom{0} \\
\hline
\end{tabular*}
\end{table}

Of these $62$ genes selected by our $\SIM$ method $\II$ with the
trimmed genes, $38$~were mapped to the official gene names (11,705 in
total) for prostate cancer with somatic mutation listed on Catalogue of
Somatic Mutation in Cancer \mbox{(COSMIC)}, supporting these genes being
putative oncogenes in prostate cancer. Notably, the top five genes,
that is, ABCA4, ABCA3, ACTG1, AADAC and ACACA, were ranked as 426, 454,
700, 780 and $848$, respectively. Particularly, the gene ACACA, known
to be involved in fatty and acid metabolism, was also identified in the
previous study \cite{lapointeetal2004}.
To integrate gene-set information from a complex system with our
experimentally-derived gene
list, a larger gene list is necessary. For this purpose, we performed
our $\SIM$ method $\II$ to the trimmed genes
at $\alpha=0.05$, which yields $331$ rejections. Among them, $102$
could be mapped to recognized genes by DAVID \cite{Huangetal2008}.
To assess the functional content of this gene list, we applied a new
approach termed as multifunctional analyzer ($\MFA$) proposed by \cite
{Wangetal2013}, in the context of gene ontology terms. Compared with
existing methods such as Fisher's exact test and model-based gene-set
analysis (MGSA) \cite{Baueretal2010}, $\MFA$ has the merit of
alleviating the redundancy problem in Fisher's exact test while
improving the statistical efficiency of MGSA. Table~\ref{Table-5}
reports the gene sets which were inferred to be activated by $\MFA$ in
prostate cancer.

%
\section{Discussion} \label{Sec-8}
This paper proposes a $\SIM$ multiple testing procedure to embed prior
information, such as the overall sample variance in a standard
two-sample $t$ test in microarray
experiments and the structurally spatial information for large-scale
imaging data, into the conventional $\FDR$ procedure, by assuming the
availability of
a bivariate \mbox{$p$-}value for each null hypothesis. We discuss the optimal
rejection region in terms of power comparison in a general bivariate
model and
project the bivariate \mbox{$p$-}value into a single-index quantified by a
projection direction $\theta$. A novel procedure is established to estimate
the optimal projection direction consistently under some mild
conditions, followed by two procedures for the estimation and control
of $\FDR$.

Although the operators $\Phi$ and $\Phi^{-1}$ in the single-index
$p(\theta)$ come from the normality assumption, generalizations, such as
$p(\theta)= \Psi(\cos(\theta) \Psi^{-1}(p_1) +\sin(\theta) \Psi
^{-1}(p_2))$, can be made, where $\Psi$ is the $\CDF$ of some random
variable. We have shown in the simulation study that the normal
operator $\Phi$ is robust
to distributions of other bivariate test statistics.
A thorough investigation of the role of the operator is beyond the
scope of this paper, but could be of interest
in the future research.

As discussed in Section~\ref{Sec-3}, the essential spirit of multiple
testing is on increasing the detection power
while maintaining the $\FDR$ rigorously.
Theoretically, the detection power is related to three quantities, that
is, $\pi_0$, $F_0(t)$ and $F_1(t)$, via the Bayesian $\mathrm{Fdr}$
formula $F_1(t)=({1}/{\alpha}-1){\pi_0}/{\pi_1}F_0(t)$.
Screening out a proportion of uninformative hypotheses by an effective filter
will enrich for nonnull hypotheses while simultaneously reducing the
number of
hypotheses to be tested at the second stage.
From this point of view, the independence filter provided by \cite
{Bourgonetal2010} aims to decrease $\pi_0$ to improve the detection power.
However, in our $\SIM$ multiple testing procedure, we project the
bivariate \mbox{$p$-}value into a single-index, which significantly changes
the true null and nonnull distributions $(F_0(t), F_1(t))$. Hence,
the power is increased by changing the structure of \mbox{$p$-}values while
keeping $\pi_0$ to be constant. Our future research will be
focused on constructing a more powerful multiple testing procedure via
reducing the proportion of true null hypothesis and
changing the structure of \mbox{$p$-}values simultaneously.

Beyond the weak dependence assumption made in \textup{(C2)}, the
sequence of the projected \mbox{$p$-}values will inevitably inherit strong
dependence from the primary test statistics, making the $\SIM$
procedure less accurate. Much published work has been developed to
handle multiple testing problem with some strong dependence structure;
see \cite{Fanetal2012} and the references therein. Much research is
needed to investigate the performance of the $\SIM$ methods for
solving multiple testing problem with strong dependence structure
across the tests.


\begin{appendix}
\section{Proofs of main results}\label{appA}

For presentational fluency, denote $\widetilde F_0(t,\theta)={ V(t,
\theta) }/ {m_0}$,
$\widetilde F_1(t,\theta)=\break \{ R(t,\theta)-V(t, \theta)\}/{m_1}$ and
$\widehat F(t, \theta)=R(t,\theta)/m$.
Analogously, define the following left-limit processes:
\begin{eqnarray*}
\widetilde F_0(t-, \theta) &=& m_0^{-1} \sum_{i=1}^{m} \ID\bigl\{ p_i(\theta) <t, H_0(i) \bigr\},
\\
\widetilde F_1(t-,\theta) &=& m_1^{-1} \sum_{i=1}^{m}\ID\bigl\{ p_i(\theta) <t, H_1(i) \bigr\},
\\
\widehat F(t-, \theta) &=& m^{-1} \sum_{i=1}^{m}\ID\bigl\{p_i(\theta) <t \bigr\}.
\end{eqnarray*}

We only prove the main results involved the nonparametric estimator
$\widehat F_0^{ \II}(t, \theta)$.
For those involved the parametric estimator $\widehat F_0^{ \I}(t,
\theta)$, all proofs will go through
as long as this estimator uniformly converges to the true null
distribution $F_0(t, \theta)$ for all $t$ and $\theta$.

We first impose some regularity conditions, which are not the weakest
possible but facilitate the technical derivations.

\subsection*{Conditions}
\begin{longlist}[(C10)]
\item[(C1)]
$\lim_{m \to \infty} m_0/m=\pi_0$ exists and $ 0<\pi_0 <1$.

\item[(C2)]
$\lim_{m \to\infty}
{m_0}^{-1}{\sum_{i=1}^{m} \ID(p_{ij} \le t, H_0(i))} =G_0^{j} (t) $
and $\lim_{m \to\infty} {m_1}^{-1}\* {\sum_{i=1}^{m} \ID(p_{ij} \le
t, H_1(i))} =G_1^{j} (t) $ almost surely,
for $j=1,2$.

\item[(C3)]
For\vspace*{2pt} any rational number $\alpha\in[0,1]$, denote by
$q_{\alpha}(\theta)$ the $100\alpha$th quantile of the distribution
function $F(t, \theta)$. Assume that $\widehat F (t, \theta)$ and\break
$F(t, \theta)$ satisfy the Lipschitz continuity as follows:
$\sup_{m} \sup_{\alpha} | \widehat F (q_{\alpha}(\theta),
\theta)-\break
\widehat F (q_{\alpha}(\theta'), \theta') | \le C_1 | \theta
-\theta' |$
and $\sup_{\alpha} |F(q_{\alpha}(\theta), \theta)- F(q_{\alpha
}(\theta'), \theta') | \le C_1 |\theta-\theta'|$, where\vspace*{1pt} $C_1$
is a generic positive constant, not depending on $\widehat F, F$ and
$\alpha$. The Lipschitz continuity conditions
also hold for $\widehat F (t-,\theta)$ and $F(t-, \theta)$.
In addition, $F_0(t, \theta)$, $F_0(t-, \theta)$, $\widetilde F_0(t,
\theta)$ and
$\widetilde F_0(t-, \theta)$ satisfy the Lipschitz continuity conditions.

\item[(C4)]
The probability density function of $(p_1,p_2)$ under the true null is
centrally symmetric with respect to $(1/2,1/2)$.

\item[(C5)] $F_1( {1}/{2}, \theta)=1$ for all $\theta$.

\item[(C6)]
$\inf_{\theta}F(\delta, \theta)>0$, for any $\delta>0$.

\item[(C7)]
$F_0(t,\theta)$ and $F(t,\theta)$ are continuous in the region $\{(t,
\theta)\dvtx  t^*_{\alpha'}(\theta) \le t \le1 \}$
and $ |F(t, \theta)-F(t^*_{\alpha'}(\theta), \theta)| \le
C_2|t-t^*_{\alpha'}(\theta) |$, where $C_2$ is a constant not
depending on $\theta$.\vspace*{2pt}

\item[(C8)]
\[
\lim_{t \to t^*_{\alpha'}(\theta) } \frac{ F_0(t, \theta)/F(t,\theta)- F_0(t^*_{\alpha'}(\theta),
\theta)/F(t^*_{\alpha'}(\theta), \theta)}{
t-t^*_{\alpha'}(\theta) }=k(\theta)
\]
uniformly for $\theta$, where $\inf_{\theta} |k(\theta)| >0$.

\item[(C9)] (Identification). Given $\delta'>0$, there exists
$\varepsilon>0$, such that
\[
\inf_{\theta\dvtx  |\theta-\theta_0(\alpha') |> \delta' } \bigl\{ F\bigl(
t^*_{\alpha'}\bigl( \theta_0(\alpha') \bigr),\theta_0(\alpha') \bigr)-F\bigl(
t^*_{\alpha'}(\theta), \theta\bigr)\bigr\}
\ge\varepsilon.
\]

\item[(C10)]
$|F(t, \theta)-F(t, \theta_0(\alpha' ) )|\le C_3 | \theta- \theta
_0(\alpha') | $ and $|F_0(t, \theta)-F_0(t, \theta_0(\alpha') )|
\le C_3 | \theta- \theta_0(\alpha') |$, where $C_3$ is a constant
not depending on $\theta$ and $t$.
\end{longlist}

Before proving the propositions and theorems, we first show Lemmas \ref
{Lemma-1} and \ref{Lemma-2}.
%
\begin{lemma} \label{Lemma-1}
$\!\!\!$Assume conditions \textup{(C1}) to \textup{(C3)}.
Let $p_i(\theta)\,{=}\,\Phi( \cos(\theta) \Phi^{-1}(p_{i1})\!+ \sin
(\theta) \Phi^{-1}(p_{i2}))$, $i=1,\ldots,m$,
where $\Phi$ is the $\CDF$ of a standard normal random variable.
Then we have
\begin{eqnarray*}
\sup_{0 \le\theta\le{\pi}/{2} } \sup_{0 \le t \le1} \Biggl|\frac{1}{m_0}
\sum_{i=1}^{m} \ID\bigl\{p_{i}(
\theta) \le t, H_0(i) \bigr\} -F_0(t, \theta) \Biggr| &
\stackrel{\mathit{a.s.}} \to & 0,
\\
\sup_{0 \le\theta\le{\pi}/{2} } \sup
_{0 \le t \le1} \Biggl|\frac{1}{m_1} \sum_{i=1}^{m}
\ID\bigl\{p_{i}(\theta) \le t, H_1(i) \bigr
\}-F_1(t, \theta) \Biggr| &\stackrel{\mathit{a.s.}} \to & 0,
\\
\sup
_{0 \le\theta\le{\pi}/{2} } \sup_{0 \le t \le1} \Biggl| \frac{1}{m} \sum
_{i=1}^{m} \ID\bigl\{p_{i}(\theta)
\le t\bigr\}-F(t, \theta) \Biggr| &\stackrel{\mathit{a.s.}} \to & 0.
\end{eqnarray*}
\end{lemma}

\begin{pf}
We first show the uniform consistency of $\widehat F(t, \theta)$. For
fixed $t$ and $\theta$,
$\{ p_i(\theta)\dvtx  i=1,\ldots,m \}$ satisfy the weak dependence:
%
\begin{eqnarray} \label{A1}
\Biggl|\frac{1}{m} \sum_{i=1}^{m} \ID\bigl
\{p_{i}(\theta) \le t\bigr\} -F(t, \theta) \Biggr| &\stackrel{\mathrm{a.s.}}
\to & 0,
\nonumber\\[-8pt]\\[-8pt]
\Biggl|\frac{1}{m} \sum_{i=1}^{m}
\ID\bigl\{p_{i}(\theta) < t\bigr\} -F(t-, \theta) \Biggr| &\stackrel{
\mathrm{a.s.}} \to & 0.\nonumber
\end{eqnarray}
This conclusion is directly implied by conditions \textup{(C1}) and
\textup{(C2)}. To prove the uniform consistency of
$\widehat F(t, \theta)$, we extend the argument in the proof of
\emph{the Glivenko--Cantelli theorem} \cite{Durrett2010}.
For $0 \le j \le k $, partitioning the domain into grid points
$(t, \theta)$ as $\{ q_{j/k}(\theta_l)\dvtx  j=0,\ldots, k;  l=0,\ldots, L_k\}$
such that $\{ \theta_l\dvtx  l=0,\ldots, L_k\}$ are equally spaced in $[0,
\pi/2]$
with unit length less than or equal to ${1}/{(C_1 k)}$, where $C_1$ is
given in condition \textup{(C3)}. The pointwise convergence (\ref
{A1}) implies that we can pick up $N_k(\omega)$ such that
%
\begin{eqnarray}\label{A2}
\bigl|\widehat F\bigl( q_{j/k}(\theta_l),
\theta_l\bigr)-F\bigl( q_{j/k}(\theta_l),
\theta_l\bigr)\bigr| &<& k^{-1}\quad\mbox{and}
\nonumber\\[-8pt]\\[-8pt]
\bigl|\widehat F\bigl(
q_{j/k}(\theta_l)-, \theta_l\bigr)-F\bigl(
q_{j/k}(\theta_l)-, \theta_l\bigr)\bigr| &<& k^{-1}\nonumber
\end{eqnarray}
for\vspace*{2pt} $ 0 \le j \le k$ and $ 0 \le l \le L_k$.
For $t \in( q_{{ (j-1) }/k} (\theta), q_{{j}/k} (\theta))$ and
$\theta\in(\theta_{l-1}, \theta_l)$ with $ 1 \le j \le k $, $ 1 \le
l \le L_k$ and $m> N_k (\omega)$, using the
monotonicity of $\widehat F$ and $F$, $F(q_{j/k}(\theta)-, \theta
)-F( q_{ {j-1}/k }(\theta), \theta)\le k^{-1}$
and condition \textup{(C3)}, we have
\begin{eqnarray*}
\widehat F(t, \theta) & \le& \widehat F\bigl( q_{j/k}(\theta)-, \theta
\bigr)
\\
&\le& \widehat F\bigl( q_{j/k}(\theta_{l-1})-,
\theta_{l-1}\bigr)+k^{-1}
\\
&\le& F\bigl( q_{j/k}(
\theta_{l-1})-, \theta_{l-1} \bigr)+2k^{-1}
\\
&
\le& F\bigl(q_{ {j-1}/k } (\theta_{l-1}), \theta_{l-1}
\bigr)+3k^{-1}
\\
& \le& F\bigl(q_{ {j-1}/k } (\theta), \theta
\bigr)+4k^{-1}
\\
& \le& F( t, \theta)+4k^{-1}.
\end{eqnarray*}
Similar arguments lead to $\widehat F(t, \theta) \ge F(t, \theta) - 4k^{-1}$.
So $\sup_{\theta}\sup_{t}| \widehat F(t, \theta)- F(t, \theta)|
\le4k^{-1}$, and we have proved the result.
The uniform convergence of $\widetilde F_0(t, \theta)$ can be derived
similarly. Combining these two
results, we obtain the uniform convergence of $\widetilde F_1(t,
\theta)$ immediately.
\end{pf}

%
\begin{lemma} \label{Lemma-2}
Under conditions \textup{(C1})--\textup{(C5)}, the nonparametric estimator
$\widehat F_0^{ \II}(t, \theta)$ uniformly converges to $F_0(t,
\theta)$ for all $t$ and $\theta$.
\end{lemma}

\begin{pf*}{Proof of Proposition \ref{Proposition-1}}
Let $\ID(\mathbf{p}\in\boldsymbolS_{\OR})$ be the indicator of
$\mathbf{p} \in \boldsymbolS_{\OR}$
and $\ID( \mathbf{p}\in \boldsymbolS)$ be the indicator of $\mathbf
{p} \in\boldsymbolS$ for any rejection region satisfying $\mathrm
{Fdr}( \boldsymbolS) \le\alpha$.
Since $\mathrm{Fdr}( \boldsymbolS_{\OR} )$ is the conditional
expectation of $\mathrm{fdr}(\mathbf{p})$
given $\mathbf{p}\in\boldsymbolS_{\OR}$ \cite{EfronTibshirani2002}, some derivations yield that
\begin{eqnarray*}
\alpha &=&\mathrm{Fdr}(\boldsymbolS_{\OR}) =E_{f}\bigl\{
\mathrm{fdr}(\mathbf{p}) \mid \mathbf{p} \in \boldsymbolS_{\OR} \bigr\}
\\
&=& \int_{ \{ \mathbf{p}\dvtx   \mathrm{fdr}(\mathbf{p}) <C \} } \mathrm{fdr}(\mathbf{p}) \,d\widetilde
\boldsymbolP + \int_{ \{ \mathbf{p}\dvtx   \mathrm{fdr}(\mathbf{p}) =C \} } \mathrm {fdr}(\mathbf{p}) \,d
\widetilde \boldsymbolP
\\
&=&\int_{ \{ \mathbf{p}\dvtx    \mathrm{fdr}(\mathbf{p}) <C \} } \mathrm{fdr}(
\mathbf{p}) \,d\widetilde \boldsymbolP<C,
\nonumber
\end{eqnarray*}
where $\widetilde \boldsymbolP$ denotes the probability measure of
$\mathbf{p}$ given $ \mathbf{p} \in \boldsymbolS_{\OR}$ and the
last equality holds by condition (\ref{36}).
As a result, $C>0$ and $1-\alpha/C>0$. By condition (\ref{36}),
there exists $\boldsymbolS'$ such that $\boldsymbolS\subseteq
\boldsymbolS'$ and $\mathrm{Fdr}(\boldsymbolS')=\alpha$.
For every~$\mathbf{p}$,
%
\begin{equation}
\label{A3} \ID\bigl(\mathbf{p}\in\boldsymbolS' \bigr) \bigl\{
1- \mathrm{fdr}( \mathbf {p}) /C \bigr\} \le \ID(\mathbf{p}\in
\boldsymbolS_{\OR}) \bigl\{ 1-\mathrm {fdr}(\mathbf{p}) /C \bigr\},
\end{equation}
where (\ref{A3}) is based on the observation that if $\mathbf
{p}\notin\boldsymbolS_{\OR}$, the left-hand side of (\ref{A3}) is less
than or equal to zero.
By taking expectation for both sides of equation (\ref{A3}),
\[
\int\ID\bigl(\mathbf{p}\in\boldsymbolS' \bigr) \bigl\{ 1-
\mathrm{fdr}( \mathbf{p}) /C \bigr\} f(\mathbf{p}) \,d\mathbf{p} \le \int\ID(
\mathbf{p}\in\boldsymbolS_{\OR} ) \bigl\{ 1-\mathrm {fdr}(\mathbf{p}) /C
\bigr\} f(\mathbf{p}) \,d\mathbf{p},
\]
we obtain the following inequality:
%
\begin{equation}
\label{A4} F\bigl(\boldsymbolS'\bigr) \bigl\{ 1-\mathrm{Fdr}
\bigl(\boldsymbolS'\bigr)/C \bigr\} \le F(\boldsymbolS_{ \OR}
) \bigl\{ 1-\mathrm{Fdr}(\boldsymbolS_{ \OR
})/C \bigr\},
\end{equation}
where $F(\boldsymbolS)=\pi_0F_0(\boldsymbolS)+\pi_1
F_1(\boldsymbolS)$.
By definition, both $1- \mathrm{Fdr}( \boldsymbolS' )/C$ and
$1-\mathrm{Fdr}(\boldsymbolS_{\OR})/C$ are equal to $1-\alpha/C>0$.
Hence, (\ref{A4}) implies that $F( \boldsymbolS' ) \le F(
\boldsymbolS_{\OR} )$. From the $\mathrm{Fdr}$ formula,
$F(\boldsymbolS')={\pi_1}/{(1-\alpha)} F_1(\boldsymbolS' )$ and
$F(\boldsymbolS_{\OR})={\pi_1}/{(1-\alpha)}F_1(\boldsymbolS_{\OR
})$. So $F_1(\boldsymbolS') \le F_1( \boldsymbolS_{\OR})$. The proof
is completed by the fact that $F_1(\boldsymbolS) \le F_1(\boldsymbolS
')$ for any $\boldsymbolS\subseteq\boldsymbolS'$.
\end{pf*}

\begin{pf*}{Proof of Proposition \ref{Proposition-2}}
By\vspace*{2pt} continuity, $t^*_{\alpha'}(\theta)$ satisfies that $F_0(t, \theta
)/\break  F(t, \theta)=\alpha'$.
From\vspace*{1pt} the $\mathrm{Fdr}$ formula, for any $\theta$, $t^*_{\alpha
'}(\theta)$ is the solution of the equation
$F_1(t, \theta)=\beta F_0(t, \theta)$.
Since $\frac{\partial F_1(t, \theta)}{\partial t}-\beta\frac
{\partial F_0(t, \theta)}{\partial t}\ne0$
for any interior point $(t, \theta, \alpha')$ in $[0,1] \times[0,
\pi/2] \times[0, 1/\pi_0]$, \emph{implicit function theorem}
implies that there exists a unique continuously differentiable function
$t=g(\alpha',\theta)$ such that
$F_1(g(\alpha', \theta), \theta)=\beta F_0(g(\alpha', \theta),
\theta)$.
By uniqueness, $t^*_{\alpha'}(\theta)=g(\alpha', \theta)$,
indicating that $t^*_{\alpha'}(\theta)$ is continuously
differentiable with respect
to $\theta$ and $\alpha'$.
Taking derivative with respect to $\theta$ for both sides of
$F_1( t^*_{\alpha'}(\theta), \theta)=\beta F_0(t^*_{\alpha'}(\theta
), \theta)$ leads to
%
\begin{eqnarray}\label{A5}
&& \biggl\{ \frac{\partial F_1(t, \theta)}{\partial t} \frac{\partial t^*_{\alpha' } (\theta) }{\partial\theta}+
\frac{\partial F_1( t, \theta ) }{\partial\theta} \biggr\} \bigg|_{t=t^*_{\alpha' }(\theta)}
\nonumber\\[-8pt]\\[-8pt]
&&\qquad = \biggl\{ \beta\frac{\partial F_0(t, \theta)}{\partial t}
\frac{\partial t^*_{\alpha' } (\theta) }{\partial\theta}+ \beta\frac{\partial F_0( t, \theta ) }{\partial\theta} \biggr\} \bigg|_{t=t^*_{\alpha' }(\theta)}.\nonumber
\end{eqnarray}
From (\ref{A5}), $\frac{ \partial t^*_{\alpha' } (\theta
)}{\partial\theta } $
can be expressed as
%
\begin{equation}
\label{A6} \frac{ \partial t^*_{\alpha'}(\theta) }{ \partial\theta}= \biggl\{ \frac{ \beta
(({\partial F_0( t, \theta ) })/({\partial\theta}))-
(({\partial F_1( t, \theta ) })/({ \partial\theta})) }{
(({\partial F_1 ( t, \theta ) })/({\partial t}))- \beta
(({\partial F_0( t, \theta ) })/({ \partial t})) } \biggr\}
\bigg|_{t=t^*_{\alpha' }(\theta)}.
\end{equation}
Since $F_1( t^*_{\alpha' }(\theta),\theta)$ achieves the maximum at
$\theta_0(\alpha' )$,
the following partial differential equation holds, that is,
%
\begin{equation}
\label{A7} \biggl[ \biggl\{ \frac{\partial F_1( t, \theta ) }{\partial t} \frac{\partial t^*_{\alpha' }(\theta)} {\partial\theta}+
\frac{\partial F_1( t, \theta ) }{\partial\theta} \biggr\} \bigg|_{t=t^*_{\alpha'}(\theta)} \biggr] \bigg|_{\theta=\theta_0(\alpha')} =0.
\end{equation}
Plugging~(\ref{A6}) into~(\ref{A7}), the partial differential
equation can be simplified as
%
\begin{eqnarray}\label{A8}
&& \biggl[ \biggl\{ \frac{\partial F_1( t, \theta) }{\partial t } \Big/ \frac{\partial F_1( t, \theta)}{\partial\theta} \biggr\}
\bigg|_{t=t^*_{\alpha'}(\theta)} \biggr] \bigg|_ {\theta
=\theta_0(\alpha')}
\nonumber\\[-8pt]\\[-8pt]
&&\qquad  = \biggl[ \biggl\{ \frac{\partial F_0( t, \theta)}{\partial t }
\Big/ \frac{\partial F_0( t, \theta) }{\partial\theta} \biggr\} \bigg|_{t=t^*_{\alpha'}(\theta)} \biggr] \bigg|_ {\theta
=\theta_0(\alpha')}.\nonumber
\end{eqnarray}
From $(\ref{A8})$, the $x$-coordinate of the point of intersection of
the solution set $\{(\theta, t)\}$
satisfying (\ref{46}) and $t=t^*_{\alpha'}(\theta)$ is $\theta
_0(\alpha')$.

``$\Leftarrow$'':
$\theta_{0}(\alpha')$ is constant for all $0<\alpha'<1/\pi_0$ if
the solution $\theta$ of $t$ of the equation (\ref{46}) is unique
and equals a constant.

``$\Rightarrow$'':
If the solution $\theta$ of $t$ of the equation (\ref{46}) is either
not unique or not equal to a constant, then there exists $t_1$ and
$t_2$ such that $\theta(t_1) \ne\theta(t_2)$. Since $t^*_{\alpha
'}(\theta)$
are continuous and nondecreasing from $[0, 1/\pi_0]$ with respect to
$\alpha'$ for any $\theta$, there exists $\alpha'_1$ and
$\alpha'_2$ such that $t^*_{\alpha'_1}(\theta(t_1))=t_1$ and
$t^*_{\alpha'_2}(\theta(t_2))=t_2$. From
(\ref{A8}), $\theta(t_1)=\theta_0(\alpha'_1)$ and $\theta
(t_2)=\theta_0(\alpha'_2)$, which implies that
$\theta_0(\alpha')$ is not constant for all $0<\alpha'<1/\pi_0$.

Under the normality assumption, $F_1(t,\theta)= \Phi( \frac{ \Phi
^{-1}(t)-\mu_1(\theta)}{\sigma_0(\theta)} ) $
and $F_0(t,\theta)= \Phi( \frac{ \Phi^{-1}(t)-\mu_0(\theta
)}{\sigma_0(\theta)} ) $,
where\vspace*{2pt} $\mu_0(\theta)=\mu_{0;1}\cos(\theta)+\mu_{0;2} \sin(\theta
)$, $\mu_1(\theta)=\mu_{1;1}\cos(\theta)+\mu_{1;2} \sin(\theta
)$ and $\sigma_0(\theta)$ is appearing in (\ref{41}).
In this case, (\ref{A8}) reduces to
$
 [ \frac{ \partial}{ \partial\theta} \{ \frac{ \mu_0(\theta)
}{ \sigma_0(\theta)} \} =
\frac{ \partial}{ \partial\theta} \{ \frac{ \mu_1(\theta) }{
\sigma_0(\theta)} \}  ]
 {|}_{\theta=\theta_0(\alpha' )}
$,
implying that $\theta_0(\alpha')$ is constant.
\end{pf*}

\begin{pf*}{Proof of Proposition \ref{Proposition-3}}
For left-sided hypotheses, the joint $\CDF$ of $(\tilde p_1,
p_2)$ under the true null can be derived as
%
\begin{eqnarray} \label{A9}
\qquad && \P ( \tilde p_1 \le\tilde t_1,
p_2 \le t_2 \mid H_0 )\nonumber
\\
&&\qquad = \P \bigl(
F_{0; X_1}(\widetilde X_1 ) \le\tilde t_1,
p_2 \le t_2 \mid H_0 \bigr)
\nonumber\\[-8pt]\\[-8pt]
&&\qquad = \P
\bigl( p_1 \le F_{0; X_1}\bigl( F_{0; X_1}^{-1}(
\tilde t_1)-\eta \bigr), p_2 \le t_2 \mid
H_0 \bigr)\nonumber
\\
&&\qquad = \int_{-\infty}^{\infty}
f_{\eta}(\eta) \biggl\{ \int_{0}^{F_{0; X_1}( F_{0; X_1}^{-1}(\tilde t_1)-\eta)}
\int_{0}^{t_2} f_{0; (p_1,p_2)}(p_1,
p_2)\,dp_1\,dp_2 \biggr\} \,d\eta,\nonumber
\end{eqnarray}
where $f_{\eta}$ is the p.d.f. of $\eta$, $f_{0; (p_1,p_2)}$
is the p.d.f. of $(p_1,p_2)$ under the true null,
and $F_{0;X_1}^{-1}$ is the inverse function of $F_{0; X_1}$.
By taking derivatives of (\ref{A9}), we obtain
%
\begin{eqnarray} \label{A10}
&& f_{0; (\tilde p_1, p_2)}(\tilde p_1, p_2)\nonumber
\\
&&\qquad = \int
_{-\infty}^{\infty} f_{\eta}(\eta)
\frac{ f_{0; X_1}( F_{0;
X_1}^{-1}(\tilde p_1)-\eta) }{f_{0; X_1}( F_{0;
X_1}^{-1}(\tilde p_1) )}
\\
&&\hspace*{52pt}{}\times  f_{0; (p_1,p_2)} \bigl( F_{0; X_1}\bigl(
F_{0; X_1}^{-1}(\tilde p_1)-\eta\bigr),
p_2 \bigr) \,d \eta,\nonumber
\end{eqnarray}
where $f_{0; (\tilde p_1, p_2)}$ is the p.d.f. of
$(\tilde p_1, p_2)$ under the true null and $f_{0; X_1}$ is the
true null p.d.f. of $X_1$. Thus,
%
\begin{eqnarray}
& & f_{0; (\tilde p_1, p_2) }(1-\tilde p_1, 1-p_2)\nonumber
\\
&&\qquad =
\int_{-\infty}^{\infty} f_{\eta}(\eta)
\frac{ f_{0; X_1}(
F_{0; X_1}^{-1}(1-\tilde p_1)-\eta) }{f_{0; X_1}( F_{0;
X_1}^{-1}(1-\tilde p_1) )}\nonumber
\\
&&\hspace*{53pt}{}\times  f_{0; (p_1, p_2)} \bigl( F_{0; X_1}\bigl(
F_{0; X_1}^{-1}(1-\tilde p_1)-\eta\bigr),
1-p_2 \bigr) \,d \eta\nonumber
\\
&&\qquad = \int_{-\infty}^{\infty}
f_{\eta}(\eta) \frac{ f_{0; X_1}(
-F_{0; X_1}^{-1}(\tilde p_1)-\eta) }{f_{0; X_1}( -F_{0;
X_1}^{-1}(\tilde p_1) )}\nonumber
\nonumber\\[-8pt]\label{A11} \\[-8pt]
&&\hspace*{53pt}{}\times f_{0; (p_1, p_2)} \bigl(
F_{0; X_1}\bigl( -F_{0; X_1}^{-1}(\tilde p_1)-\eta\bigr), 1-p_2 \bigr) \,d \eta \nonumber
\\
&&\qquad = \int_{-\infty}^{\infty} f_{\eta}(\eta)
\frac{ f_{0; X_1}(
F_{0; X_1}^{-1}(\tilde p_1)+\eta) }{f_{0; X_1}( F_{0;
X_1}^{-1}(\tilde p_1) )}
\nonumber\\[-8pt] \label{A12}\\[-8pt]
&&\hspace*{53pt}{}\times f_{0; (p_1, p_2)} \bigl( 1-F_{0; X_1}\bigl(
F_{0; X_1}^{-1}(\tilde p_1)+\eta\bigr),
1-p_2 \bigr) \,d \eta\nonumber
\\
&&\qquad = \int_{-\infty}^{\infty} f_{\eta}(\eta)
\frac{ f_{0; X_1}(
F_{0; X_1}^{-1}(\tilde p_1)+\eta) }{f_{0; X_1}( F_{0;
X_1}^{-1}(\tilde p_1) )}
\nonumber\\[-8pt] \label{A13} \\[-8pt]
&&\hspace*{53pt}{}\times f_{0; (p_1, p_2)} \bigl( F_{0; X_1}\bigl(
F_{0; X_1}^{-1}(\tilde p_1)+\eta\bigr),
p_2 \bigr) \,d \eta\nonumber
\\
&&\qquad = \int_{-\infty}^{\infty} f_{\eta}(-\eta)
\frac{ f_{0; X_1}(
F_{0; X_1}^{-1}(\tilde p_1)-\eta) }{f_{0; X_1}( F_{0;
X_1}^{-1}(\tilde p_1) )}\nonumber
\\
&&\hspace*{53pt}{}\times f_{0; (p_1, p_2)} \bigl( F_{0; X_1}\bigl(
F_{0; X_1}^{-1}(\tilde p_1)-\eta\bigr),
p_2 \bigr) \,d \eta\nonumber
\\
&&\qquad = \int_{-\infty}^{\infty}
f_{\eta}(\eta) \frac{ f_{0; X_1}(
F_{0; X_1}^{-1}(\tilde p_1)-\eta) }{ f_{0; X_1}( F_{0;
X_1}^{-1}(\tilde p_1) )}
\nonumber\\[-8pt]\label{A14}\\[-8pt]
&&\hspace*{53pt}{}\times f_{0; (p_1, p_2)} \bigl(
F_{0; X_1}\bigl( F_{0; X_1}^{-1}(\tilde p_1)-\eta\bigr), p_2 \bigr) \,d \eta, \nonumber
\end{eqnarray}
where (\ref{A11}) and (\ref{A12}) are due to the fact that $f_{0;
X_1}$ is symmetric with respect to 0,
(\ref{A13}) is satisfied by using the symmetry property assumption on
$f_{0; (p_1, p_2)}(p_1, p_2)$, and (\ref{A14}) holds under the
assumption that the p.d.f. of $\eta$
is symmetric. (\ref{A10}) together with (\ref{A14}) yields that
\[
f_{0; (\tilde p_1, p_2)}(\tilde p_1, p_2)=f_{0; (\tilde p_1, p_2)}(1-
\tilde p_1, 1-p_2),
\]
for any $\tilde p_1$ and $p_2$ in $[0,1] \times[0,1]$. The case
for right-sided hypotheses can be derived in a similar way. These
complete the proof.
\end{pf*}

\begin{pf*}{Proof of Theorem \ref{Theorem-1}}
Before proving Theorem~\ref{Theorem-1}, we first provide Lemma~\ref{Lemma-3}.
%
\begin{lemma} \label{Lemma-3}
Under conditions \textup{(C1}) to \textup{(C8)},
\[
\sup_{\theta} \bigl|\hat{t}{}_{\alpha'}^{* \II } (
\theta )-t^*_{\alpha'}(\theta) \bigr| \stackrel{\mathit{a.s.}} \to 0.
\]
\end{lemma}

\begin{pf}
Fix $\delta_1>0$, and let $\bar {t}(\theta)$ be any curve such
that $t^*_{\alpha'}(\theta)+\delta_1 \le\bar {t}(\theta)\le1$.
Then
\begin{eqnarray*}
&& \frac{ \widehat {F }_0^{ { \II} }( \bar {t}(\theta),
\theta) }{
 \{R( \bar {t} (\theta), \theta)\vee1\}/m }
\\
&&\qquad \ge \frac{ F_0( \bar {t}(\theta), \theta)-
| \widehat F_0^{ { \II} }( \bar {t}(\theta), \theta)-F_0(
\bar {t}(\theta), \theta) | }{
 F( \bar {t}(\theta), \theta) + | \{ R( \bar {t}(\theta
), \theta)\vee1\}/m-F( \bar {t}(\theta), \theta) | }
\\
&&\qquad \ge \frac{ \inf_{\theta} F_0( \bar {t}(\theta), \theta
)/F( \bar {t}(\theta), \theta)
-\epsilon_1 }{1+ \epsilon_2},
\end{eqnarray*}
where $\epsilon_1= \inf_{\theta} \inf_{ t \ge\delta_1 }
|\widehat F_0^{ { \II} }(t, \theta)-F_0(t, \theta) |/ F( t, \theta
) $,
and
\[
\epsilon_2= \sup_{\theta} \sup_{ t \ge\delta_1}
\bigl| \bigl\{R( t, \theta)\vee1\bigr\}/m-F(t, \theta) \bigr|/ F( t, \theta).
\]
By Lemmas~\ref{Lemma-1}, \ref{Lemma-2} and condition \textup{(C6)},
$\epsilon_1 \stackrel{\mathrm{a.s.}} \to 0$ and $\epsilon_2
\stackrel{\mathrm{a.s.}} \to 0$.
Note that $F_0( \bar {t}(\theta), \theta)/\break  F( \bar
{t}(\theta), \theta)>\alpha'$;
otherwise it contradicts $t^*_{\alpha'}(\theta)$ being supremum. By
condition~\textup{(C7)},
$\inf_{\theta}F_0( \bar {t}(\theta), \theta)/F( \bar
{t}(\theta), \theta)>\alpha' $.
Hence, for a sufficiently large $M_1(\delta_1)$, when $m>M_1(\delta
_1)$, if follows that
\[
m \widehat F_0^{ { \II} }\bigl( \bar {t}(\theta), \theta
\bigr)/ \bigl\{ R\bigl( \bar {t}(\theta), \theta\bigr)\vee1\bigr\} >
\alpha'
\]
with probability $1$, which implies that
$ \hat{t}{}_{\alpha'}^{ *\II}(\theta) \le t^*_{\alpha'}(\theta
)+\delta_1$ almost surely.\vspace*{1pt}

On the other hand, by condition \textup{(C8)}, since ${F_0(t, \theta
)} / {F(t, \theta) }$ has a nonzero derivative
$k(\theta)$ at $t^*_{\alpha'}(\theta)$, it must be positive;
otherwise $t^*_{\alpha'}(\theta)$
cannot be the true supremum for all $t$ such that ${F_0(t, \theta)} /
{F(t, \theta) } \le\alpha'$.
For any $\varepsilon>0$, there exists $\xi>\delta_1$ such that, for
$| \tilde t(\theta)-t^*_{\alpha'}(\theta) |\le\xi$,
\[
\biggl| \frac{ F_0( \tilde t(\theta), \theta)/F(\tilde t(\theta), \theta)-
F_0( t^*_{\alpha'}(\theta), \theta)/F( t^*_{\alpha'}(\theta),
\theta) }{
 \tilde t(\theta)-t^*_{\alpha'}(\theta) }-k(\theta) \biggr|<\varepsilon.
\]
For a truncated area with $t^*_{\alpha'}(\theta)-\xi\le\tilde t(\theta) \le t^*_{\alpha'}(\theta)- \delta_1$,
$\sup_{\theta} {F_0( \tilde t(\theta), \theta) } /\break { F(
\tilde t(\theta), \theta) }< \alpha'$.
When $\tilde t(\theta) \in[t^*_{\alpha'}(\theta)-\xi,
t^*_{\alpha'}(\theta)-\delta_1]$, some derivation yields that
\begin{eqnarray*}
&& \frac{ \widehat F_0^{ { \II} } ( \tilde t(\theta),\theta)
} { \{R( \tilde t(\theta), \theta)\vee1 \}/m }
\\
&&\qquad \le \frac{ F_0( \tilde t(\theta), \theta)/F(\tilde t(\theta), \theta) +
|F_0( \tilde t(\theta), \theta) -\widehat F_0^ {\II}(
\tilde t(\theta), \theta)|/F( \tilde t(\theta), \theta) }{
 1- |F( \tilde t(\theta), \theta)-\{R( \tilde t(\theta),
\theta)\vee1 \}/m|/F( \tilde t(\theta), \theta) }
\\
&&\qquad \le \frac{ \sup_{\theta} F_0( \tilde t(\theta), \theta)/
F(\tilde t(\theta), \theta) +\epsilon_3 }{1-\epsilon_4},
\end{eqnarray*}
where\vspace*{1pt} $\epsilon_3= \sup_{\theta} \sup_{t \ge\delta^{+} } |F_0(t,
\theta) -\widehat F_0^ {\II}(t, \theta) |/F(t, \theta)$,
$\epsilon_4= \inf_{\theta} \inf_{t \ge\delta^{+} } |F(t, \theta
)-\{R( t, \theta)\vee1 \}/m |/F( t, \theta) $ and
$\delta^{+}=\inf_{\theta} \{ t^*_{\alpha'}(\theta)-\xi\}$.
By\vspace*{1pt} Lemmas~\ref{Lemma-1}~and~\ref{Lemma-2}, and condition \textup{(C6)}, it follows that
$\epsilon_3 \stackrel{\mathrm{a.s.}} \to 0$ and $\epsilon_4
\stackrel{\mathrm{a.s.}} \to 0$. Thus, for another sufficiently
large $M_2(\delta_1)$, when $m>M_2(\delta_1)$,
\[
{ m \widehat F_0^{ { \II} } \bigl( \tilde t(\theta),\theta
\bigr) } / { \bigl\{ R\bigl( \tilde t(\theta), \theta\bigr)\vee1 \bigr\} }<
\alpha'
\]
with probability $1$, which implies that $\hat{t}{}_{\alpha'}^{ *
\II}(\theta) \ge t^*_{\alpha'}(\theta)-\delta_1 $ almost surely.
Combining this and previous result, we obtain that \mbox{$ \sup_{\theta}
| \hat{t}{}_{\alpha'}^{ *\II} (\theta) - t^*_{\alpha'}(\theta) |
\stackrel{\mathrm{a.s.}}{\longrightarrow} 0 $}.
\end{pf}

Now, we prove Theorem \ref{Theorem-1}. First, we show the uniform
consistency of $\widehat F( \hat{t}{}_{\alpha'}^{ *\II} (\theta
), \theta)$, that is,
%
\begin{equation}
\label{A15} \sup_{\theta} \bigl| \widehat F \bigl( \hat{t}{}_{\alpha'}^{ *\II} (\theta), \theta\bigr) - F\bigl(
t^*_{\alpha'}(\theta), \theta\bigr) \bigr| \stackrel{\mathrm{a.s.}} \to 0.
\end{equation}
The left-hand side of (\ref{A15}) can be decomposed as
\begin{eqnarray*}
&& \bigl|\widehat F \bigl( \hat{t}{}_{\alpha'}^{* \II} (\theta), \theta
\bigr) -F\bigl( t^*_{\alpha'}(\theta), \theta\bigr) \bigr|
\\
&&\qquad \le \bigl| \widehat F
\bigl( \hat{t}{}_{\alpha'}^{* \II} (\theta), \theta\bigr)- F\bigl(
\hat{t}{}_{\alpha'}^{* \II} (\theta), \theta\bigr)\bigr| +\bigl|F\bigl(
\hat{t}{}_{\alpha'}^{* \II} (\theta), \theta\bigr)-F\bigl(
t^*_{\alpha'}(\theta), \theta\bigr) \bigr|
\\
&&\qquad \le \sup_{\theta}
\sup_t \bigl| \widehat F(t, \theta) -F(t, \theta )\bigr|+
C_2 \bigl| \hat{t}{}_{\alpha'}^{* \II} (
\theta)-t^*_{\alpha
'}(\theta)\bigr|.
\end{eqnarray*}
Thus, (\ref{A15}) is obtained by condition \textup{(C7)}, Lemmas~\ref
{Lemma-1}~and~\ref{Lemma-3} directly.

For presentational fluency, denote
$\hat{\theta}{}^{\II}(\alpha') $ by $\hat{\theta}_m(\alpha
') $. For each subsequence
$\{ \hat{\theta}_{m_k}(\alpha')\dvtx  k=1,\ldots\}$,
there exists a subsequence $\{ \hat{\theta}_{m_{k,l} }(\alpha
')\dvtx  l=1,\ldots\}$ such that $ \lim_{l \to\infty}
\hat{\theta}_{m_{k,l} }(\alpha') =\theta_{+}(\alpha') $
almost surely.
The next step is to show
%
\begin{equation}
\label{A16} F\bigl( t^*_{\alpha'} \bigl(\theta_{+}\bigl(
\alpha'\bigr) \bigr), \theta_{+}\bigl(\alpha'
\bigr) \bigr)\ge F\bigl( t^*_{\alpha'}\bigl(\theta_0\bigl(
\alpha'\bigr) \bigr), \theta_0\bigl(\alpha'
\bigr) \bigr).
\end{equation}
Thus, $\theta_{+}(\alpha' )=\theta_0(\alpha' )$ by condition
\textup{(C9)}. This completes
the proof.

If\vspace*{2pt} (\ref{A16}) is violated, we have
$F( t^*_{\alpha'}(\theta_{+}(\alpha') ), \theta_{+}(\alpha') ) <
F( t^*_{\alpha'}(\theta_0(\alpha') ), \theta_0(\alpha') )$.
To get contradiction, we partition
$
\widehat F( \hat{t}{}_{\alpha'}^{* \II} ( \hat {\theta}_{m_{k,l}}(\alpha') ), \hat {\theta}_{m_{k,l}}(\alpha') )-
\widehat F ( \hat{t}{}_{\alpha'}^{* \II} ( \theta_0(\alpha')
),\break  \theta_0(\alpha') )
$ as $ A_1+A_2+A_3+A_4$, where
\begin{eqnarray*}
A_1&=& \widehat F \bigl( \hat{t}{}_{\alpha'}^{*\II}
\bigl( \hat {\theta }_{m_{k,l}}\bigl(\alpha'\bigr) \bigr),
\hat {\theta}_{m_{k,l}}\bigl(\alpha'\bigr) \bigr) -F
\bigl( t^*_{\alpha'}\bigl( \hat {\theta}_{m_{k,l}}\bigl(
\alpha'\bigr) \bigr), \hat {\theta}_{m_{k,l}}\bigl(
\alpha'\bigr) \bigr),
\\
A_2&=& F \bigl(
t^*_{\alpha'}\bigl( \hat {\theta}_{m_{k,l}}\bigl(
\alpha'\bigr) \bigr), \hat {\theta}_{m_{k,l}}\bigl(
\alpha'\bigr) \bigr) - F \bigl( t^*_{\alpha'}\bigl(
\theta_{+}\bigl(\alpha'\bigr) \bigr),
\theta_{+}\bigl(\alpha'\bigr) \bigr),
\\
A_3&=& F \bigl( t^*_{\alpha'}\bigl( \theta_{+}
\bigl(\alpha'\bigr) \bigr), \theta_{+}\bigl(
\alpha'\bigr) \bigr) -F \bigl( t^*_{\alpha'}\bigl(
\theta_0\bigl(\alpha'\bigr) \bigr),
\theta_0\bigl(\alpha'\bigr) \bigr),
\\
A_4&=& F \bigl( t^*_{\alpha'}\bigl( \theta_0\bigl(
\alpha'\bigr) \bigr), \theta_0\bigl(\alpha'
\bigr) \bigr) -\widehat F \bigl( \hat{t}{}_{\alpha'}^{ *\II}
\bigl( \theta_0\bigl(\alpha'\bigr) \bigr),
\theta_0\bigl(\alpha'\bigr) \bigr).
\end{eqnarray*}
The term $A_1$ can be bounded by $ \sup_{\theta} |\widehat F (
\hat{t}{}_{\alpha'}^{*\II} ( \theta), \theta)
-F ( t^*_{\alpha'}( \theta), \theta) |$, which is $o(1)$ by (\ref
{A15}). Similarly,\vspace*{1pt} $A_4=o(1)$.
By continuous mapping theorem, $A_2$ is $o(1)$. Thus,
$\widehat F ( \hat{t}{}_{\alpha'}^{*\II} ( \hat {\theta
}_{m_{k,l}}(\alpha') ), \hat {\theta}_{m_{k,l}}(\alpha') ) -
\widehat F ( \hat{t}{}_{\alpha'}^{*\II} ( \theta_0(\alpha') ),
\theta_0(\alpha') ) <0$ almost\vspace*{1pt} surely,
which contradicts the fact that
$\widehat F ( \hat{t}{}_{\alpha'}^{*\II} ( \hat {\theta
}_{m_{k,l}}(\alpha') ), \hat {\theta}_{m_{k,l}}(\alpha') ) \ge
\widehat F ( \hat{t}{}_{\alpha'}^{*\II} ( \theta_0(\alpha') ),
\theta_0(\alpha') )$ obtained from (\ref{47}).
\end{pf*}
\begin{pf*}{Proof of Theorem \ref{Theorem-2}}
To justify Theorem~\ref{Theorem-2}, we first provide Lemmas~\ref{Lemma-4}~and~\ref{Lemma-5} below.
%
\begin{lemma} \label{Lemma-4}
Let $ \hat {\pi}_{0\#}^{ \II} (\lambda, \theta)=\frac{ \sum_{i} \ID\{ p_i(\theta) >\lambda, H_0(i) \} }{
 m\{ 1-\widehat F_0^{ \II} (\lambda, \theta)\}}
=\frac{ m_0-V(\lambda, \theta) }{ m\{ 1-\widehat F_0^{ \II}
(\lambda, \theta)\}}$,
where $ 0< \lambda \le1/2$.
Then under conditions \textup{(C1}) to \textup{(C5)},
\begin{eqnarray*}
&&\lim_{m \to \infty} \sup_{\theta} \sup
_{ 0 < \lambda\le1/2} \bigl| \hat {\pi}_{0 \#}^{ \II} (
\lambda, \theta)-\pi_0 \bigr| \stackrel{\mathit{a.s.}} \to 0.
\end{eqnarray*}
\end{lemma}

\begin{pf}
By decomposing,
\begin{eqnarray*}
\bigl| \hat {\pi}_{0 \#}^{\II} (\lambda, \theta)-
\pi_0 \bigr| & \le& \biggl| \frac{ m_0}{m}-\pi_0 \biggr| \biggl|
\frac{ 1-V(\lambda, \theta)/m_0}{ 1- \widehat F_0^{\II} (\lambda, \theta) } \biggr| + \pi_0 \biggl| \frac{ V(\lambda, \theta)/m_0- \widehat F_0^{ \II
}(\lambda, \theta) }{ 1- \widehat F_0^{\II} (\lambda, \theta) } \biggr|
\\
&=&
\Pi_1(\lambda, \theta)+\Pi_2(\lambda, \theta).
\end{eqnarray*}
Uses of
\[
\sup_{0 < \lambda\le1/2} \sup_{\theta} \bigl|1- V(\lambda, \theta)
/m_0\bigr|\le2\quad\mbox{and}\quad
\inf_{0 < \lambda\le1/2} \inf_{\theta} \bigl|1-\widehat F_0^{\II
}(\lambda, \theta) \bigr|\ge1/2
\]
yield that
\[
\lim_{m \to \infty} \sup_{\theta} \sup
_{0 < \lambda\le1/2} \Pi_1(\lambda, \theta) \le\lim
_{m \to\infty} 4 \biggl| \frac{m_0}{m}-1 \biggr| =0\qquad\mbox{almost surely}.
\]
For the term $\Pi_2(\lambda, \theta)$, it suffices to show that
\begin{eqnarray*}
&&\lim_{m \to\infty} \sup_{ 0 < \lambda\le1/2} \sup
_{\theta} \bigl| V(\lambda, \theta)/m_0- \widehat
F_0^{ \II}(\lambda, \theta) \bigr| \stackrel{\mathrm{a.s.}} \to
0,
\end{eqnarray*}
which is completed by using Lemmas~\ref{Lemma-1}~and~\ref{Lemma-2}.
\end{pf}

%
\begin{lemma} \label{Lemma-5}
Suppose conditions \textup{(C1}) to \textup{(C6)} hold. Then, for
each $\delta>0$,
%
\begin{eqnarray}
\lim_{m \to \infty} \inf_{t \ge\delta} \inf
_{\theta}\bigl\{ \widehat { \FDR}_{\lambda}^{ \II}
(t, \theta)-\FDR(t, \theta)\bigr\} &\ge& 0, \label{A17}
\\
\lim_{m \to \infty} \inf_{t \ge\delta} \inf
_{\theta} \biggl\{ \widehat { \FDR}{}^{ \II}_{\lambda}(t,
\theta)- \frac{ V(t, \theta
) }{R(t, \theta)\vee1} \biggr\} &\ge& 0 \label{A18}
\end{eqnarray}
with probability $1$, where $\widehat { \FDR}{}^{\II}_{\lambda}(t,
\theta)=
\frac{ \hat {\pi}{}^{\II}_0(\lambda, \theta) \widehat F_0^{\II
}(t, \theta) }{ \{R(t, \theta)\vee1\}/m }$,
for fixed $\lambda$.
Furthermore, the estimator $\widehat { \FDR}{}^{ \II}_{ \lambda^{*}
}(t, \theta)$ with $\lambda^*$ arbitrarily selected from
the sequence of values $\{ \lambda_j\dvtx  j=1,\ldots, n\}$ of a finite
size is simultaneously conservatively consistent for $\FDR(t, \theta
)$ or $\frac{V(t, \theta)}{ R(t, \theta) \vee1}$ for all $t \ge
\delta$ and $\theta$.
\end{lemma}

\begin{pf}
By Lemma~\ref{Lemma-1}, we have
%
\begin{eqnarray}
\lim_{m \to \infty} \sup_{t} \sup
_{\theta} \biggl| \frac{ V(t,
\theta) }{m}- \pi_0
F_0(t, \theta) \biggr| &\stackrel{\mathrm{a.s.}} \to & 0, \label{A19}
\\
\lim_{m \to \infty} \sup_{t} \sup
_{\theta} \biggl| \frac{ R(t,
\theta)\vee1 }{m}- \bigl\{ \pi_0
F_0(t, \theta)+ \pi_1 F_1(t, \theta)\bigr
\} \biggr| &\stackrel{\mathrm{a.s.}} \to & 0. \label{A20}
\end{eqnarray}
To show (\ref{A18}), we observe that
\begin{eqnarray*}
&& \widehat { \FDR}_{\lambda} (t, \theta)-\frac{ V(t, \theta) }{R(t,
\theta)\vee1 }
\\
&&\qquad =
\frac{ \hat {\pi}_0 ^{\II} (\lambda, \theta) \widehat
F_0^{ \II}(t, \theta)-\pi_0 F_0(t, \theta) } { \{ R(t, \theta)
\vee1\}/m }
-\frac{ V(t, \theta)/m-\pi_0 F_0(t, \theta) }{ \{ R(t, \theta)
\vee1\}/m }
\\
&&\qquad = I_1(t,
\theta)-I_2(t, \theta).
\end{eqnarray*}
For the term $I_2(t, \theta)$, applying (\ref{A19}), (\ref{A20})
and condition \textup{(C6)} yields that
%
\begin{eqnarray}\label{A21}
\qquad & &\lim_{m \to\infty} \sup_{t \ge\delta} \sup
_{\theta} \bigl|I_2(t, \theta)\bigr|\nonumber
\\
&&\qquad \le \lim
_{m \to\infty} \sup_{\theta} \frac{m}{ R(\delta, \theta
)\vee1 } \times
\lim_{m \to\infty} \sup_{t \ge\delta} \sup
_{\theta} \biggl| \frac{ V(t, \theta) }{m}- \pi_0
F_0(t, \theta) \biggr|
\\
&&\qquad \stackrel{\mathrm{a.s.}} \to  0.\nonumber
\end{eqnarray}
For the term $I_1(t, \theta)$, using the fact that $\hat {\pi}{}^{
\II}_0(\lambda, \theta) \ge
\hat {\pi}{}^{ \II }_{ 0 \# }( \lambda, \theta) $, we have
%
\begin{equation}
\label{A22} \qquad\quad \lim_{m \to\infty} \inf_{t \ge\delta} \inf
_{\theta} I_1(t, \theta) \ge \lim
_{m \to\infty} \inf_{t \ge\delta} \inf_{\theta}
\bigl\{ \hat {\pi}{}^{ \II}_{0 \#} (\lambda, \theta) \widehat
F^{\II
}_0(t, \theta)-\pi_0 F_0(t,
\theta)\bigr\}.
\end{equation}
To show that the right-hand side of (\ref{A22}) converges to $0$
almost surely, it suffices to verify
%
\begin{equation}
\label{A23} \lim_{m \to\infty} \sup_{t \ge\delta} \sup
_{\theta} \bigl| \hat {\pi}_{0 \#}^{\II} (
\lambda, \theta) \widehat F^{ \II
}_0(t, \theta)-
\pi_0 F_0(t, \theta) \bigr| \stackrel{\mathrm{a.s.}} \to 0,
\end{equation}
which can be achieved by Lemmas~\ref{Lemma-2}~and~\ref{Lemma-4}.
Combining (\ref{A21}), (\ref{A22}) and (\ref{A23}) completes the
proof of (\ref{A18}).

To show (\ref{A17}), it suffices to show that
%
\begin{equation}
\label{A24} \lim_{m \to\infty} \sup_{t \ge\delta}\sup
_{\theta} \biggl| \frac
{ V(t, \theta) }{ R(t, \theta) \vee1}- \FDR (t, \theta) \biggr| \stackrel{
\mathrm{a.s.}} \to 0.
\end{equation}
Since $\inf_{\theta} F(\delta, \theta) >0 $ and $\{ R(t,\theta)$,
$F(t, \theta) \}_{\theta}$ are nondecreasing functions for $t$,
it is straightforward to show that
\[
\lim_{m \to\infty} \sup_{t \ge\delta} \sup
_{\theta} \biggl| \frac{m}{R(t, \theta) \vee1}-\frac{1}{F(t, \theta)} \biggr| \stackrel{
\mathrm{a.s.}} \to 0.
\]
Using this, inequality (\ref{A21}) and the triangle inequality, we obtain
%
\begin{eqnarray}\label{A25}
&& \lim_{m \to \infty} \sup_{t \ge\delta} \sup
_{ \theta} \biggl| \frac{ V(t, \theta) }{ R(t, \theta) \vee1 }-\frac{ \pi_0
F_0(t, \theta) }{ F(t, \theta) } \biggr|\nonumber
\\
&&\qquad \le
\lim_{m \to \infty} \sup_{t \ge\delta} \sup
_{ \theta} \biggl| \frac{ V(t, \theta) }{ R(t, \theta) \vee1 }-\frac{ m \pi_0
F_0(t, \theta) }{ R(t, \theta) \vee1 } \biggr|
\nonumber\\[-8pt]\\[-8pt]
&&\quad\qquad{} +
\lim_{m \to \infty} \sup_{t \ge\delta} \sup
_{ \theta} \biggl| \frac{ m \pi_0 F_0(t, \theta) }{ R(t, \theta) \vee1 }- \frac{ \pi_0 F_0(t, \theta) }{ F(t, \theta)} \biggr|\nonumber
\\
&&\qquad \stackrel{\mathrm{a.s.}} \to  0. \nonumber
\end{eqnarray}
By (\ref{A25}), (\ref{A24}) is implied if we can show that
%
\begin{equation}
\label{A26} \lim_{m \to \infty} \sup_{t \ge\delta} \sup
_{ \theta} \biggl| \FDR(t, \theta)- \frac{\pi_0 F_0(t, \theta) }{ F(t, \theta) } \biggr| \stackrel{
\mathrm{a.s.}} \to 0.
\end{equation}
Combining (\ref{A25}) and the fact that $
 | { V(t, \theta) }/{ \{R(t, \theta) \vee1 \} } - {\pi_0 F_0(t,
\theta) }/\break { F(t, \theta) }  | \le2 $, we have
\begin{eqnarray*}
&&\lim_{m \to \infty} \sup_{t \ge\delta} \sup
_{ \theta} \biggl| \FDR(t, \theta)- \frac{\pi_0 F_0(t, \theta) }{ F(t, \theta) } \biggr|
\\
&&\qquad \le \lim
_{m \to\infty} E \biggl\{ \sup_{t \ge\delta} \sup
_{\theta} \biggl| \frac{V(t, \theta)}{ R(t, \theta) \vee1 } -\frac{\pi_0 F_0(t, \theta) }{ F(t, \theta) } \biggr| \biggr\}
\\
&&\qquad \le E \biggl\{ \lim_{m \to\infty} \sup_{t \ge\delta}
\sup_{\theta} \biggl| \frac{V(t, \theta)}{ R(t, \theta) \vee1 } -\frac{\pi_0 F_0(t, \theta) }{ F(t, \theta) } \biggr| \biggr
\}
\\
&&\qquad = 0.
\end{eqnarray*}
This completes the proof of (\ref{A26}).\vadjust{\goodbreak}

Now we turn to show the second part of the lemma.
Let $ \widehat { \FDR}_{*}^{\II} (t, \theta)={ \hat {\pi
}_{0*}^{\II}(\theta) \widehat F_0^{\II}(t, \theta) } /
{( \{R(t, \theta)\vee1 \}/m ) }$,
where $\hat {\pi}_{0*} ^{\II} (\theta)=\min_{j} \hat {\pi
}_0 ^{\II} (\lambda_j, \theta) $.
By Lem\-ma~\ref{Lemma-4} and a slight
modification of the proof in first part, the simultaneously
conservative control of
$ \widehat { \FDR}_*^{ \II} (t, \theta)$ is also satisfied, that is,
%
\begin{eqnarray}
\lim_{m \to \infty} \inf_{t \ge\delta} \inf
_{\theta} \bigl\{ \widehat { \FDR}_*^{ \II} (t, \theta)-
\FDR(t, \theta)\bigr\} &\ge& 0, \label{A27}
\\
\lim_{m \to \infty} \inf_{t \ge\delta} \inf
_{\theta} \biggl\{ \widehat { \FDR}_*^{ \II} (t, \theta)-
\frac{ V(t, \theta
) }{R(t, \theta)\vee1} \biggr\} &\ge& 0. \label{A28}
\end{eqnarray}
The conclusion for $\widehat {\FDR}_{\lambda^{*} }^{\II}(t, \theta
)$ is implied by (\ref{A27}) and (\ref{A28}).
\end{pf}

Now, we show Theorem~\ref{Theorem-2}.
The proof of this theorem is implied by the following inequalities:
%
\begin{eqnarray}
\lim_{m \to \infty} \inf_{t \ge\delta} \bigl\{ \widehat {
\FDR}_{*
}^{\II} \bigl(t, \hat {\theta}{}^{\II}
\bigl(\alpha'\bigr) \bigr) -\FDR\bigl(t, \theta_0\bigl(
\alpha'\bigr) \bigr)\bigr\} &\ge& 0, \label{A29}
\\
\lim_{m \to \infty} \inf_{t \ge\delta} \biggl\{ \widehat {
\FDR }_{ * }^{\II} \bigl(t, \hat {\theta}{}^{\II}
\bigl(\alpha'\bigr) \bigr) -\frac{ V(t, \theta_0(\alpha') ) }{ R( t, \theta_0(\alpha') )\vee
1} \biggr\} &\ge& 0
\label{A30}
\end{eqnarray}
with probability $1$.

To verify (\ref{A29}), it suffices to show that
%
\begin{eqnarray}
\lim_{m \to \infty} \inf_{t \ge\delta} \bigl\{ \widehat {
\FDR }_{*}^{\II} \bigl(t, \hat {\theta}{}^{ \II}
\bigl(\alpha'\bigr) \bigr) - \FDR \bigl(t, \widehat {
\theta}{}^{ \mathrm{II} }\bigl(\alpha'\bigr) \bigr)\bigr\} &\ge& 0,
\label{A31}
\\
\lim_{m \to \infty} \sup_{t \ge\delta} \bigl| \FDR \bigl( t,
\hat {\theta}{}^{ \II}\bigl(\alpha'\bigr) \bigr) - \FDR
\bigl( t, \theta_0\bigl(\alpha'\bigr) \bigr) \bigr| &
\stackrel{\mathrm{a.s.}} \to & 0. \label{A32}
\end{eqnarray}
Note that (\ref{A31}) is readily implied by (\ref{A27}). By using
(\ref{A26}), (\ref{A32}) is implied by
%
\begin{equation}
\label{A33} \lim_{m \to\infty} \sup_{t \ge\delta} \biggl|
\frac{ \pi_0 F_0(t, \hat{\theta}{}^{ \mathrm{II} }(\alpha') ) }{
 F(t, \hat{\theta}{}^{ \mathrm{II} }(\alpha') ) } - \frac{ \pi_0 F_0(t, \theta_0(\alpha') ) } { F(t, \theta_0(\alpha
') ) } \biggr| \stackrel{\mathrm{a.s.}} \to 0.
\end{equation}
The proof of (\ref{A33}) is completed by using condition \textup{(C10)} and Theorem~\ref{Theorem-1}.
For~(\ref{A30}), directly applying (\ref{A24}) and (\ref{A29})
completes the proof.
\end{pf*}

\begin{pf*}{Proof of Theorem \ref{Theorem-3}}
First, we will show the uniform consistency of
$\widehat { \FDR} _{\lambda}^{ \II} (t, \theta) $ for fixed
$\lambda$, that is,
\[
\lim_{m \to\infty} \sup_{ t \ge\delta} \sup
_{\theta} \bigl| \widehat { \FDR} _{\lambda}^{ \II} (t,
\theta)- \widehat { \FDR}_{\lambda}^{\infty}(t, \theta) \bigr| \stackrel{
\mathrm{a.s.}} \to 0\qquad\mbox{for any } \delta>0.
\]
This can be completed by a slight modification of Lemma~\ref{Lemma-5}.
Following the similar arguments of (\ref{A29}) and (\ref{A30}), we obtain
%
\begin{equation}
\label{A34} \lim_{m \to\infty} \sup_{ t \ge\delta} \bigl|
\widehat { \FDR} _{\lambda} \bigl(t, \hat {\theta}{}^{ \II}\bigl(
\alpha '\bigr) \bigr)- \widehat { \FDR}_{\lambda}^{\infty}
\bigl(t, \theta_0\bigl(\alpha'\bigr) \bigr) \bigr|
\stackrel{\mathrm{a.s.}} \to 0.
\end{equation}

Abbreviate $t_{\alpha}(\widehat { \FDR}_{\lambda}^{ \II} (:,
\hat {\theta}{}^{ \II}(\alpha') ) ) $ by $t_{\alpha}^{\lambda}$.
According to the condition, for each $\lambda_j$, there is $t_j$ such that
$\alpha-\widehat { \FDR}_{\lambda_j}^{\infty} ( t_j, \theta
_0(\alpha') ) = \varepsilon_j >0$. By (\ref{A34}), we can
take $m$ sufficiently large that
$| \widehat { \FDR}_{\lambda_j}^{\infty} (t_j,\theta_0(\alpha')
)- \widehat { \FDR}_{\lambda_j}^{ \II}
(t_j, \hat {\theta}{}^{ \II}(\alpha') ) |
< \varepsilon_j$, which implies that $\widehat { \FDR}_{\lambda
_j}^{\II}
(t_j, \hat {\theta}{}^{ \II}(\alpha') ) < \alpha$ and $t_{\alpha
}^{\lambda_j} \ge t_j $. Therefore,\break
$\liminf_{m \to\infty} t_{\alpha}^{\lambda_j} \ge t_j $ with
probability $1$.
For $\delta_j=t_j/2$,
\begin{eqnarray*}
&&\liminf_{m \to\infty} \biggl\{ \widehat { \FDR}_{\lambda_j}^{\II
}\bigl( t_{\alpha}^{\lambda_j}, \hat {\theta}{}^{\II}\bigl(
\alpha'\bigr) \bigr) -\frac{ V( t_{\alpha}^{\lambda_j}, \hat {\theta}{}^{\II}(\alpha
') ) }{
 R( t_{\alpha}^{\lambda_j},\hat {\theta}{}^{\II}(\alpha') )
\vee1 } \biggr\}
\\
&&\qquad  \ge \lim
_{ m \to\infty} \inf_{ t\ge\delta_i } \biggl\{ \widehat { \FDR
}_{\lambda_j}^{\II} \bigl( t, \hat {\theta }{}^{ \II}
\bigl(\alpha'\bigr) \bigr) -\frac{ V( t, \hat {\theta}{}^{\II}(\alpha') ) }{
 R( t,\hat {\theta}{}^{\II}(\alpha') ) \vee1 } \biggr\}
\\
&&\qquad  \ge
\lim_{ m \to\infty} \inf_{ t\ge\delta_i }\inf
_{\theta} \biggl\{ \widehat { \FDR }_{\lambda_j}^{\II}
( t, \theta) -\frac{ V( t, \theta) }{ R( t, \theta) \vee1 } \biggr\}
\\
&&\qquad \ge 0,
\end{eqnarray*}
where the last inequality is due to (\ref{A18}).
By the definition of $t_{\alpha}^{\lambda_j}$,
$\widehat { \FDR}_{\lambda_j}^{\II} ( t_{\alpha}^{\lambda_j},\break
\hat {\theta}{}^{ \II }(\alpha') ) \le\alpha$,
and it follows that
\begin{eqnarray*}
&&\limsup_{ m \to\infty} \biggl\{ \frac{ V( t_{\alpha}^{\lambda_j},
\hat {\theta}{}^{\II}(\alpha') ) }{
 R(t_{\alpha}^{\lambda_j}, \hat {\theta}{}^{\II}(\alpha') )
\vee1} \biggr\} \le
\alpha
\end{eqnarray*}
with probability $1$. Let $\lambda^*$ be determined by the
algorithm in Section~\ref{Sec-43}. Then
\begin{eqnarray*}
&& \limsup_{ m \to\infty} \biggl\{ \frac{ V( t_{\alpha}^{\lambda
^{*}}, \hat {\theta}{}^{\II}(\alpha') ) }{
 R(t_{\alpha}^{\lambda^{*}}, \hat {\theta}{}^{\II}(\alpha') )
\vee1 } \biggr\}
\\
&&\qquad \le
\limsup_{ m \to\infty} \biggl\{ \max_{ 1 \le j \le n}
\frac{
V( t_{\alpha}^{ \lambda_j }, \hat {\theta}{}^{\II}(\alpha') ) }{
 R(t_{\alpha}^{ \lambda_j }, \hat {\theta}{}^{\II}(\alpha') )
\vee1 } \biggr\}
\le\alpha
\end{eqnarray*}
with probability $1$. Following Fatou's lemma,
\begin{eqnarray*}
&& \limsup_{ m \to\infty} E \biggl\{ \frac{ V( t_{\alpha}^{\lambda
^{*}}, \hat {\theta}{}^{\II}(\alpha') ) }{
 R(t_{\alpha}^{\lambda^{*}}, \hat {\theta}{}^{\II}(\alpha') )
\vee1 } \biggr\}
\\
&&\qquad \le E
\biggl[ \limsup_{ m \to\infty} \biggl\{ \frac{ V( t_{\alpha
}^{\lambda^{*}}, \hat {\theta}{}^{\II}(\alpha') ) }{
 R(t_{\alpha}^{\lambda^{*}}, \hat {\theta}{}^{\II}(\alpha') )
\vee1 } \biggr\}
\biggr]
\le \alpha.
\end{eqnarray*}\upqed
\end{pf*}
%

\section{Density of the bivariate $p$-value when the bivariate test
statistic under the true null is a bivariate normal or $t$ distribution}\label{appB}
Assume that we are interested in testing the left-sided hypotheses,
%
\begin{equation}
H_0\dvtx  \mu=\mu_0\quad\mbox{versus}\quad
H_1\dvtx  \mu< \mu_0, \label{B1}
\end{equation}
where $\mu$ is the parameter involved in some population and $\mu_0$
is given.
The right-sided hypotheses can be discussed similarly.
Suppose that ${\boldsymbolX}=(X_1,X_2)$ are the preliminary and
primary test statistics with the true null joint $\CDF$ $F_{0; (X_1,X_2)}$.
Denote by $F_{0; X_1}$ and $F_{0; X_2}$ the marginal $\CDF$s of $X_1$
and $X_2$ under the true null, respectively.
The joint $\CDF$ of $\mathbf{p}=(p_1,p_2)$
under the true null hypothesis of (\ref{B1}) has the following form:
\begin{eqnarray*}
F_0 ( \boldsymbolt) &=& \P( \mathbf{p}\le\boldsymbolt \mid
H_0 )
\\
&=& \P\bigl( F_{0; X_1}(X_1) \le
t_1, F_{0; X_2}(X_2) \le t_2 \mid
H_0 \bigr)
\\
&=& F_{0; (X_1,X_2) }\bigl(F_{0; X_1}^{-1}(t_1),
F_{0; X_2}^{-1}(t_2) \bigr),
\end{eqnarray*}
with\vspace*{1pt} $F_{0; X_1}^{-1}$ and $F_{0; X_2}^{-1}$ being the inverse
functions of $F_{0; X_1}$ and $F_{0; X_2}$, respectively.

If ${\boldsymbolX}=(X_{1},X_{2})$ under the true null follows a
bivariate normal distribution with
mean zero and covariance matrix $\Sigma_0$ given by (\ref{38}) in
Section~\ref{Sec-33},
then direct calculations yield that
%
\begin{eqnarray}\label{B2}
\qquad F_0( \mathbf{p}) &=& \int_{-\infty}^{ \Phi^{-1}(p_1) }
\int_{-\infty}^{ \Phi
^{-1}(p_2) } \frac{1}{2\pi\sqrt{1-\rho_0^2}} \exp \biggl\{-
\frac{x^2-2\rho_0 xy +y^2}{2(1-\rho_0^2)} \biggr\}\,dx \,dy,\nonumber
\\
f_0( \mathbf{p}) &=& \frac{1}{\sqrt{1-\rho_0^2}}
\exp \bigl(- \bigl( \rho_0^2\bigl\{\Phi^{-1}(p_1)\bigr\}^2
\nonumber\\[-8pt]\\[-8pt]
&&\hspace*{72pt}{} -2\rho_0 \Phi^{-1}(p_1)\Phi
^{-1}(p_2)+\rho_0^2\bigl\{\Phi^{-1}(p_2)\bigr\}^2\bigr)\nonumber
\\
&&\hspace*{206pt}{} /\bigl({ 2\bigl(1-\rho_0^2\bigr) }\bigr) \bigr), \nonumber
\end{eqnarray}
where $\Phi$ is the standard normal $\CDF$.

If $\boldsymbolX$ under the true null has a bivariate $t$ distribution
with $v$ degrees of freedom and correlation coefficient $\rho_0$,
then derivations similar to (\ref{B2}) imply that
%
\begin{eqnarray}\label{B3}
F_0( \mathbf{p}) &=& \int_{-\infty}^{T_v^{-1}(p_1)}
\int_{-\infty}^{T_v^{-1}(p_2)} \frac{1}{2\pi\sqrt{1-\rho_0^2}} \biggl\{ 1 +
\frac{x^2-2\rho_0 xy +y^2}{v(1-\rho_0^2)} \biggr\}^{-(v+2)/2 } \,dx \,dy,\nonumber\hspace*{-5pt}
\\
f_0(
\mathbf{p}) &=& \frac{\{\Gamma({v}/2)\}^2v}{2\{\Gamma((v+1)/2)\}^2 \sqrt
{1-\rho_0^2}}
\nonumber\hspace*{-5pt}\\[-4pt]\\[-12pt]
&&{}\times  \biggl( \biggl[ 1+ \frac{\{T_v^{-1}(p_1)\}^2-2\rho_0
T_v^{-1}(p_1)T_v^{-1}(p_2) + \{T_v^{-1}(p_2)\}^2}{v(1-\rho_0^2)}
\biggr]^{-(v+2)/2 } \biggr)\hspace*{-5pt}\nonumber\hspace*{-5pt}
\\
&&\hspace*{14pt}{}\Big/
\biggl(\biggl[1+\frac{\{T_v^{-1}(p_1)\}^2}{v}\biggr]^{-(v+1)/2 } \biggl[1+\frac{\{
T_v^{-1}(p_2)\}^2}{v}\biggr]^{-(v+1)/2}\biggr), \nonumber\hspace*{-5pt}
\end{eqnarray}
where $T_v(x)$ is the $\CDF$ of $t$ distribution with $v$ degrees of
freedom.\eject

\subsection*{Derivation of \texorpdfstring{$\Delta(\theta)$}{Delta(theta)} in Section~\protect\ref{Sec-46}}
%
\begin{eqnarray} \label{B4}
\Delta(\theta) &=& \frac{\partial}{\partial\theta} \bigl\{ F_{0}
\bigl(t_{\alpha}(\theta), \theta\bigr) \bigr\}\bigg|_{ \theta
={\pi/2} } \times
\frac{(\theta-\pi/2) }{F_0(t_{\alpha}(\pi/2), \pi/2 )}\nonumber
\\[3pt]
&=& \biggl[ \biggl\{ \frac{\partial F_0 ( t, \theta ) }{\partial t}
\frac{\partial t_{\alpha}(\theta)} {\partial\theta}+ \frac{\partial F_0 ( t, \theta ) }{\partial\theta} \biggr\} \bigg|_{t=t_{\alpha}(\theta)} \biggr]
\bigg|_{\theta=\pi/2 }
\\[3pt]
&&{} \times\frac{(\theta-\pi/2) }{F_0(t_{\alpha}(\pi/2), \pi/2 )}.\nonumber
\end{eqnarray}
Derivation similar to (\ref{A6}) yields that
%
\begin{equation}
\label{B5} \frac{ \partial t_{\alpha}(\theta) }{ \partial\theta}= \biggl\{ \frac{ \beta' (({\partial F_0( t, \theta ) })/({\partial\theta}))-
(({\partial F_1( t, \theta ) })/({ \partial\theta})) }{
(({\partial F_1 ( t, \theta ) })/({\partial t}))- \beta'
(({\partial F_0( t, \theta ) })/({ \partial t})) } \biggr\}
\bigg|_{t=t_{\alpha}(\theta)}.
\end{equation}
Plugging (\ref{B5}) into (\ref{B4}), $\Delta(\theta)$ can be
expressed explicitly as
%
\begin{eqnarray}\label{B6}
\qquad \Delta(\theta) &=& \biggl({  \biggl[  \biggl\{
\frac{\partial F_0 ( t, \theta ) }{\partial\theta}
\frac{\partial F_1 ( t, \theta ) }{\partial t}-
\frac{\partial F_0 ( t, \theta ) }{\partial t}
\frac{\partial F_1 ( t, \theta ) }{\partial\theta}
 \biggr\}  \bigg|_{t=t_{\alpha}(\theta)}  \biggr]  \bigg|_{\theta=\pi/2 } }\biggr)\nonumber
\\[3pt]
&& {} \Big/\biggl({
\biggl[
 \biggl\{ \frac{\partial F_1 ( t, \theta ) }{\partial t}-\beta'\frac
{\partial F_0 ( t, \theta ) }{\partial t}
 \biggr\}
 \bigg|_{t=t_{\alpha}(\theta)}
 \biggr]  \bigg|_{\theta=\pi/2}
}\biggr)
\\[3pt]
&&{}  \times\frac{(\theta-\pi/2) }{F_0(t_{\alpha}(\pi/2), \pi/2 )}.\nonumber
\end{eqnarray}
Now consider
%
\begin{eqnarray}\label{B7}
\hspace*{-5pt}&& \biggl[ \biggl\{ \frac{ \partial F_1(t, \theta) }{ \partial t} \biggr\} \bigg|_{t=t_{\alpha}(\theta)} \biggr]
\bigg|_{\theta=\pi/2}\nonumber
\\[3pt]
\hspace*{-5pt}&&\qquad = \biggl[ \biggl\{ \frac{\partial}{\partial t} \int
_{0}^{1}\!\int_{0}^{ \Phi ( ({ \Phi^{-1}(t)-\Phi
^{-1}(p_1)\cos(\theta) })/{\sin(\theta)}  ) }
\!f_{1; (p_1, p_2)}(p_1, p_2) \nonumber
\\[3pt]
\hspace*{-5pt}&&\hspace*{239pt}{}\times  d{p_2}
\,d{p_1} \biggr\} \bigg|_{t=t_{\alpha}(\theta)} \biggr] \bigg|_{\theta=\pi/2}
\nonumber\\[-8pt]\hspace*{-5pt}
\\[-6pt]
\hspace*{-5pt}&&\qquad =
\biggl[ \biggl\{ \int_{0}^{1} \frac{ \phi ( ({ \Phi^{-1}(t)-\Phi^{-1}(p_1)\cos(\theta)
})/{\sin(\theta)}  ) }{
 \phi(\Phi^{-1}(t)) \sin(\theta) }
f_{1; (p_1, p_2)}\nonumber
\\[3pt]
\hspace*{-5pt}&&\hspace*{72pt}{}\times \biggl( p_1, \Phi \biggl( \frac{ \Phi^{-1}(t)-\Phi
^{-1}(p_1)\cos(\theta) }{\sin(\theta)}
\biggr) \biggr) \,d{p_1} \biggr\} \bigg|_{t=t_{\alpha}(\theta)} \biggr]
\bigg|_{\theta=\pi/2}\nonumber
\\[3pt]
\hspace*{-5pt}&&\qquad =\int_{0}^{1}
f_{1; (p_1, p_2)}\bigl(p_1, t_{\alpha}(\pi/2)
\bigr)\,d{p_1}, \nonumber
\end{eqnarray}
where $f_{1; (p_1, p_2)}$ is the p.d.f. of $(p_1, p_2)$ under
the nonnull. Analogously,
%
\begin{eqnarray}
\qquad && \biggl\{ \frac{ \partial F_1(t, \theta) }{ \partial\theta} \bigg|_{t=t_{\alpha}(\theta)} \biggr\} \bigg|_{\theta=\pi/2}
\nonumber\\[-8pt]\label{B8} \\[-8pt]
&&\qquad =
\phi\bigl[ \Phi^{-1} \bigl\{ t_{\alpha}(\pi/2) \bigr\} \bigr] \int
_{0}^{1} \Phi ^{-1}(p_1)
f_{1; (p_1, p_2)}\bigl(p_1, t_{\alpha}(\pi/2)
\bigr)\,d{p_1}, \nonumber
\\
&& \biggl\{ \frac{ \partial F_0(t, \theta) }{ \partial t } \bigg|_{t=t_{\alpha}(\theta)} \biggr\} \bigg|_{\theta=\pi/2}
\nonumber\\[-8pt]\label{B9}\\[-8pt]
&&\qquad = \int_{0}^{1} f_{0; (p_1, p_2)}
\bigl(p_1, t_{\alpha}(\pi/2) \bigr)\,d{p_1},\nonumber
\\
&& \biggl\{ \frac{ \partial F_0(t, \theta) }{ \partial\theta} \bigg|_{t=t_{\alpha}(\theta)} \biggr\} \bigg|_{\theta=\pi/2}
\nonumber\\[-8pt] \label{B10}\\[-8pt]
&&\qquad =
\phi\bigl[ \Phi^{-1} \bigl\{ t_{\alpha}(\pi/2) \bigr\} \bigr] \int
_{0}^{1} \Phi ^{-1}(p_1)
f_{0; (p_1, p_2)}\bigl(p_1, t_{\alpha}(\pi/2)
\bigr)\,d{p_1}.\nonumber
\end{eqnarray}
Plugging (\ref{B7}), (\ref{B8}), (\ref{B9}) and (\ref{B10})
into (\ref{B6}), we have
\begin{eqnarray*}
\Delta(\theta) &=& \Biggl( \phi\bigl[ \Phi^{-1}\bigl\{ t_{\alpha}(\pi/2) \bigr\} \bigr]
\\
&&\hspace*{5pt}{}\times \int_{0}^{1}f_{0;
(p_1, p_2)}\bigl(p_1, t_{\alpha}(\pi/2) \bigr)\,d{p_1} \int_{0}^{1}f_{1; (p_1,
p_2)}\bigl(p_1, t_{\alpha}(\pi/2) \bigr)\,d{p_1} \Biggr)
\\
&&{} \Big/
\Biggl(
 \int_{0}^{1}f_{1; (p_1, p_2)}\bigl(p_1, t_{\alpha}(\pi/2) \bigr)\,d{p_1}-\beta
' \int_{0}^{1}f_{0; (p_1, p_2)}\bigl(p_1, t_{\alpha}(\pi/2) \bigr)\,d{p_1} \Biggr)
\\
&&{}  \times \frac{( \theta-\pi/2 ) }{F_0(t_{\alpha}(\pi/2), \pi/2 )}
\\
&&{} \times
\bigl[ E_{H_0} \bigl\{ \Phi^{-1}(p_1) \mid
p_2=t_{\alpha}(\pi/2) \bigr\} -E_{H_1} \bigl\{
\Phi^{-1}(p_1) \mid p_2=t_{\alpha}(\pi/2)
\bigr\} \bigr]
\\
&=& \frac{ \phi[ \Phi^{-1} \{ t_{\alpha}(\pi/2) \} ]
f_{1;p_2}(t_{\alpha}(\pi/2) )f_{0; p_2}(t_{\alpha}(\pi/2) ) }{
 f_{1; p_2}(t_{\alpha}(\pi/2) )-\beta' f_{0; p_2}(t_{\alpha}(\pi
/2) )} \times \frac{(\theta-\pi/2) }{F_0(t_{\alpha}(\pi/2), \pi/2 )}
\\
& &
\times \bigl[ E_{H_0} \bigl\{ \Phi^{-1}(p_1) \mid
p_2=t_{\alpha}(\pi/2) \bigr\} -E_{H_1} \bigl\{
\Phi^{-1}(p_1) \mid p_2=t_{\alpha}(\pi/2)
\bigr\} \bigr],
\end{eqnarray*}
where $f_{0; p_2}$ and $f_{1; p_2}$ are the p.d.f.s of $p_2$
under true null and nonnull, respectively.
\end{appendix}

\section*{Acknowledgments}
The comments of two referees, the Associate Editor and the Co-Editor,
Peter Hall, were greatly appreciated. We thank Daisy Phillips and
Debashis Ghosh for sending the real data.



\printaddresses
\end{document}